\documentclass{iopart}
\usepackage{graphics}
\newtheorem{definition}{Definition}[section]
 \newtheorem{theorem}[definition]{Theorem}
\newtheorem{proposition}[definition]{Proposition}

\newtheorem{corollary}[definition]{Corollary}

\newtheorem{example}[definition]{Example}

\eqnobysec

\begin{document}

\title[An Introduction to  Conformal Ricci Flow]{An Introduction to  Conformal Ricci Flow}

\author{Arthur E.\ Fischer}

\address{Department of Mathematics\\ University of California\\
Santa Cruz, CA 95064, USA}

 \ead{aef@ucsc.edu}

\begin{abstract}We introduce a 
 variation of the classical Ricci flow  equation that modifies the  unit  volume    constraint   of that equation to  a    
  scalar curvature    constraint. The resulting equations are named  the {\it conformal Ricci flow   equations} because  of the role
that 
  conformal geometry plays in constraining  the scalar curvature   and because these   equations are the  vector field      sum 
of a   conformal flow  equation and   a Ricci flow  equation. These new equations are given by
 
 \begin{eqnarray}\nonumber  \frac{\partial g}{\partial t}+2\big({\rm Ric}(g){\textstyle+\frac{1}{n}}g\big)=-p   g         
\\\nonumber 
 \;\qquad\qquad\quad  R(g)=-1   
\end{eqnarray} 

\medskip\noindent
 for a dynamically evolving metric $g$ and a scalar    non-dynamical  field $p$.
The conformal Ricci flow equations are analogous to  the Navier-Stokes equations of fluid mechanics,

 \begin{eqnarray}\nonumber  \frac{\partial v}{\partial t}+ \nabla _vv + \nu\Delta v &=&-{\rm grad}\,\,p        
\\ \nonumber 
  \; \qquad  \qquad  {\rm div}\,\,v&=&0  \,.
\end{eqnarray} 

\medskip\noindent
 Because of this analogy, the  time-dependent scalar  field $p$  is called a {\it conformal pressure} and, as for the  real 
physical pressure in fluid mechanics that serves to maintain the incompressibility of the fluid,  the conformal pressure serves as a
Lagrange multiplier to conformally deform  the metric flow so as to maintain the scalar curvature constraint.

The equilibrium points of the conformal Ricci flow equations  are Einstein metrics with Einstein constant $- {\textstyle\frac{1}{n}}
$.
 Thus the term  $-2 ({\rm
Ric}(g){\textstyle+\frac{1}{n}}g ) $   measures the
deviation of the flow from an equilibrium point and acts as a nonlinear restoring force. The
  conformal pressure     $p\ge0$ is 
     zero
at an equilibrium point  
and   positive  otherwise. The constraint force $-pg$ 
 acts pointwise orthogonally to   the nonlinear restoring force
  $-2{  ({\rm Ric}(g) +  \textstyle\frac{1}{n}} g)$ and      conformally deforms $g$     so that the
 scalar curvature is preserved.

A variety of  properties of the  conformal Ricci flow are discussed, including the reduced conformal Ricci flow, local existence and
uniqueness,  a variational formulation using a quasi-gradient of the Yamabe functional,   strictly monotonically decreasing global
and  local volume results,  and    applications   to
 3-manifold geometry. The geometry of the conformal Ricci flow is discussed as well as the remarkable analytic fact that the
constraint force  does not lose derivatives   and thus analytically the conformal Ricci equation is a bounded 
  perturbation of the classical unnormalized Ricci equation.
 Lastly, we discuss potential applications to Perelman's  proposed implementation of Hamilton's program to prove Thurston's  
    3-manifold  geometrization conjectures.
\end{abstract}

 \vspace{-.25in}

\vfill
 \noindent
 $\rule{2.65in}{.2mm}$ 

 \vspace{-.05in}
 
\noindent
{This paper will  appear 
in {\it 
 A Spacetime Safari: Papers in Honour of Vincent Moncrief},   a Special Issue of   {\it Classical and Quantum Gravity}, {\it
Institute of Physics Publishing},  January
2004. Online at http://www.iop.org.}



\maketitle

\begin{center}{\bf
 Dedicated to Professor Vincent Moncrief\\
on the  occasion of his 60th birthday  

\medskip
To Vince: Scholar, educator, colleague, friend
}\end{center}

{
 \tableofcontents
}
 \renewcommand{\contentsname}{Vincent4}
\newpage
\section{Introduction and background}
 \label{introduction}
  \subsection{The conformal Ricci flow  equations}
\label{crs}
 
We introduce a 
 variation of the {\bf classical Ricci flow equation} of Hamilton  \cite{ham82}  that modifies the   volume constraint  ${\rm
vol}(M,g_t)=1$  of the evolving metric $g_t$   to a 
   scalar curvature constraint $R(g_t)=-1$.
The resulting modified Ricci flow equations are named  the {\bf conformal Ricci flow  equations} because  of the role that 
  conformal geometry plays in maintaining the scalar curvature. Moreover, as we shall see in Sections~\ref{existenceCRF} and
\ref{nonortho}, the conformal Ricci flow equations are literally   the  vector sum  of a  smooth conformal   evolution equation
and   a densely defined Ricci evolution equation.

Since the volume of a Riemannian manifold $(M,g)$ is a  positive real number and since the scalar curvature is a real-valued function
on $M$,  the constraint on the scalar curvature is considerably more drastic than the volume constraint of the classical Ricci
flow equation. Thus the configuration space of the conformal Ricci flow equations    is considerably smaller than the configuration
space of  the classical Ricci flow equation.   We discuss why working on a smaller configuration space may
have
 advantages over working on a larger configuration space (see   Section~\ref{summary9}).

We take a unified point of view and show the similarities and differences between   the classical and conformal Ricci flow 
equations.   From this point of view, both  systems   have  the structure of a constrained dynamical system with    a Lagrange
multiplier that enters as a  time-dependent parameter field and  whose purpose is  to preserve  the constraints.
We also discuss various analogies between the conformal Ricci flow
equations and the Euler and Navier-Stokes equations for incompressible fluid flow. In particular, a time-dependent parameter (i.e.,
non-dynamical)  scalar field
$p\ge0$ naturally arises which we call the {\bf conformal pressure} 
and which has the   property that it     is zero
at an equilibrium point  
and   positive  otherwise.

A variety of  properties of the  conformal Ricci flow equations  are discussed, including the reduced conformal Ricci flow equation,
local existence and uniqueness,  a variational formulation using a quasi-gradient of the Yamabe functional,   strictly monotonically
decreasing global and  local volume results,  and  possible applications   to
 3-manifold geometry. 
The geometry of the conformal Ricci flow is discussed, leading to the remarkable analytic fact that the
constraint force  does not lose derivatives   and thus analytically the conformal Ricci equation is a bounded 
  perturbation of the classical unnormalized Ricci equation.
 Lastly, we discuss potential applications to Perelman's  proposed implementation of Hamilton's program to prove Thurston's  
    3-manifold  geometrization conjectures.

 Throughout this paper,       $M$ will denote  a smooth ($C^\infty$) closed (compact without boundary)  connected oriented
$n$-manifold, 
$n\ge3$. Let ${\cal D} = {\cal D} (M)={\rm Diff}(M)$ denote the infinite-dimensional group of smooth diffeomorphisms of $M$,  
$S_2=S_2(M)=C^\infty(T^*M\otimes_{\rm sym}T^*M)$   the space of smooth symmetric 2-covariant tensor fields on
$M$,   
   ${\cal M}={\rm Riem}(M) \subset S_2$   the
space of smooth Riemannian (positive-definite) metrics on $M$, and  
  $ {\cal F} =C^\infty(M, \mbox{\boldmath$R$} )$   the space of smooth real-valued functions on $M$, where $ \mbox{\boldmath$R$} $
denotes the real numbers.
 The group $ {\cal D} $ is an   infinite-dimensional
{\sc ilh}   (inverse limit Hilbert) Lie group, $ S_2 $ and  $ {\cal F} $    are infinite-dimensional
{\sc ilh}    linear spaces, and $ {\cal M} $ is  an 
open {\sc ilh} 
submanifold of $ S_2$ (i.e., $ {\cal M} $ is open in  $S_2$)
  with tangent space at $g\in {\cal M} $ given by  \begin{equation}
\label{LABEL} T_g {\cal
M}=\{g\}\times S_2 
\approx S_2\,,
\end{equation} 
  which we identify with $S_2$ (for more information about {\sc ilh} spaces, see Omori
\cite{omo70}, Ebin \cite{ebin70}, and Ebin-Marsden \cite{em70}).

For each   $g\in {\cal M} $  there exists a unique   {\bf volume element}
$d
\mu _g$  on $M$ 
 determined  by
$g$ and  the orientation of $M$.
 In local coordinates $(x^1,\ldots,x^n)$,   $d \mu _g=\sqrt{\det{g_{ij}}}\,\,dx^1\wedge\ldots\wedge dx^n$.
   Let ${\rm
vol}(M,g)=\int_M
d\mu _g$ denote the {\bf volume} of the Riemannian manifold 
$(M,g)$, let $ \mbox{\boldmath$R$} ^+=(0,\infty)$,   let
\begin{eqnarray*}
{\rm vol}: {\cal M} \longrightarrow \mbox{\boldmath$R$} ^+ \;,\quad g \longmapsto \int_M d\mu _g= {\rm vol}(M,g)
\end{eqnarray*}
denote the {\bf volume functional} on $ {\cal M} $, and let
   \[
{\cal M} ^1=\{\,g\in {\cal M}
\mid{\rm vol}(M, g)= 1\,
\}= {\rm vol}  ^{-1} (1)
\] denote   the {\bf space of  Riemannian metrics with unit volume}.

For a Riemannian metric $g$ in $ {\cal M} $,    let $K(g)$ denote the {\bf sectional curvature} and let  ${\rm Riem}(g)$ denote the 
{\bf Riemann-Christoffel curvature tensor} defined on vector fields $X,Y$ and $Z$ on $M$ by  \begin{equation}
\label{LABEL} {\rm Riem}(g)(X,Y,Z)=\nabla_X\nabla_YZ-\nabla_Y\nabla_XZ-\nabla_{[X,Y]} Z \,.
\end{equation} Let
 ${\rm Ric}(g)\in
S_2$ denote  the {\bf Ricci curvature tensor},    $ R(g)={\rm tr}_g({\rm Ric}(g)) =g ^{-1}({\rm Ric}(g)) \in {\cal F}  $    the {\bf
scalar curvature},   where $ {\rm tr}_g$ is the metric trace and $ g ^{-1} $ is the inverse of $g$,
 ${\rm  Weyl} (g)$   the {\bf Weyl conformal curvature  tensor},
 ${\rm Ric}^T\! ( g )={\rm Ric}  ( g
 )-{\textstyle \frac{1}{n}} R(g ) g\in S_2$      the {\bf traceless part of the Ricci tensor},
and  
  ${\rm Ein}(g)=  {\rm Ric}(g)- {\textstyle\frac{1}{2}}R(g)g
\in S_2$   the  {\bf Einstein tensor} of $g$. In local coordinates, $ ({\rm Riem}(g))^i{}_{jkl}=  R ^i{} _{jkl}$, 
  $({\rm Ric}(g))_{ij}=  R_{ij}=R^a{} _{iaj}$,
  $R(g)= g^{ij}R_{ij}  $, $ ({\rm Weyl}(g))^i{}_{jkl}=  W ^i{} _{jkl}$, 
$ ({\rm Ric}^T\! ( g ))_{ij}=R_{ij}  - {\textstyle\frac{1}{n}}R g_{ij} $, and
$ ({\rm
Ein}(g))_{ij}=R_{ij}  - {\textstyle\frac{1}{2}}R g_{ij} $.
(The $T$ designating traceless in 
$(\,\cdot\,)^T$   is unrelated to the $T$ that occurs in  the semi-open intervals $[0,T)$ introduced later.)

We
  define a {\bf hyperbolic metric} on $M$ as a Riemannian metric $g\in {\cal M} $  with  constant negative  sectional curvature,
  not necessarily
$-1$.   This definition is slightly more general than the usual definition in which   $K(g)=-1$. By    allowing
$K(g)$ to float, we can normalize it by other conditions, such as requiring that the scalar curvature  
$R(g)=-1$ in which case   $K(g)=-\frac{1}{n(n-1)}\,$,
in contrast to having $K(g)=-1$ and   $R(g)= -n(n-1)$. We remark that with this definition of hyperbolic metric, Mostow rigidity
only determines a hyperbolic metric up to isometry and homothety. However, by normalizing the {\it scalar} curvature of hyperbolic
metrics to $-1$, rigidity of hyperbolic metrics up to isometry is recovered.

  For $ g\in {\cal M} $ and 
 $  \rho $   a tensor field on  $M$, let $|\rho |_g=( \rho \cdot \rho )^{1/2}\in {\cal F} $ denote the pointwise $g$-metric norm of
$
\rho
$ where center dot ``$\,\,\cdot\,$" denotes the $g$-metric contraction to a scalar.  For example,   
$ |{\rm Ric}( g ) |^2  _{
g }={\rm Ric}( g )\cdot {\rm Ric}( g )= g^{ik}g^{jl}R_{ij}R_{kl}$
is the   pointwise squared norm of the Ricci tensor and   
  $ | {\rm Riem}(g) |^2_g= g_{ai}g^{bj}g^{ck}g^{dl}R^a{}_{bcd}R^i{}_{j k l }$ is the pointwise squared  norm of the
Riemann-Christoffel curvature tensor of
$g$.

On $ {\cal M} $ there is a natural $L_2$ weak infinite-dimensional Riemannian metric $ {\cal G} $, first constructed by Palais
\cite{pal61} on November 10, 1961 and first published in   Ebin  
\cite{ebin70}, given by
\begin{equation}
\label{naturalmetric}  {\cal G} _g:T_g {\cal M} \times T_g {\cal M}\approx S_2\times S_2  \longrightarrow \mbox{\boldmath$R$}
\;,\quad (h,h) \longmapsto \int_Mh\cdot h\, d \mu _g\,.
\end{equation} 
 Let $( {\cal M} , {\cal G} )$ denote this weak infinite-dimensional Riemannian manifold.

Let
 \begin{equation}
\label{eq:Put equation label here}{\cal M} _{-1}=\{\,g \in {\cal M} \mid R(g) =-1\,\}
\end{equation} denote the {\bf subspace of Riemannian metrics with constant scalar
curvature} $-1$. Then     $ {\cal M}_{-1} \subset {\cal M} $ with the induced weak Riemannian metric  $ {\cal G} _{-1} $ is a weak
Riemannian 
 submanifold of
$ ( {\cal M} , {\cal G} )$ which we denote $( {\cal M}_{-1} , {\cal G} _{-1} )\subset ( {\cal M} , {\cal G} )$.
 The differential topology (but not the geometry) of this  space has been studied in detail in
Fischer-Marsden \cite{fm75b} and we summarize the needed information in Section~\ref{space1}.

The conformal Ricci flow equations (with smooth initial conditions) are   defined as follows.
\begin{definition}  
\label{parabolic31} Let $M$ be a smooth   closed   connected oriented $n$-manifold, $n\ge3$.
 The {\bf conformal Ricci flow  equations}  on $M$ are
defined by the 
  equations
 \begin{eqnarray}\label{parabolic} \frac{\partial g}{\partial t}+2\big({\rm
Ric}(g){\textstyle+\frac{1}{n}}g\big)=-p   g          \\
\label{elliptic} \;\;\;\qquad\qquad\quad  R(g)=-1    
\end{eqnarray} for       curves 
 \[
\label{LABEL} g:[0,T) \longrightarrow {\cal M} \;,
\quad t \longmapsto  g(t) \,,\quad \mbox{with}\quad g(0)=g_0\in {\cal M}_{-1}\,,    \;\mbox{and}  
\]
 \[
\label{LABEL} p:[0,T) \longrightarrow {\cal F}  \;,
\quad\;\, t \longmapsto  p (t)\,,
\] continuous on the semi-open  interval $[0,T)$,  
    $0<T\le\infty$, and smooth on the open interval $(0,T)$.
  \hfill$\rule{2mm}{3mm}$ 
 \end{definition} 

\noindent {\bf Remarks: }
\begin{enumerate}
\item
 Because of the constraint equation $R(g)=-1$ and the initial condition $g_0\in {\cal M}_{-1}\, $, the curve
$g:[0,T) \longrightarrow {\cal M}$ actually lies in $ {\cal M}_{-1} $. However, if $g$ is viewed as such, then to  assert that $g$
restricted to  $(0,T)$ is a smooth curve requires the result given below (Theorem~\ref{main1})  that $ {\cal M}_{-1} $ is a smooth
submanifold of $ {\cal M} $. Thus for now, for differentiability purposes, we take $g$ to be smooth as a curve in  $ {\cal M} $ 
and   constrained to  the subspace
$ {\cal M}_{-1}
$ of $ {\cal M} $.
 \item
As we shall see in Proposition~\ref{maxprin}  (see also the Remark following that Proposition), for non-static flows $g$, the
curve  
 $p:[0,T)
\longrightarrow {\cal P}
$ actually lies in the open subspace of positive real-valued functions $ {\cal P} =C^\infty(M, \mbox{\boldmath$R$} ^+)\subset
C^\infty(M, \mbox{\boldmath$R$})= {\cal F} $.
  \end{enumerate} 
\medskip

 We refer to $(\ref{parabolic})$  as the {\bf  evolution  equation} (of the conformal Ricci flow equations), 
  $(\ref{elliptic})$   as the {\bf   constraint equation} (of the conformal Ricci flow equations),
 the curve 
$ g:[0,T) \rightarrow {\cal M} $ as  the {\bf conformal Ricci flow},   and the curve
 $ p:[0,T) \rightarrow {\cal F} $ as the {\bf conformal pressure}.
  Because of the {\bf quasi-parabolic} nature of the
equations (see Section~\ref{existenceCRF}), we do not expect any backward time of existence  and so we refer to a solution with a 
  positive semi-infinite time of existence
$[0,\infty)$
  as   an {\bf
all-time solution}.
 
 \subsection{A brief survey of the geometry of the conformal Ricci flow equations}

  As a consequence of the constraint equation, 
  the conformal Ricci flow takes place in $ {\cal M}_{-1}$ which   
 is    a   non-empty  infinite-dimensional closed submanifold of  $ {\cal M} $ (see Theorem~\ref{main1}).
 As a consequence of this constraint, it is 
 of importance to note  that the conformal pressure
$p
$ is {\it not}  on the same footing as the  dynamical  metric field
$g$ and thus
        {\it no} initial value of $p$ is  given. Rather, $p$ serves as a time-dependent Lagrange multiplier  and
the term $-p g$ acts as the   constraint force       necessary to preserve the scalar curvature constraint (see also the discussion
in  Section~\ref{nonortho}). Consequently,  
$p
$ must solve  a time-dependent elliptic partial differential equation  as the metric evolves. This situation is completely analogous
to the incompressible Euler or Navier-Stokes equations where the real physical pressure serves as   a
Lagrange multiplier in order to maintain the divergence-free constraint of the dynamical vector field
$v$ (see Ebin-Marsden \cite{em70}).  Because of this analogy, we   think of the  non-dynamical   scalar field $p$ as a {\bf
conformal pressure} because it pointwise conformally deforms $g$ so as to maintain the scalar curvature constraint.

The conformal Ricci flow equations share  characteristics with both the quasilinear  elliptic-hyperbolic  incompressible Euler
equations and the semilinear elliptic-parabolic  incompressible  Navier-Stokes equations, sharing      the  elliptic characteristics
of both equations,  the quasilinear characteristics of the Euler equations,  and   the parabolic characteristics of the 
Navier-Stokes equation (see also Section~\ref{existenceCRF}).

In fact, one of the design criteria for the conformal Ricci flow equations  was that     the structural similarities to the
Euler and Navier-Stokes equations for an  incompressible  fluid were apparent. For example,  Euler's equations for a  unit density
incompressible ideal fluid with
 time-dependent   velocity vector field $v$   on a fixed
Riemannian manifold $(M,g)$  are 
 \begin{eqnarray}\label{euler1} \frac{\partial v}{\partial t}+\nabla_vv=-{\rm grad}\,\, p 
\;, \qquad  {\rm div}\,v=0 \;,\qquad v(0)=v_0 \,, 
\end{eqnarray} where $p$ is the time-dependent   pressure of  the fluid. 
These equations can also be written as
\begin{equation}   
\label{projectiona}  \frac{\partial v}{\partial t}+P(\nabla_vv)=0\;,\qquad  {\rm div}\, v_0=0\;,\qquad v(0)=v_0\,,
\end{equation} where $P$ is the projection of  a vector field onto its divergence-free part.
Similarly, for a viscous fluid,  the incompressible Navier-Stokes equations are
  \begin{eqnarray}\nonumber  \frac{\partial v}{\partial t}+ \nabla _vv + \nu\Delta v &=&-{\rm grad}\,\,p        
\\ \nonumber 
  \; \qquad  \qquad  {\rm div}\,\,v&=&0  \,,
\end{eqnarray} 
where the Laplacian on (1-forms metrically associated with)  vector fields is $ \Delta=  \delta d+ d \delta $
and  the kinematic viscosity is $ 
 \nu>0$ 
(see   Ebin-Marsden \cite{em70}, p.\ 161, and Taylor \cite{tay96c}, p.\ 493, for the ``correct"  Navier-Stokes equations on a
Riemannian manifold; here the simpler form is sufficient for our purposes).  These equations can be thought of as a vector
field sum of the heat equation
\[
 \frac{\partial v}{\partial t}  =-  \nu\Delta v
\]
and the Euler equation.  A nonlinear Trotter product formula is then  applied  to prove existence and uniqueness
of solutions (Ebin-Marsden
\cite{em70}). In Section~\ref{nonortho},  we use a similar approach to prove existence and uniqueness to the 
conformal Ricci flow equations.

  Equation 
(\ref{euler1})  
  should be compared with   (\ref{parabolic}--\ref{elliptic}) and (\ref{projectiona}) should be compared with  the {\bf projection
form}  (\ref{proj3})  of the conformal Ricci flow equations  
 \begin{eqnarray}
\label{parabolic54} \frac{\partial g}{\partial t}+2\widetilde P_g(  {\rm Ric}(g) ) =0  \,,
     \end{eqnarray} where $\widetilde P_g$ is a non-orthogonal projection from $T_g {\cal M} $ to 
  $T_g {\cal M}_{-1} $ 
   described in  Section~\ref{nonortho}.
 
 As we shall see in Proposition~\ref{equilprop}, the
equilibrium points of the conformal Ricci flow equations occur at Einstein metrics with Einstein constant $ -{\textstyle\frac{1}{n}}
$, i.e., when 
$ {\rm Ric}(g) +  {\textstyle\frac{1}{n}} g=0$. Thus the term    $ -2 ({\rm Ric}(g) + {\textstyle\frac{1}{n}} g)$  
  can be viewed as a measure of the deviation from an equilibrium point or, in a rough mechanical analogy, can be thought of as a
nonlinear restoring force (of course, since
 these equations are first order in time, this analogy is only valid for Aristotelian, rather than Newtonian,   mechanics).
On the other hand, the    conformal pressure term 
$-pg$ is the constraint force necessary to counter-balance the   conformal component  $pg$ of the 
nonlinear restoring force 
 $ -2 ({\rm Ric}(g) + {\textstyle\frac{1}{n}} g)$. Thus the constraint force    acts by   conformally deforming $g$     so that
 the scalar curvature is preserved. The constraint force  $-pg$  acts pointwise orthogonally to  $ -2 ({\rm Ric}(g) +
{\textstyle\frac{1}{n}} g)$ and    keeps the flow in the submanifold
$  {\cal M}_{-1}
$  (see   Figure~1 in Section~\ref{nonortho} and the ensuing discussion).  
  Moreover, 
  as we shall see in Proposition~\ref{maxprin},  the conformal pressure  $p$ is, as expected (or, as designed), zero
at an equilibrium point  
and strictly  positive  otherwise, and is  calculated in    (\ref{global}).

\medskip
\noindent {\bf Remarks: }
\begin{enumerate}
\item
 When we refer to   the conformal Ricci flow  equations, we shall
always assume all of the above conditions regarding $M$, $n$,  $g$, and $p$.
For the reasons for the restriction $\dim M=n\ge3$, see Remark~(i) following Proposition~\ref{equilprop}. 
\item
For the most part $g$ and $p$ will denote time-dependent objects and we sometimes 
write this explicitly as $g=g(t)=g_t$ or $p=p(t)=p_t$. A  subscript $t$ will never mean the partial
$t$-derivative,  which we will always write explicitly. 
\item
The $2$ that appears in (\ref{parabolic}) as well as in (\ref{parabolic22})
below comes about because in local coordinates  the second order quasilinear terms of the Ricci tensor 
$R_{ij}=({\rm
Ric}(g))_{ij}$ are given by
\begin{equation}
\label{riccilocal3}
\fl 
 \quad R_{ij}=-{\textstyle\frac{1}{2}} g^{kl}\frac{\partial^2g_{ij}}{\partial x^k\partial x^l}
+{\textstyle\frac{1}{2}} g^{kl}\frac{\partial^2g_{ki}}{\partial x^l\partial x^j}
+{\textstyle\frac{1}{2}} g^{kl}\frac{\partial^2g_{lj}}{\partial x^k\partial x^i}
-{\textstyle\frac{1}{2}} g^{kl}\frac{\partial^2g_{kl}}{\partial x^i\partial x^j}+Q_{ij} \,,
\end{equation} 
 where 
$Q_{ij}=Q_{ij}(g_{kl},\frac{\partial g_{mn}}{\partial x^p })$ represents the  lower order terms of $R_{ij}$,  which are
quadratic in the first  order derivatives of $g$ and rational in $g$. Locally,  the evolution equation (\ref{parabolic}) can be
written as 
 \[\fl  \frac{\partial g_{ij}}{\partial t}-   g^{kl}\frac{\partial^2g_{ij}}{\partial x^k\partial x^l}
+  g^{kl}\frac{\partial^2g_{ki}}{\partial x^l\partial x^j}
+ g^{kl}\frac{\partial^2g_{lj}}{\partial x^k\partial x^i}
-  g^{kl}\frac{\partial^2g_{kl}}{\partial x^i\partial x^j}+Q_{ij}   +{\textstyle\frac{2}{n}}  g_{ij}
 =-pg_{ij}\;.\] 
Thus  if we ignore for the moment  the three second order {\bf index  mixing terms}, then
\begin{equation} 
\label{lochar}  \frac{\partial g_{ij}}{\partial t}-   g^{kl}\frac{\partial^2g_{ij}}{\partial x^k\partial x^l}
 +Q_{ij}   +{\textstyle\frac{2}{n}}  g_{ij}
\approx -pg_{ij}\,,
\end{equation} so that the left hand side is a  quasilinear strictly parabolic (or non-mixing)  heat operator for $g_{ij}$.
In fact, in a system of time-dependent   coordinates  which we call {\bf moving harmonic coordinates},    (\ref{lochar}) 
is locally  the   evolution equation of the  conformal Ricci flow  equations.
Moreover, these time-dependent coordinates can be used as  a basis for proving  existence and uniqueness of solutions   for  
the classical Ricci  flow equation   (see Fischer \cite{fis04}).
 \item
Since the conformal Ricci flow satisfies the constraint $R(g)=-1$, $ {\rm Ric}^T\! (g) = {\rm Ric}(g)  -\frac{1}{n}R(g)g={\rm
Ric}(g)  +\frac{1}{n} g$  and so (\ref{parabolic}) can also be written in the useful form \begin{equation}
\label{parabolicv} \frac{\partial g}{\partial t}+2 {\rm Ric}^T \!(g)=- p  
g        \,.
\end{equation}  
 \end{enumerate}

\subsection{Some comparisons with the   classical Ricci flow equation}
\label{comparison}

For $ g\in {\cal M} $, let $\bar R(g)=  \int_MR(g) d \mu _g/ {\rm vol}(M,g) $ denote  the   volume-averaged total scalar curvature
of 
$g$. To put the conformal Ricci flow equations into perspective, we consider the 
  {\bf classical (normalized)   Ricci flow  equation} of Hamilton \cite{ham82}, 
  \begin{equation}\label{parabolic100}\;\; \quad\frac{\partial g}{\partial t}=-2 {\rm Ric}(g)+{\textstyle\frac{2}{n}}  \bar R( g  
)  g    \,,
 \end{equation}
 from a slightly different point of view than that in which  it
is usually considered. Equation (\ref{parabolic100})  is called normalized because the volume is conserved;
Hamilton's {\bf unnormalized Ricci flow equation}, where the volume is not preserved, is
\begin{equation}
\label{unnormalRicci} \frac{\partial g}{\partial t}=-2 {\rm Ric}(g)\,. 
\end{equation} 

We augment 
  the classical Ricci flow equation (\ref{parabolic100}) and write it  as a constrained dynamical system 
 \begin{eqnarray}\label{parabolic22}\;\; \quad\frac{\partial g}{\partial t} &=&-2\,{\rm Ric}(g)+ c g          \\
\label{elliptic22}{\rm vol}(M, g) & = & 1 
\end{eqnarray} for      curves 
\[
\label{curve1} \label{LABEL} g:[0,T) \longrightarrow {\cal M} \;,
\quad t \longmapsto  g(t)\quad\mbox{with} \quad  g(0) \;  =  g_0\in {\cal M} ^1\, ,\mbox { and}      
\]\[
\label{curve2} \label{LABEL} c:[0,T) \longrightarrow  \mbox{\boldmath$R$}   \;,
\quad \;\;t \longmapsto c (t)\,,
\] continuous on the semi-closed interval $[0,T)$ and smooth on 
  $(0,T)$, 
 $0<T\le\infty$.
 
The curve
$c:[0,T) \rightarrow  \mbox{\boldmath$R$} $ acts as  a time-dependent Lagrange multiplier to preserve the constraint equation
${\rm vol}(M, g) = 1   $.  As with any Lagrange multiplier, it can formally and in this case actually be eliminated. Taking the
time
derivative of the constraint equation (\ref{elliptic22})  and using the evolution   equation  (\ref{parabolic22})  yields
\begin{eqnarray} \nonumber 0&=& \frac{d }{dt}\,{\rm vol}(M,g ) =\frac{d  }{dt}    \int_M d\mu _{g } = 
\int_M \frac{\partial }{\partial t}d\mu
_{g }=\int_M {\textstyle\frac{1}{2}}{\rm tr}_{g }\left(\frac{\partial g }{\partial
t}
\right)d \mu _{g } \nonumber 
\\&=&\int_M {\textstyle\frac{1}{2}}{\rm tr}_{g }\left(  -2\,{\rm Ric}(g)+c g        
\right)d \mu _{g } \nonumber   =-\int_M   R(g)d \mu _{g }+{\textstyle\frac{n}{2}}  c \int_M        
 d \mu _{g } \nonumber \\  &=&-   R_{\rm total}(g) +{\textstyle\frac{n}{2}}  c\, {\rm vol}(M,g )  \,,  
  \label{neww}
\end{eqnarray} 
where $ R_{\rm total}(g)=\int_M   R(g)d \mu _{g }$ is the  {\bf total} (or integrated) scalar curvature. Thus  
\begin{equation}
\label{cval0} c=  {{\textstyle\frac{2}{n}}  } \frac{ R_{\rm total}(g)}{{\rm vol}(M,g ) }=
{{\textstyle\frac{2}{n}}  }  \bar R (g) 
\end{equation}where $\bar R (g) =\frac{  R_{\rm total}(g)}{{\rm vol}(M,g ) }    $ is the {\bf volume-averaged total scalar
curvature of}
$g$.
Since $c $  involves integrals over $M$, it is a  global,
or nonlocal,   function of
$g$ and thus depends on
$g$ at every point of
$M$ (see also Remark~(iii) after Proposition~\ref{maincons}).

 Using  (\ref{cval0}),   both   $c$  and the constraint equation can be eliminated from the   system
(\ref{parabolic22}--\ref{elliptic22})    which 
  can then be  written in its usual {\bf reduced  form}, 
\begin{eqnarray} \label{usual} \frac{\partial g}{\partial t}\quad&=&-2\,{\rm Ric}(g)+  {{\textstyle\frac{2}{n}}  }
\bar  R (g) g \;,\qquad g_0\in {\cal M} ^1 \,,         
\end{eqnarray} where now the constraint ${\rm vol}(M,g_0)    =  1$  applies only to the initial condition $g_0$ rather than the
entire flow $g_t$,  $t\in [0,T)$.

Note that the reduced form (\ref{usual}) is 
achieved by using only the time derivative of the constraint equation,  (\ref{neww}), and not the constraint equation
${\rm vol}(M, g) = 1   $  itself (see also
Section~\ref{observation}  below).

 The constraint on the entire flow  $g_t$ is   recovered from (\ref{usual}) since  
\begin{eqnarray} \fl \frac{d}{dt} \,{\rm vol}(M,g ) & = & {\textstyle\frac{1}{2}}\int_M {\rm tr}_{g }\left(\frac{\partial g
}{\partial t}
\right)d \mu _{g } 
={\textstyle\frac{1}{2}}\int_M {\rm tr}_{g }\left(  -2\,{\rm Ric}(g)+  {{\textstyle\frac{2}{n}}  }
\bar  R (g) g            
\right)d \mu _{g } \nonumber \nonumber           \\\fl  & = &   -\int_M   R(g)d \mu _{g }+{\textstyle\frac{n}{2}}  
{{\textstyle\frac{2}{n}}  }
\bar  R (g)        \int_M        
 d \mu _{g }   =
-   R_{\rm total}(g) + R_{\rm total}(g)=0\,,
\label{cancel0}
\end{eqnarray} so that from  the initial constraint $g_0\in {\cal M} ^1$,    
${\rm vol}(M,g_t  ) =1$.
 Thus a solution of the    reduced equation $(\ref{usual})$   stays in $ {\cal M} ^1$ if it starts in $ {\cal M}
^1$ and thus the   constrained  system $(\ref{parabolic22}$--$\ref{elliptic22}) $   is equivalent to
the reduced equation $(\ref{usual})$  with initial
conditions in $ {\cal M} ^1$. 
 This same theme with the appropriate spaces will appear for     the conformal Ricci flow equations 
   (see Section~\ref{elimin}).

\subsection{A further simplification of the classical Ricci  flow equation}
\label{observation} As we have seen in Section~\ref{comparison}, to get to   the usual
formulation of the classical Ricci flow,   the Lagrange multiplier $c$ is   eliminated by   
using the derivative of the constraint equation, $\frac{d}{dt} {\rm vol}(M,g) =0$. However, the undifferentiated constraint equation
itself 
$ {\rm vol}(M,g) =1$ is not    explicitly
used. 
  If we do   apply the constraint, then
\begin{equation}
\label{fullyred}   c= {{\textstyle\frac{2}{n}}  }  \bar R (g)= {{\textstyle\frac{2}{n}}  } \frac{ R_{\rm total}(g)}{{\rm vol}(M,g )
}= {{\textstyle\frac{2}{n}}  }  R_{\rm total}(g) 
\label{eliminatec}
\end{equation} 
 and      the classical Ricci equation  (\ref{usual}) further reduces to   the {\bf fully reduced equation}
\begin{eqnarray} \label{usual0} \frac{\partial g}{\partial t}\quad&=&-2\,{\rm Ric}(g)+  {{\textstyle\frac{2}{n}}  }
  R_{\rm total} (g) g        \;,\qquad g_0\in {\cal M} ^1 \,.
\end{eqnarray} Since 
 the constraint equation itself   and not just its time derivative have   been used in deriving 
 (\ref{usual0}), solutions to  (\ref{usual0}) no longer satisfy 
  $\frac{d}{dt} {\rm vol}(M,g) =0$.
However, the constraint on the entire flow is still recoverable as follows. 
Let 
$\tilde v(g) = {\rm vol}(M,g) -1$. 
Then if $g$    satisfies
(\ref{usual0}),   
\begin{eqnarray}\fl\frac{d\tilde v(g)}{dt}  & = & \frac{d}{dt} \,{\rm vol}(M,g )= {\textstyle\frac{1}{2}}\int_M {\rm tr}_{g
}\left(\frac{\partial g }{\partial t}
\right)d \mu _{g } 
={\textstyle\frac{1}{2}}\int_M {\rm tr}_{g }\Big(  -2\,{\rm Ric}(g)+  {{\textstyle\frac{2}{n}}  }
  R_{\rm total} (g) g            
\Big)d \mu _{g } \nonumber \nonumber           \\\fl   & = &   -\int_M   R(g)d \mu _{g }+{\textstyle\frac{n}{2}}  
{{\textstyle\frac{2}{n}}  }
 R_{\rm total} (g)        \int_M        
 d \mu _{g }   =
-   R_{\rm total}(g) + R_{\rm total}(g) {\rm vol}(M,g)   \nonumber \\\fl
& =&
   R_{\rm total}(g)(  {\rm vol}(M,g) -1) =  R_{\rm total}(g)\tilde v(g)\nonumber 
\label{cancel04}\,.
\end{eqnarray} 
Thus if we let $\tilde v(t)= \tilde v(g(t)) $ and  $R_{\rm total}(t)= R_{\rm total}(g(t))$, we see that
$\tilde v(t)$ solves the linear first-order   non-autonomous  ordinary differential equation
 \begin{equation}
\label{ode} \frac{d\tilde v( t )}{dt}   = R_{\rm total}( t )\tilde v ( t )\,.
\end{equation}
Since the initial metric $g_0$ has  ${\rm vol}(M,g_0 ) =1$,   
$\tilde v( 0 )=\tilde v_0=0$, so  by uniqueness of  solutions of  (\ref{ode}) (or by the explicit solution $\tilde v(t)=\tilde v
_0e^{\int_0^tR_{\rm total}(t')dt'}
 $),
  $\tilde v( t )= {\rm vol}(M,g(t))
-1\equiv0$. Thus the initial constraint is    maintained by the fully reduced equation (\ref{usual0})  and thus  
   (\ref{parabolic22}--\ref{elliptic22}) with $g(0)=g_0\in {\cal M} ^1$,    (\ref{usual}),  and (\ref{usual0}) are all equivalent.

As we shall see, the  fully   reduced   equation (\ref{usual0}) is the one  most similar to the  fully   reduced conformal
Ricci equation where a linear (partial differential)  equation  similar in
structure to (\ref{ode}) appears  (see (\ref{kac}) and Remark~(i) following) to insure that the constraint is maintained on the
entire flow if it is satisfied by the initial condition.

\subsection{The conformal Ricci flow on manifolds of Yamabe type $-1$}
\label{conformal4}
 In order to explain why we have chosen the name   conformal Ricci flow, we
 recall the following terminology introduced by     Fischer-Moncrief \cite{fm94b}. 
 
  \begin{definition}[The Yamabe type of a manifold]
\label{defyam} Let  $M$  be a    closed  connected       $n$-manifold,
$n\geq3$.  Then

 $(1)$  $M$ is of {\bf  Yamabe-type $-1$} 
   if   
$M$ admits no metric with $R(g)=0;$

$(2)$  $M$ is of {\bf Yamabe-type $0$} 
  if   
$M$ admits a metric with $R(g)=0$  but no    metric with $R(g)=1 ;$

 $(3)$ $M$ is of {\bf   Yamabe-type}  $1$  
   if   
$M$ admits a metric with $R(g)=1  $. 
\end{definition}

The   definition of Yamabe type  partitions the class of closed  $n$-manifolds, $n\ge3$, into 
three classes that are mutually exclusive and exhaustive   (see~\cite{fm94b} for details).
  We note in particular that the three categories of closed manifolds that admit   metrics of  
constant positive, zero, and constant negative sectional curvature  fall into the categories of
Yamabe type $1$, 0, and $-1$, respectively.

The pointwise multiplicative Abelian group     $ {\cal P}  =C^\infty(M, \mbox{\boldmath$R$} ^+)$    
of smooth positive real valued functions on $M$  acts freely, smoothly, and properly on $ {\cal M} $. The
resulting orbit space
$ {\cal M} / {\cal P} $ is the space of {\bf (pointwise) conformal classes} on $M$.  In particular, we note that $ {\cal M} /
{\cal P} $ is a contractible manifold (see Fischer-Moncrief
\cite{fm96} for more information).

The key  fact is that if
  $M$ is of Yamabe type $-1$, 
then every
Riemannian metric on
$M$ is uniquely pointwise conformally deformable to a metric with scalar curvature
$=-1$.
Thus  
  $ {\cal M} / {\cal P} $ and $ {\cal M} _{-1}$ are in bijective correspondence and in 
fact are  {\sc ilh} diffeomorphic.
  Thus   $ {\cal M} _{-1}$ is a
 representation 
of the space of pointwise 
  conformal structures  $ {\cal M} / {\cal P} $ and results regarding the conformal Ricci flow on $ {\cal M}_{-1} $ can   be
interpreted in terms of conformal geometry.

When $M$ is not of Yamabe type $-1$, then $ {\cal M}_{-1} $ represents only one component of  $ {\cal M} / {\cal P} $  and thus the
conformal interpretation needs modification.  We emphasize, however, that as we have done in  this paper, one can work    
directly on
$ {\cal M}_{-1}
$   as a space in its own right  without       restricting to manifolds of Yamabe type $-1$. If, however, we wish to interpret our
results
  on the space of (pointwise) conformal geometries  $ {\cal M} / {\cal P} $, then we would need to either restrict to manifolds of
Yamabe type
$-1$ or augment our results to include other subspaces of $ {\cal M} $.
 
Since  the conformal Ricci equation was designed with conformal geometry in mind,        applications to conformal geometry  
  will, as expected,  be of more interest on manifolds of 
  Yamabe type~$-1$ (see for example Proposition~\ref{increasescale} and also  Sections~\ref{conclusions} and \ref{end0}).

\subsection{Summary}
\label{summary9}
 Comparing the  classical    and conformal Ricci flow equations,   we have the 
following:
  \begin{enumerate}
\item
The constraint equation  changes from 
  ${\rm vol}(M,g)=1 $ for  the classical    Ricci flow to $R(g)=-1$ for the conformal Ricci flow with the concomitant
change of the   configuration space from $ {\cal M} ^1 $ to  ${\cal M}_{-1} $. Since $ {\cal M} ^1 $ is a codimension-1 submanifold
of
$ {\cal M} $ whereas $ {\cal M}_{-1}
$ is a
codimension-$C^\infty(M, \mbox{\boldmath$R$} )$ submanifold 
 of $ {\cal M} $, $ {\cal M}_{-1}
$ is a much smaller configuration space than $ {\cal M} ^1$.
\item To compensate for this smaller configuration space, the Lagrange multiplier in the 
  evolution equation adjusts  from  being a time-dependent constant $c=c(t)$   (independent of $x$)
in the classical Ricci system to a time-dependent real-valued  
function  $p=p(t,x)$   
 in the conformal Ricci system.
 \item In the augmented classical Ricci flow equations (\ref{parabolic22}--\ref{elliptic22}),   solving for the Lagrange
multiplier and eliminating the constraint equation involves solving a time-dependent linear inhomogeneous  algebraic equation
(\ref{neww}) whereas in the conformal Ricci flow equations solving for the Lagrange multiplier and eliminating the constraint
equation involves solving a time-dependent linear inhomogeneous  elliptic partial differential equation~(\ref{elliptic6}). Thus,
although  both the classical and conformal Ricci equations are quasilinear quasi-parabolic systems (see
Section~\ref{existenceCRF}),  the conformal Ricci equations have an  additional   linear   elliptic equation and thus overall are a
linear-elliptic quasilinear  quasi-parabolic system.

\item For the classical  and conformal Ricci flow equations, 
    the volume and scalar curvature behave somewhat oppositely. For the classical Ricci flow equation, the volume
is preserved,
  $  {\rm vol}(M,g) =1$,  but for non-static flows the scalar curvature is not preserved,  whereas for
 the conformal Ricci
equation, the  scalar curvature is preserved,
$
 R(g) =-1 $, but for  non-static flows  the volume is not preserved. For more specific information,
see Proposition~\ref{increasescale}, Table~1 following that Proposition,    and   Proposition~\ref{volel}. 

\item As we shall see  in Example~\ref{ex2},     $ {\rm Ric}(g) $ is orthogonal to $ {\cal M}_{-1} $.  Thus the classical
unnormalized Ricci flow (\ref{unnormalRicci}), at a metric $g\in {\cal M}_{-1} $, is orthogonal to
$  {\cal M}_{-1} $, whereas  the conformal Ricci flow  is
tangential to $ {\cal M}_{-1} $ (see Figure~\ref{figure1} in Section~\ref{nonortho}). 
Thus in some sense the  classical
unnormalized Ricci flow acts    like a Riemannian  gradient  as it is orthogonal to the level hypersurface $ {\cal M}_{-1} $, 
whereas the conformal Ricci flow acts like a symplectic gradient as it  is tangent  to the level hypersurface $ {\cal M}_{-1} $.
  \end{enumerate}

We now discuss why 
from the point of view of geometry having a smaller configuration space is  potentially better.
The classical Ricci flow was  designed to search out Einstein metrics. But if one is looking for something, it makes sense
to look in as small a space as possible if it  is known that the object  of interest is in that smaller space.  Since an Einstein
metric with negative Einstein constant can always be scaled so as to lie in  $ {\cal M}_{-1} $, we know that what we are looking
for must lie in
$ {\cal M}_{-1}
$. Thus    
  it makes sense at the outset to restrict the problem  to the smaller space $ {\cal M}_{-1} $ and in
  fact the  conformal Ricci flow equations were   designed   so as to  restrict to this smaller configuration space $ {\cal
M}_{-1}
$ from the beginning.

On the other hand, since 
  {\it any} metric (and not just an Einstein metric) can be scaled to lie in $ {\cal M} ^1$,
restricting the search for Einstein metrics to  $ {\cal M} ^1$ does not significantly reduce the search space.

As a second  motivation for restricting to the smaller  space $ {\cal M}_{-1} $,
 we expect
that the conformal Ricci flow equations will be of use in  searching for curves of  metrics that asymptotically achieve  the 
$
\sigma
$-constant of
$M$ (see Section~\ref{conclusions}).  Computing $ \sigma (M)$ involves a minimax procedure which attempts to obtain  $
\sigma (M)$ in two steps, first by minimizing  the Yamabe functional in a fixed conformal class  and then by maximizing  over all
conformal classes (Anderson \cite{and97}).   The first step of this procedure, corresponding to the {\it min} part,  is the Yamabe
problem and has been solved.  The unsolved second step then involves a maximization over the conformal classes of $M$.  In the case
that 
$M$ is of 
  Yamabe type $-1$, the conformal classes are represented by  $ {\cal M}_{-1} $ (see
Section~\ref{conformal4}). Thus  the   conformal Ricci flow, which takes place  on the space 
$  {\cal M}_{-1}
$,   exploits the fact that the first part of the problem of realizing the $ \sigma $-constant has been solved and then moves on to
the appropriate smaller space $ {\cal M}_{-1} $  in which to tackle the second part of the problem of finding curves of Riemannian
metrics that asymptotically realize $ \sigma (M)$.

We also remark  that the conformal Ricci flow equation 
was designed as a parabolic model for understanding the reduced Hamiltonian formulation of Einstein's equations of general
relativity     (\cite{fm97},\cite{fm00},\cite{fm01}), similar   to the fact that the
classical Ricci  flow equation is somewhat of  a parabolic model for the unreduced  Hamiltonian formulation of Einstein  equations.  

In comparing the conformal Ricci flow with the classical Ricci flow, one should also make some remarks concerning the  $n=3$ case. 
At first,  the classical Ricci flow  equation had its most spectacular success in the case of strictly positive Ricci tensors in its
attempt to prove the Poincar\'{e} conjecture. Indeed, the classical Ricci flow equation ``prefers" positive curvature (Hamilton
\cite{ham82}, p.\ 276, 279) since it preserves both positive scalar and positive Ricci curvature (see also
Proposition~\ref{increasescale}). The conformal Ricci flow equation, on the other hand, was designed more to deal with the
complementary situation of negative Yamabe type manifolds (see Section~\ref{conformal4}) and hyperbolic geometry rather than
positive Ricci
 tensors and spherical space forms. Thus in this  sense       the classical Ricci flow and the conformal Ricci  flow are concerned
with complementary problems. However, we hasten to add that the classical Ricci flow equation  is now intensely used to study general
questions regarding 3-manifold topology (see for example Cao-Chow~\cite{cc99},  Hamilton (\cite{ham95},\cite{ham99}), Perelman
(\cite{per02},\cite{per03a},\cite{per03b}) and  Ye~\cite{ye93}) and thus the classical and conformal Ricci flows are best viewed
as complementary tools when examining    similar problems.

 \section{The space $ {\cal M}_{-1} $}  \setcounter{equation}{0}
\label{space1}

The space
 \begin{equation}
 {\cal M} _{-1}=\{g \in {\cal M} \mid R(g )=-1\}
\end{equation}
   plays an important role in geometry (Fischer-Marsden \cite{fm75b}), in general relativity
 (Fischer-Moncrief
\cite{fm97}), in  Teichm\"{u}ller theory (Fischer-Tromba \cite{ft84a}), and in the conformal Ricci flow as developed in this
paper. The structure of this  space has been studied in detail in \cite{fm75b}. The main
  viewpoint pioneered  there is to consider the {\bf scalar curvature function}
\begin{equation}
\label{eq:Put equation label here} R: {\cal M}  \longrightarrow  {\cal F} \;,\qquad g \longmapsto R(g)\,,
\end{equation}    as  a map between
infinite-dimensional manifolds and to study the properties of this map. 
Since $ {\cal M} $ is open in $S_2$ and since $ {\cal F} $ is a linear space, the derivative $D\!R(g)$ of $R$ at $g\in {\cal M} $ 
maps
$S_2$ to $ {\cal F} $ and is given by
\begin{equation}\fl\qquad
\label{linearR}D\!R(g ): S_2  \longrightarrow {\cal F} \;,\qquad h \longmapsto D\!R(g )h=
 \Delta_g{\rm
tr}_gh+\delta_g\delta_gh- {\rm Ric}(g) \cdot h\,,
\end{equation}  
 where  $\Delta_g$ is the Laplace operator acting on scalar-valued functions, ${\rm tr}_g$ is the metric trace,
and  $ \delta _g\delta
_g$ is the double covariant divergence acting on 2-covariant symmetric tensors.
In local
coordinates,
\begin{eqnarray} D\!R(g )h & = &  -g^{kl}(g^{ij}h_{ij})_{|k|l}+  g^{ik}g^{jl}h_{ij|k|l} -g^{ik}g^{jl}R_{ij}h_{kl}    
\nonumber      \\   & = & -g^{kl}( h^j{}_{ j})_{|k|l}+   h^{kl}{}_{|k|l} - R^{kl}h_{kl}\nonumber \,,
\end{eqnarray} where the vertical bar denotes covariant derivative with respect to $g$.
We note explicitly that  both the Laplacian  on scalar functions $\Delta_g$  and the covariant divergence $  \delta _g$ are taken
with embedded negative signs,
$  \Delta_g\phi=-g^{ij}\phi_{|i|j}$,  $ (\delta _gh)_j=- g^{ik}h_{ij|k} = 
-h^k{}_j{}_{|k}$, and for a 1-form
$
\alpha
$,
$ 
\delta _g
\alpha   =-g^{ij}\alpha _{i|j}$, so that $ \delta _g \delta _gh= h^{kl}{}_{|k|l}$ and   $ 
\Delta_g\phi = -g^{ij}(\phi_{|i})_{|j}=\delta _g  d\phi$.  Thus in this sign convention,
   the Laplacian  $\Delta_g= \delta _gd$ is a positive operator.

Let
\begin{equation}\fl\qquad
\label{adjoint43}D\!R(g )^*: {\cal F} \longrightarrow S_2\;,\qquad  \phi \longmapsto
D\!R(g )^*\phi=  {\rm Hess}_g  \phi +g
\Delta_g \phi -  {\rm Ric}(g  ) \phi 
\end{equation}  
denote the $L_2$-adjoint of $ D\!R(g )$ where  ${\rm Hess}_g  \phi$ denotes the Hessian of the scalar function $\phi$. In local
coordinates,  
 \begin{equation}
(D\!R(g )^*\phi)_{ij}= \phi_{|i|j}-g_{ij}
(g^{kl} \phi_{|k|l} )-  R_{ij}\phi\,.
\end{equation}  

The second order operator
  $D\!R( g  )^*$ has injective symbol since for
  $\xi\in T_x^*M$ and $1\in \mbox{\boldmath$R$}$,  its symbol   
$ \sigma _\xi(D\!R( g  )^*)1=\xi\otimes\xi-g\|\xi\|^2$   is injective if $\xi\not=0$ since its trace is $(1-n)\|\xi\|^2   $.
Thus using 
  the operators $D\!R( g )$ and $D\!R( g  )^*$  we have the following $L_2$-orthogonal splitting of $S_2$ (see
Berger-Ebin \cite{be69}),
  \begin{equation}  
\label{eqn:splitting1} S_2 =\ker D\!R( g )\oplus {\rm range}\, D\!R(g )^*\;,\qquad h=\bar h+ D\!R(g )^*\phi\,,
\end{equation} 
where $\phi= \big(D\!R(g) D\!R(g)^*\big) ^{-1} \big( 
  D\!R(g)h   
\big)$ and
$\bar h=    h-  D\!R(g)^*\phi  \, $ so that \begin{equation}
 \label{splitnew} 
  \fl  \qquad h=      ( h-  D\!R(g)^*\phi )+  
 D\!R(g)^*\phi=\bar h+ D\!R(g)^*\phi=\bar P_g(h)+ D\!R(g)^*\phi \, ,
\end{equation} 
 where \begin{equation}
\label{projectionortho}  \bar P_g:S_2 \longrightarrow \ker D\!R(g) \;,\qquad h \longmapsto \bar h=\bar P_gh\,,
\end{equation} denotes the  $L_2$-orthogonal   projection   
 onto $\ker D\!R(g)$.

 Note that \begin{equation}
\label{LABEL} D\!R(g) D\!R(g)^* : {\cal F} \longrightarrow {\cal F} 
\end{equation} is an elliptic fourth order operator with  
$
   \ker (D\!R(g) D\!R(g)^*)=\ker   D\!R(g)^* 
$ 
 which   is generically zero as a function of $g$, as shown in Fischer-Marsden \cite{fm75b}.
Also note that in the splitting (\ref{eqn:splitting1}),  there are no curvature assumptions on $g$.

We now consider $D\!R(g)h$ on two types of deformations $h\in S_2$ that we shall need in Section~\ref{maintaining}. 
We first  recall the following. For  any metric  $g\in {\cal M} $, 
 the   {\bf doubly  contracted (differential) Bianchi identity}    asserts that the divergence of the 
  Einstein tensor
   ${\rm Ein}(g)=  {\rm Ric}(g)- {\textstyle\frac{1}{2}}R(g)g
 $  vanishes, 
\begin{equation}
\label{LABEL}0= \delta _g({\rm Ein}(g)) =\delta _g\Big({\rm Ric}(g)   -  {\textstyle\frac{1}{2}} R(g)g\Big) 
=\delta _g({\rm Ric}(g)  ) +  {\textstyle\frac{1}{2}}d R(g)\,,\end{equation} 
where in local coordinates $( \delta _g( R(g)g))_i = -g^{jk}(R(g)g_{ij})_{|k}=-R(g)_{|i}=-(dR(g))_{ i} $.  
Thus the doubly contracted Bianchi identity  is a third-order 
differential identity (i.e., an equation that is  true for all
metrics
$g\in {\cal M}
$) that
``simplifies"       the divergence of the Ricci tensor by expressing it
as the
    derivative of the scalar curvature.
 The {\bf double divergence} of ${\rm Ein}(g)$ then  leads to the fourth order differential identity, 
 \begin{equation}\fl\qquad
\label{eindouble0} 0 =\delta _g\delta _g({\rm Ein}(g)) 
=\delta _g\delta _g({\rm Ric}(g)  ) + {\textstyle\frac{1}{2}}\delta _g d R(g) 
 =\delta _g\delta _g({\rm Ric}(g)  ) +{\textstyle\frac{1}{2}} \Delta_g R(g) \,. 
\end{equation}   

 We consider   examples of the splitting 
(\ref{eqn:splitting1})  for two explicit deformations that we shall need later.
 \begin{example}
\label{ex1}
{\rm
 For $g\in {\cal M} $ and $\varphi\in {\cal F} $, $h=\varphi g\in S_2$ is an {\bf infinitesimal pointwise conformal
deformation of
$g$}. On  a deformation $\varphi g$, 
 \begin{eqnarray} \nonumber 
 D\! R(g)( \varphi g  )  
  & = &  \Delta_g {\rm tr}_g(\varphi g)  + \delta _g \delta _g  (\varphi g) -
{\rm Ric}(g)\cdot (\varphi g)  \\
\fl \nonumber&=& 
   n\Delta_g  \varphi     - \Delta_g  \varphi   -
R(g) \varphi    
\\ \fl  &=& 
 ( n-1)\Delta_g  \varphi      -
R(g) \varphi   \,.
\label{last4d}
\end{eqnarray}
Thus we define an operator $L_g: {\cal F} \rightarrow {\cal F} $ by \begin{equation}
\label{operator0}  L_g \varphi= D\! R(g)( \varphi g  ) \,, 
\end{equation}
so that from (\ref{last4d}), 
\begin{equation}
\label{operator}  L_g    
  : {\cal F} \longrightarrow {\cal F} \;,\quad \varphi \longmapsto L_g\varphi= (n-1)\Delta_g\varphi-R(g)\varphi\,.
\end{equation} 

 We remark that, interestingly, \begin{equation}
\label{LABEL}{\rm tr}_g( D\!R(g )^*\varphi)=   D\!R(g )( \varphi g)= L_g\varphi\,,
\end{equation} since $ {\rm tr}_g( {\rm Hess}_g  \varphi +g
\Delta_g \varphi -  {\rm Ric}(g  ) \varphi)= (n-1)\Delta_g\varphi-R(g)\varphi$.

As a special case of (\ref{operator}), note that for a constant $c\in \mbox{\boldmath$R$} $, \begin{equation}
\label{drg}  L_gc=  D\!R(g)(cg)=-cR(g) \,.
\end{equation} Thus if
$R(g)=-1$,
$L_gc= D\!R(g)(cg)=c$.  Thus   we  have the following $L_2$-orthogonal splitting of $h=g $  under the assumption that $R(g)=-1$,
 in which case $  D\!R(g)g=1$,
\begin{equation} 
\label{splitortho33} g = \bar g + D\!R(g)^*\big( (D\!R(g) D\!R(g)^* ) ^{-1} ( 
1
 ) 
\big) \,, 
\end{equation} where
 $\bar g=\bar P_g(g)= g -  D\!R(g)^*\big( (D\!R(g) D\!R(g)^* ) ^{-1} ( 
 1
 )   
\big)$.
}
\end{example}

\begin{example}
\label{ex2}
{
\rm(See also Examples~\ref{splitRicci1}  and \ref{splittingexample}.)
  For $g\in {\cal M} $,   consider the deformation $h= {\rm Ric}(g) $. 
From (\ref{adjoint43}), \begin{equation}
\label{hessperp}   D\!R(g)^*(-1)= 
  {\rm Hess}_g (-1) +g
\Delta_g(-1)-  {\rm Ric}(g  ) (-1) ={\rm Ric}(g) \,,
\end{equation}  so that $ {\rm Ric}(g) \in {\rm range}\, D\!R(g)^*$.
Thus $\bar P_g ({\rm Ric}(g) )=\overline{ {\rm Ric}(g) }=0$ since  $ {\rm Ric}(g) $ is $L_2$-orthogonal to
$ \ker D\!R(g)$.    Thus $ {\rm Ric}(g) $ does not split (non-trivially)  $L_2$-orthogonally  (see 
Example~\ref{splitRicci1}  for a non-trivial non-orthogonal splitting of $ {\rm Ric}(g) $).

This fact is also ``discoverable"
by integrating (\ref{linearR})  over $M$, to give 
\[\fl\qquad\int_MD\!R(g )h\,d \mu _g=\int_M
( \Delta_g{\rm
tr}_gh+\delta_g\delta_gh- {\rm Ric}(g) \cdot h)\,d \mu _g=
-\int_M
 {\rm Ric}(g) \cdot h\,d \mu _g\,.\]
Thus for all $h\in \ker  D\!R(g) $, 
$\int_M
 {\rm Ric}(g) \cdot h\,d \mu _g=0$, so that $ {\rm Ric}(g) $ is $L_2$-orthogonal to $\ker  D\!R(g)$.

Analytically, if we try to  split  ${\rm Ric}(g) $, we first note that from
  (\ref{eindouble0}),
 \begin{eqnarray} 
  D\! R(g)  {\rm Ric}(g)    
 \nonumber 
  & = &  \Delta_g {\rm tr}_g{\rm Ric}(g) + \delta _g \delta _g  {\rm Ric}(g)-
{\rm Ric}(g)\cdot  {\rm Ric}(g) \\
\fl \nonumber&=& 
  \Delta_g R(g) - {\textstyle\frac{1}{2}}  \Delta_g R(g) -
|{\rm Ric}(g)|^2_g  
\\  &=& 
 {\textstyle\frac{1}{2}}     \Delta_g R(g) -
|{\rm Ric}(g)|^2_g \label{ricci2} \,, 
\end{eqnarray}which is an expression that we shall need frequently.
Thus, from  (\ref{splitnew}) and (\ref{ricci2}), 
\begin{eqnarray} \fl\qquad
{\rm Ric}(g) &=&\overline {{\rm Ric}(g)} + D\!R(g)^*\Big(\big(D\!R(g) D\!R(g)^*\big) ^{-1} \big( 
  D\! R(g)  {\rm Ric}(g) )   
\Big) \nonumber 
\\
\label{splitortho3} 
\fl\qquad
 &=&\overline {{\rm Ric}(g)} + D\!R(g)^*\Big(\big(D\!R(g) D\!R(g)^*\big) ^{-1} \big( 
{\textstyle\frac{1}{2}}
\Delta _g R(g)-  |{\rm Ric}(g) |^2\big)   
\Big)\,.
\end{eqnarray} 
Then, from   (\ref{hessperp})  and (\ref{ricci2}),  \[ 
   D\!R(g) (D\!R(g)^*(-1))= D\!R(g)   {\rm Ric}(g)  ={\textstyle\frac{1}{2}}
\Delta _g R(g)-  |{\rm Ric}(g) |^2\,, 
\] so that 
$(D\!R(g) D\!R(g)^* ) ^{-1} ({\textstyle\frac{1}{2}}
\Delta _g R(g)-  |{\rm Ric}(g) |^2 )=-1$. Thus (\ref{splitortho3})  reduces to 
${\rm Ric}(g) =\overline {{\rm Ric}(g)} + D\!R(g)^*(-1)= \overline {{\rm Ric}(g)} + {\rm Ric}(g) $ so that 
$ \overline {{\rm Ric}(g)} =0$   (see also the remarks after Theorem~\ref{main1}).
}  \hfill$\rule{2mm}{3mm}$
\end{example}
 
\medskip  

The information that we shall need about
$ {\cal M}_{-1} $ is contained in the following theorem. We first recall that 
Aubin \cite{aub76} has shown that every closed $n$-manifold $M$, $n  \ge 3$,  has a Riemannian metric with  constant negative
scalar curvature 
  and thus there are no topological obstructions to the    constraint $R(g)=-1$.
  In particular, $ {\cal M}_{-1}
\ne\emptyset$. 
   \begin{theorem}{\bf(The structure of the space ${\cal M}_{-1} $ (Fischer-Marsden \cite{fm75b}))}
\label{main1}  Let $M$ be a closed connected orientable $n$-manifold with either $n\ge3$ or $n=2$ and ${\rm {\rm genus}}(M)\ge2$. Then 
$ {\cal M} _{-1}$ is a non-empty     
 closed (as a subset)  smooth  infinite-dimensional  {\sc
ilh}-submanifold of
$ {\cal M}
$.
 Moreover, $ {\cal M}_{-1} $  is ${\cal D} $-invariant   and has codimension   $C^\infty (M, \mbox{\boldmath$R$} )$
 in $ {\cal M} $.

  For $ g\in {\cal M} _{-1}$, the tangent space to $ {\cal M}_{-1} $ at $g$  is given by 
\begin{equation}\fl\quad
\label{eqn:splitting2} T_g{{\cal M}_{-1} } \approx  \ker D\!R( g  )  
 =\{\,h\in S_2 \mid  \Delta_g{\rm
tr}_gh+\delta_g\delta_gh- {\rm Ric}(g)\cdot h=0\,\} \,,
\end{equation}
which is the first summand in the 
$L_2$-orthogonal splitting, 
\begin{equation}
\label{eqn:splitting3} T_g{\cal M} =T_g  {\cal M} _{-1}\oplus T_g ^\perp {\cal M} _{-1} \,.
\end{equation} 
Equivalently, 
\begin{equation}
\label{eqn:splitting4} S_2 =\ker D\!R(g  )\oplus {\rm range}\, D\!R(g )^* \,,
\end{equation} 
where 
   $T_g ^\perp {\cal M} _{-1} \approx {\rm range}\, D\!R(g )^*  $ 
 is the  $L_2$-orthogonal complement of $T_g{\cal M} _{-1}$ in the natural $L_2$-Riemannian metric
 $ {\cal G} $ (see $(\ref{naturalmetric})$)  on $ {\cal M} $.

If $M$ is of Yamabe type $-1$, then the space of pointwise conformal classes  $ {\cal M} / {\cal P} $  on $M$ and $ {\cal M} _{-1}$
are {\sc ilh} diffeomorphic.

 \end{theorem} 
 {\bf Proof (sketch):}
 It is   shown in \cite{fm75b} that the scalar curvature  map  $R: {\cal M} \rightarrow {\cal F} $  is  a smooth {\sc ilh}
  submersion
almost everywhere. In  particular, $-1$ is a regular value of $R$ and thus  from the inverse function theorem adapted to  the {\sc
ilh} topology,  the level set  
\begin{equation}
\label{LABEL}  {\cal M}_{-1}  =R^{-1}(-1)
\end{equation} is a  smooth closed {\sc ilh} submanifold of $ {\cal M} $ with the indicated tangent space.

If   
$g\in  {\cal M} 
$ and
$f\in {\cal D}
$,  let $f^*g$ denote the pullback of $g$ by $f$, defined 
 by 
  $f^*g(x)(v_1,v_2)=g(f(x))(T_xf(v_1), T_xf(v_2))$  
for $x\in M$, 
$v_1, v_2\in T_xM$, the tangent space to $M$ at $x\in M$.

By the {\bf covariance of the scalar curvature operator}
\begin{equation}
\label{coveqn1}  R(f^*g)=f^*(R(g))=R(g)\circ f\,.
\end{equation}   The $   {\cal D}
$-invariance of
$ {\cal M}_{-1} $ then follows easily, 
for if $g\in  {\cal M}_{-1} $, \[
\label{LABEL}
 R(f^*g)=f^*(R(g))=f^*(-1)=-1\,,
\]so that $f^*g\in {\cal M}_{-1} $.

 If $M$ is  of Yamabe type $-1$, then every
Riemannian metric on
$M$ is uniquely pointwise conformally deformable to a metric with scalar curvature
$=-1$. An analysis of the proof shows that this bijection is an {\sc ilh} diffeomorphism, so 
 the 
space of conformal structures  $ {\cal M} / {\cal P} $ is {\sc ilh} diffeomorphic to $ {\cal M} _{-1}$.
 \hfill$\rule{2mm}{3mm}$ \\

Thus 
in the case of Yamabe
type
$-1$ manifolds, the
  space $ {\cal M} _{-1}$  is an important  representation of $ {\cal M} / {\cal P} $    and thus plays an important role
in conformal geometry.

As discussed in Example~(ii) above, for any  $g\in {\cal M} $,   $ {\rm Ric}(g) \in {\rm
range}\,D\!R(g)^* =(\ker D\!R(g) )^\perp 
$ (in the $L_2$-metric). Thus we have the following interesting geometrical fact about the   Ricci tensor:
 If $\, {\rm Ric}(g) \not=0$, then in an $L_2$-sense, the direction 
    ${\rm Ric}(g) \in  (\ker D\!R(g) )^\perp $     infinitesimally deforms   $g$   in a direction orthogonal to 
those directions preserving the scalar curvature of $g$.

This can be formalized as follows.
 Let \[{\cal M} _ \rho =R ^{-1} (\rho  )=\{\,g'\in
{\cal M}  \mid R(g')= \rho \,\}\]denote the space of Riemannian metrics with scalar curvature $ \rho $. Then, generically, $ {\cal
M} _
\rho
$  is a    $C^\infty (M, \mbox{\boldmath$R$} )$-codimensional submanifold of
$ {\cal M}  $ (Fischer-Marsden \cite{fm75b}) and $T_g {\cal M} = T_g {\cal M} _ \rho  \oplus  T_g ^\perp{\cal M} _\rho $, in which
  case     $ {\rm Ric}(g) \in {\rm range}\,D\!R(g)^* \approx T_g ^\perp{\cal M}  _   \rho $  and so  is $L_2$-orthogonal to $
{\cal M} _ \rho $  (see also Figure~1 in Section~\ref{nonortho}).

Somewhat similarly, for $g\in {\cal M} $ with $ {\rm vol}(M,g) =c>0$, let 
\[ {\cal M} ^c=\{\,g'\in {\cal M} \mid {\rm vol}(M,g') =c \,\}\] denote the 
  codimension-1
submanifold of metrics with volume $c$. Then 
\[
\label{m1} T_g {\cal M}
^c\approx\{\,h\in S_2\mid \int_M {\rm tr}_gh\,d \mu _g= \int_M g\cdot h\,d \mu _g=0\,\} \,,
\] so that
\[
\fl\quad
\label{LABEL} T_g^\perp {\cal M}^c \approx \{\,h'\in S_2\mid \int_Mh'\cdot h\, d \mu _g=0\;\mbox{for all}\;h\in T_g {\cal M}
^c\,\}  =
\{\, \lambda g\mid \lambda \in
\mbox{\boldmath$R$}\,\}=\mbox{\boldmath$R$} g\,.
\]
Thus, in an $L_2$-sense, any   metric $g\in {\cal M} $, if viewed as a    ``vector" $g\in \mbox{\boldmath$R$} g  \approx
T_g^\perp {\cal M}^c $,   is the  direction 
 which   maximally  increases the volume $ {\rm vol}(M,g)  $.

  An interesting and useful kinematical observation regarding  $ {\cal M}_{-1} $ is the following. 
 \begin{proposition} 
\label{kinematical2} 
For      $g \in {\cal M}  $,   \begin{eqnarray} \label{ricidentity} |{\rm Ric}( g )|^2_{ g } & = &  |{\rm Ric}^T\! ( g )|^2_{ g
}+{\textstyle
\frac{1}{n}}      R^2(g)      \\  |{\rm Riem}( g )|^2_{ g }  & = &  |{\rm Weyl}( g )|^2_{ g }+ {\textstyle\frac{4}{n-2}}
 |{\rm Ric} ( g )|^2_{ g }-{\textstyle
\frac{2}{ (n-1)(n-2)}}      R^2(g)     \\ & = &
  |{\rm Weyl}( g )|^2_{ g }+ {\textstyle\frac{4}{n-2}}
 |{\rm Ric}^T\!( g )|^2_{ g }+{\textstyle
\frac{2}{n(n-1)}}      R^2(g)     
\end{eqnarray}(where $R^2(g)=(R(g))^2$).  Consequently, for $g \in {\cal M}_{-1} $, 
 \begin{eqnarray} \label{ricidentity1} |{\rm Ric}( g )|^2_{ g } & = &  |{\rm Ric}^T\! ( g )|^2_{ g }+{\textstyle
\frac{1}{n}}    \ge  {\textstyle
\frac{1}{n}}          \\ |{\rm Riem}( g )|^2_{ g }  & = &  |{\rm Weyl}( g )|^2_{ g }+ {\textstyle\frac{4}{n-2}} |{\rm Ric}^T\! ( g
)|^2_{ g }+{\textstyle
\frac{2}{n(n-1)}}       \ge {\textstyle
\frac{2}{n(n-1)}}   \,.    
\end{eqnarray}
 Thus,    kinematically,  any    set (or curve) of   metrics in $ {\cal M} _{-1} $  
  has   Ricci and Riemann curvature norm bounded away from zero.  Thus
the constraint
$R(g)=-1$  prevents the norms $|{\rm Ric}( g )|_g$ and  $|{\rm Riem}( g )|_g$
  from collapsing to zero.
 \end{proposition} {\bf Proof:}  For $ g\in {\cal M}  $,  
   $ {\rm Ric}  ( g  )  = {\rm Ric}^T\!( g )
+{\textstyle
\frac{1}{n}}R(g) g $. Squaring and using the facts  that 
${\rm Ric}^T\!( g )$ and $  g$  are pointwise orthogonal  and that $|  g  |^2_{ g }=
g^{ik}g^{jl}g_{ij}g_{kl}=g^{kl} g_{kl}={\rm tr}_gg=n$ (not $n^2$) gives
\[  \fl\quad
  |{\rm Ric}( g )|^2_{ g }=
 |{\rm Ric}^T\!( g )+{\textstyle \frac{1}{n}}R(g) g |^2_{ g }= |{\rm Ric}^T\!( g )|^2_{ g }+{\textstyle
\frac{1}{n^2}}R^2(g)| g |^2_{ g }
=|{\rm Ric}^T\!( g )|^2_{ g }+{\textstyle \frac{1}{n}}R^2(g) 
\]

In a similar (but longer) manner,  the expressions for $ |{\rm Riem}( g )|^2_{ g }$ follows by squaring the local coordinate
expression 
 which   defines the Weyl conformal curvature tensor,\begin{eqnarray}R_{ijkl} & = &W_{ijkl}+ {\textstyle\frac{1}{n-2}}
(g_{jk}R_{il}+ g_{il}R_{jk}-g_{jl}R_{ik}-g_{ik}R_{jl})  
\nonumber     \\  &   & \qquad{}+{\textstyle\frac{1}{(n-1)(n-2)}}R(g) ( g_{ik}g_{jl}-g_{jk}g_{il} ) \,.   \nonumber 
\end{eqnarray} 
 \hfill  \rule{2mm}{3mm}

 \section{The conformal Ricci flow equations}   \setcounter{equation}{0}

\subsection{The equilibrium points of the conformal Ricci flow equations}

\begin{proposition}{\bf(The equilibrium points of the conformal Ricci flow equations)} 
\label{equilprop}
  A  metric $g \in {\cal M}_{-1} $ is an equilibrium point
of the conformal Ricci flow equations  if and only if
 $ g$ is an Einstein  metric with  
negative Einstein constant $-\frac{1}{n}$, 
 \begin{equation}
\label{equilpts} {\rm Ric}(g)= -{\textstyle\frac{1}{n}}  g \,,
\end{equation}
if and only if $ {\rm Ric}^T\!(g) =0$,
in which case  the conformal pressure  \begin{equation}
\label{LABEL} p=0\,.
\end{equation}

 If $n=3$, then $g\in {\cal M}_{-1} $ is an equilibrium point if and only if $ g$ is a
hyperbolic metric with constant negative sectional curvature $K(g)=-\frac{1}{6} $ and  ${\rm Ric}(g)=-\frac{1}{3}g$,  
 in which case
$g$ is unique up  to isometry.
\label{prop:Put label here}  
\end{proposition}
  {\bf Proof: }Suppose $ g\in {\cal M}_{-1} $ is an    equilibrium point  of the conformal Ricci flow equations with corresponding
 conformal pressure $p=p(g)\in {\cal F} $. Then 
\begin{equation}
\label{equil}  2 {\rm Ric}(g)+ ( {\textstyle\frac{2}{n}}+p)  g  =0 
\end{equation} and the trace yields $ 2R(g)+2+ np= -2+2+np=np=0  $.
 Substituting $p=0$ 
back into (\ref{equil}) gives \begin{equation}
\label{LABEL}  {\rm Ric}(g)= -{\textstyle \frac{1}{n}}\, g \,.
\end{equation} 

 Conversely, if   ${\rm Ric}(g)
=-{\textstyle
\frac{1}{n}}\, g
$, 
from (\ref{global}) below,
 $p=p(g)=0 $ so that  $ 2 {\rm Ric}(g)+ ( {\textstyle\frac{2}{n}}+0)  g   =0$ so that $g$ is an equilibrium point (see also
Remark~(ii) following Proposition~\ref{maincons}).
  For $g\in {\cal M}_{-1} $, $ {\rm Ric}^T\!(g) = {\rm Ric}(g) -\frac{1}{n}R(g)g= {\rm Ric}(g) +\frac{1}{n} g$ so that $g$ is an
equilibrium point if and only if $ {\rm Ric}^T\!(g)=0$.

 If $n=3$, then $g\in {\cal M}_{-1} $ is an equilibrium point if and only if ${\rm Ric}(g)= 
 -{\textstyle \frac{1}{3}}\, g$ if and only if  $ g$ is a hyperbolic metric with constant   sectional curvature
$K(g)=-\frac{1}{6}
\,$, in which case,  by Mostow rigidity, $g$ is unique up  to isometry 
 when restricted to $ {\cal M}_{-1} $  (see definition of hyperbolic metric in  Section~\ref{crs}).
\hfill$\rule{2mm}{3mm}$ \\

\medskip 
\noindent {\bf Remarks: }
\begin{enumerate}\item 
If $M= \mbox{\boldmath$S$}
^1$, $ \mbox{\boldmath$S$} ^2$, or $ \mbox{\boldmath$T$} ^2$,  $ {\cal M}_{-1} $ is empty  for   these
manifolds and so the conformal Ricci flow equations (\ref{parabolic}--\ref{elliptic})   have no solutions. If    
$n=2 $, {\rm genus}$\,(M)\ge2$,  
  every  $g\in {\cal M}_{-1}\not=\emptyset $ satisfies  ${\rm Ric}(g) =-\frac {1}{2}\, g$ and thus every $g\in {\cal M}_{-1}$
 is an equilibrium point  of (\ref{parabolic}--\ref{elliptic})  with  
$p=0 $.
Thus in this case the
evolution equation
$(\ref{parabolic})$ reduces to the trivial equation
$\frac{\partial g}{\partial t}\equiv0 $.
Thus  we    restrict to 
   $n\ge3$.  
\item
At an equilibrium point $g\in {\cal M}_{-1} $, the conformal pressure is zero   since in that case trivially no constraint force is
necessary to constrain the flow to $ {\cal M}_{-1} $.
  In Proposition~\ref{maxprin} below  we will see that,  conversely, if at any time $t_0$ the conformal pressure $ p(t_0, x_0)=0$
for some     $x_0\in M$,   then 
$p $  is identically zero and 
   $g$
must be an equilibrium point,
$  {\rm Ric}(g)=-  {\textstyle\frac{1}{n}} g
$.

\item A   3-manifold $M$
has an equilibrium point $g\in {\cal M}_{-1} $ if and only if $M$ supports a hyperbolic metric. Contrapositively,  if $M$ does not
support a hyperbolic metric,  then there are no equilibrium points for the conformal Ricci flow. For   example, if $M$ is
diffeomorphic to any spherical space form, flat manifold, handle   $
\mbox{\boldmath$S$} ^1\times
\mbox{\boldmath$S$} ^2$,  or to any non-trivial  connected sum
$M_1\#M_2$ of 3-manifolds ($M_1\not =
\mbox{\boldmath$S$} ^3$, $M_2\not =
\mbox{\boldmath$S$} ^3$), then $M$ would  not support a hyperbolic metric and thus there would be no equilibrium points for the
conformal Ricci flow on $M$.
\item
We   note the  significant differences  between the cases $n=3$ and $n\ge4$. If $n=3$, equilibrium points
are hyperbolic metrics and thus, since they are normalized by $R(g)=-1$,  are  by Mostow rigidity   unique up to isometry. On the
other hand, for 
$n\ge4$, there may exist continuous families of {\it non-isometric} Einstein metrics with 
${\rm Ric}(g)=-\frac{1}{n}g$ which would be equilibrium points of the conformal Ricci flow. A similar situation also arises for the
higher dimensional reduced Einstein equations (see Fischer-Moncrief \cite{fm02b} for more
information regarding this situation). 
 \item 
 If we let $g_e$ denote an equilibrium point, where the ``$e$" stands for  equilibrium (or einstein), then the static
(or equilibrium, or trivial) curve
$g(t)\equiv g_e$ is a solution of the conformal Ricci system  with initial value $g_e$ and $p=0$.  On the other
hand, if
$g(t)$ is a non-static (or non-equilibrium, or non-trivial) solution, then no equilibrium point can be on the image of
$g(t)$. Thus for a non-static solution   $g$, ${\rm Ric}^T\!(g(t))\not =0$ for any $t\in [0,T)$.
\item We remark that  the equilibrium points for both the  classical    and conformal Ricci flow  equations are Einstein
metrics.
\end{enumerate}

\subsection{Maintaining the constraint}
\label{maintaining}
The   conformal  pressure $p$    can be thought of as a time-dependent
Lagrange multiplier necessary to maintain the constraint equation $R(g)=-1$.  Here we find the elliptic  equation satisfied by
$p$.
In Section~\ref{elimin}  we  show how $p$   
and  the constraint equation itself  can be eliminated.

 \begin{proposition}[Maintaining the constraint equation]
\label{maincons}  Let    $g:[0,T)\rightarrow {\cal M}_{-1}  $ be a   solution of the
conformal Ricci flow equations with conformal pressure
  $p:[0,T) \rightarrow {\cal F} $.
 Then  for each $t\in [0,T)$, 
\begin{equation}
\label{LABEL} \frac{\partial g}{\partial t}\in \ker D\!R(g)\,,
\end{equation} 
 and
 $ p $   satisfies the linear inhomogeneous elliptic  equation \begin{equation}
\label{elliptic6}  L_g  p   =  ( n-1)\Delta_g  p    +
p = 
   2
(|{\rm Ric} (g)|^2_g -  {\textstyle\frac{1}{n}}) = 2
(|{\rm Ric}^T \!(g)|^2_g  ) \,,
\end{equation}
where the operator $L_g     =  ( n-1)\Delta_g       +
1 : {\cal F} \rightarrow {\cal F} $ is an isomorphism and is the operator $L_g$ defined in $(\ref{operator0}$-- $\!\ref{operator})$ 
when
$g\in {\cal M}_{-1} $.
Thus for each
$t\in[0,T)
$,   the  conformal pressure $p$ is uniquely determined by \begin{equation}
\label{global} \fl\qquad p= 2 L_g ^{-1} \Big(
|{\rm Ric} (g)|^2_g -  {\textstyle\frac{1}{n}}\Big)
=2 L_g ^{-1} \Big(
|{\rm Ric}^T \!(g)|^2_g \Big)= 2 L_g ^{-1} \Big(
|{\rm Ric} (g)|^2_g \Big)-  {\textstyle\frac{2}{n}}\,.
\end{equation}  
   \end{proposition}
  {\bf Proof:} Differentiating the  constraint equation $R(g)=-1$ with respect to time
  yields 
\begin{eqnarray} 
\label{test}  \frac{  \partial R(g)}{\partial t}  = D\! R(g)  \frac{\partial g}{\partial t}  =0\,,
 \end{eqnarray}so that $ \frac{\partial g}{\partial t}\in \ker D\!R(g)$.
 Then from    (\ref{parabolic}), (\ref{last4d}), and
  (\ref{ricci2}), 
 \begin{eqnarray}\fl 
0 & = &  
     D\! R(g)  \frac{\partial g}{\partial t} =D\! R(g)   \left(-2\,{\rm Ric}(g)- ( {\textstyle\frac{2}{n}}+p)   
g \right )=-2D\! R(g)    {\rm Ric}(g)-D\! R(g)( ( {\textstyle\frac{2}{n}}+p) g  ) \nonumber \\
\fl \nonumber 
  & = & -2 \Big( \Delta_g {\rm tr}_g{\rm Ric}(g) + \delta _g \delta _g  {\rm Ric}(g)-
{\rm Ric}(g)\cdot  {\rm Ric}(g)\Big) \\
\fl
\label{example}&&
\qquad\qquad\qquad{}-\Big( \Delta_g {\rm tr}_g((p+{\textstyle\frac{2}{n}})  g)  + \delta _g \delta _g 
(( {\textstyle\frac{2}{n}}+p)  g) - {\rm Ric}(g)\cdot (( {\textstyle\frac{2}{n}}+p)  g) \Big)\nonumber  \\
\fl \nonumber&=& 
 -2 \Big( {\textstyle\frac{1}{2}} \Delta_g R(g)  -
|{\rm Ric}(g)|^2_g\Big) 
-\Big( (n-1)\Delta_g  ( {\textstyle\frac{2}{n}}+p)       -
R(g) ( {\textstyle\frac{2}{n}}+p)  \Big)  
\\   \fl  &=& 
 -    \Delta_g R(g) +2
|{\rm Ric}(g)|^2_g 
- ( n-1)\Delta_g p + (p+{\textstyle\frac{2}{n}}) R(g) \,.
   \label{last4}
\end{eqnarray} 
In  local
coordinates,  (\ref{test})  and (\ref{last4})  are given by \[0=  \frac{\partial R }{\partial t}= g^{ij}    R _{|i|j} +2
R^{ij}R_{ij}    
+ ( n-1)g^{ij}  p   _{|i|j}    +(p+{\textstyle\frac{2}{n}}) R \,.
    \] 
  Applying   
  the  constraint   $R(g)=-1$,    (\ref{last4}) reduces  to  
\[
\label{LABEL}  0=  2
|{\rm Ric}(g)|^2_g 
- ( n-1)\Delta_g p - (p+{\textstyle\frac{2}{n}})  =
 2
|{\rm Ric}(g)|^2_g 
  -L_g  p -{\textstyle\frac{2}{n}}  \,,
\]so that using the equality part of  (\ref{ricidentity1}),
\begin{equation}
\label{maintainConstraint} L_g  p = 2(
|{\rm Ric}(g)|^2_g 
-{\textstyle\frac{1}{n}})   =
 2 
|{\rm Ric}^T\!(g)|^2_g \,,
\end{equation} where  
$L_g   
 = (n-1)\Delta_g  +1=L_g^*$ 
is a  linear self-adjoint elliptic operator.
Since $\ker L_g=\ker L^*_g= 0$, $L_g =L^*_g$ is  injective and thus  by ellipticity of $L_g$,
 \[{\cal F} ={\rm range}\, L_g\oplus \ker L_g^*={\rm range}\, L_g\,.\]Thus $L_g$ is surjective and thus an
     isomorphism. 
Let 
  $ L_g ^{-1} 
  : {\cal F} \rightarrow {\cal F}$ denote its inverse.
Then $ p=  2 L_g ^{-1}  (
|{\rm Ric}(g)|^2_g 
-{\textstyle\frac{1}{n}} )= 2 L_g ^{-1}  (
|{\rm Ric}^T\!(g)|^2_g  )$ uniquely solves (\ref{maintainConstraint}). Lastly, since for any constant $c$, $L_gc=  (n-1)\Delta_g c
+c  =c$, $L^{-1} _g(c
  )= c $, so that  $p= 2 L_g ^{-1}  (
|{\rm Ric}(g)|^2_g 
-{\textstyle\frac{1}{n}} )= 2 L_g ^{-1}  (
|{\rm Ric}(g)|^2_g )
-{\textstyle\frac{2}{n}}  $\,\,.
 \hfill$\rule{2mm}{3mm}$ \\

\noindent {\bf Remarks:}
\begin{enumerate}
\item 
   Equation (\ref{global}), $p(t)=  2 L_{g(t)} ^{-1}  (
|{\rm Ric}^T\!({g(t)})|^2_{g(t)}  )$, links the conformal pressure   $p(t)$ to the flow  ${g(t)}$  after imposing the scalar
curvature constraint
$R(g(t))=-1$. This equation should be compared with 
$c(t)= 
  {{\textstyle\frac{2}{n}}  }    R_{\rm total} ({g(t)}) $ for the classical Ricci flow
which links 
$c(t)$ to $g(t)$, also {\it after}  imposing the volume constraint $ {\rm vol}(M,g) =1$ (see the  remark above (\ref{usual0})).  
\item
If  $g$ is an equilibrium point, then from Proposition~\ref{equilprop},  ${\rm Ric}^T
\!(g)=0$  and  thus
$   p  =   2 L_{g } ^{-1} (0)=0
$ as also  known   from that Proposition.  Conversely, from (\ref{global}), if $ {\rm Ric}^T\!(g)=0$, then $p=0$ and so $g$ is an
equilibrium point. 
\item  
The conformal pressure   $p= 
 2 L_g ^{-1}  (
|{\rm Ric}^T\!(g)|^2_g  )$ is   an 
integral or nonlocal function of $g$ in contrast to a local, or differential function of $g$, such as $R(g)$.
Thus 
  at a point  $x\in M$, $p(x)$
depends globally on  
$g$ 
whereas $R(g)(x)$ depends  on
$g$  only locally in a neighborhood $U_x\subset M$ of $x$.
  This global behavior of $p$  is completely analogous to the pressure that appears  in the incompressible Euler or Navier-Stokes
equations and is also 
      analogous to the nonlocal  Lagrange multiplier $c $ in the classical Ricci equation (\ref{parabolic22}).
 Similarly,   
the nonlocal equation (\ref{global})   that determines $p$ is completely analogous to the nonlocal equation
(\ref{eliminatec})
that determines $c$.
 \hfill$\rule{2mm}{3mm}$ 
 \end{enumerate}
\medskip

Using the strong maximum principle, we can find  additional important and useful information about  the conformal pressure 
$p$.

By a {\bf domain} $D\subseteq M$ we mean an open connected set of $M$.

\begin{proposition}{\bf(The strong maximum principle and the conformal pressure)}
\label{maxprin}  For    $g\in{\cal M}_{-1}  $, let 
  $p\in {\cal F} $ be a solution of 
\begin{equation}
\label{elliptic66} L_gp= ( n-1)\Delta_g  p      +
p  = 2
|{\rm Ric}^T\!(g)|^2_g  \,.
\end{equation}Then  $p \ge 0$, and more particularly,  either
\begin{description}\item[{$(a)$}]
$p=0 $ and ${\rm Ric}(g)
=-{\textstyle\frac{1}{n}}g$, or 
 {$(b)$} $
  p>0 $
 and $ {\rm Ric}(g)
\not=-{\textstyle\frac{1}{n}}g
$.
\end{description}
 \end{proposition} {\bf Proof:} Equation (\ref{elliptic66})   
  is in a form in which we can apply the {\bf strong maximum principle} (also known as the maximum principle of
H.~Hopf;   see e.g.\ Protter-Weinberger
\cite{pw84}, pp.\ 61, 64) to $-p$. 
The corresponding
strong minimum  principle    asserts   in our notation  that if 
$ \phi \in {\cal F} $ satisfies the inequality \begin{equation}
\label{LABEL}  - \Delta_g  \phi      -u
 \phi   \le0\;,\quad\qquad u\ge0\,,
\end{equation} on any domain $D\subseteq M$ 
and if   $\phi_{\rm min}$ denotes the minimum value of $\phi $ on $D$, then if  $\phi_{\rm min}\le0$, then 
 $\phi\equiv \phi_{\rm min}$ is a constant function on $D$.
Applied to the  negative of  (\ref{elliptic66}) with $u=1$,  
 \begin{equation}
\label{elliptic9}  -( n-1)\Delta_g  p     -
p = -2
 |{\rm Ric}^T\!(g)|^2_g \le0\,,
\end{equation} and
letting $p_{\rm min}$ denote  the minimum value of $p $ on $M$,      
if  $p_{\rm min}\le0$, then
$p=p_{\rm min} $ is
a constant function.
  But if $p $ is   constant, then from
  (\ref{elliptic9}),
 \begin{equation}
\label{elliptic109}    
p=    2
 |{\rm Ric}^T\!(g)|^2_g \ge0\,,
\end{equation} 
and since 
$p=p_{\rm min}\le0$, it follows that $p=p_{\rm min}=0$. Thus from   
(\ref{elliptic109}),   
\begin{equation}\fl\quad
\label{LABEL} 0=    2
 |{\rm Ric}^T\!(g)|^2_g    \Longrightarrow  {\rm Ric}^T\!(g) =0\Longrightarrow  {\rm Ric} (g) = {\textstyle\frac{1}{n}}
R(g)g\Longrightarrow
 {\rm Ric} (g) =-{\textstyle\frac{1}{n}}g\,.
\end{equation}

If on the other hand  $p_{\rm min}>0$, then  $p\ge p_{\rm min}>0$ and thus from (\ref{global}), 
${\rm Ric}(g)\not =-{\textstyle\frac{1}{n}}g$.
  Thus overall, either
(1) $p=0$ and  ${\rm Ric}(g) =-{\textstyle\frac{1}{n}}g$,  or (2) $p>0$ 
and ${\rm Ric}(g)\not =-{\textstyle\frac{1}{n}}g$.
\hfill$\rule{2mm}{3mm}$ 

\bigskip
\noindent
{\bf Remark:}  As a consequence of this proposition, if   $g:[0,T) \rightarrow {\cal M}_{-1} $ 
is a non-static  conformal Ricci flow,  then its conformal      pressure
$p:[0,T) \rightarrow {\cal P}  $ must  lie in the space of positive functions $ {\cal P} = C^\infty(M,
\mbox{\boldmath$R$} ^+)$. To put this another way, if $p(t_0,x_0)=0$  vanishes at a ``spacetime"  point    $(t_0,x_0)\in
[0,T)\times M$, 
 then  $p(t ,x )=0$  for any $(t,x)\in[0,T)\times M$.

\subsection{Elimination of  the constraint and the reduced conformal Ricci flow}
\label{elimin}Since the conformal Ricci flow equations are    a  constrained dynamical system with the conformal pressure $p$ acting
as  a Lagrange multiplier,  both   $p$ and the constraint equation $R(g)=-1$   
    can  be formally eliminated to give the   {\bf (fully) reduced conformal Ricci equation}.  This is analogous to the
elimination of   $c$ to give the classical Ricci flow equation as discussed in Section~\ref{comparison}  and the
elimination of the fluid  pressure $p$ in the incompressible Euler or Navier-Stokes equations (see Ebin-Marsden \cite{em70}).

 \begin{proposition}{\bf(Elimination of the constraint   and the reduced conformal Ricci flow equation)}
\label{elimination1}  
A smooth curve 
 $ 
 g:[0,T) \rightarrow {\cal M} 
 $ with initial value $g_0\in {\cal M}_{-1} $ is a solution of the  
 {\bf (fully) reduced conformal Ricci flow equation}

\begin{equation}
\label{reducednew}   \frac{\partial g}{\partial t} +2\,{\rm Ric}(g)=-  
      2\Big( L ^{-1} _g( |{\rm Ric} (g)|_g^2 )\Big) g\,,
\end{equation} 
 where 
$ L_g 
 = (n-1)\Delta_g +1 : {\cal F} \rightarrow {\cal F}
$
is a  curve of   invertible   elliptic operators   with inverses
$  L_{g }^{-1}: {\cal F} \rightarrow {\cal F} $  associated with the curve $g:[0,T) \rightarrow {\cal M} $, 
if and only if $g$ is a solution of the  conformal Ricci flow equations 
 \begin{eqnarray}\label{parabolicnew} \frac{\partial g}{\partial t}+2\big({\rm
Ric}(g){\textstyle+\frac{1}{n}}g\big)=-p   g          \\
\label{ellipticnew} \;\;\;\qquad\qquad\quad  R(g)=-1   
\end{eqnarray}   
with  conformal pressure $p=  2    L ^{-1} _g( |{\rm
Ric} (g)|_g^2  )- {\textstyle\frac{2}{n}}  
$.  
  \end{proposition}
 \noindent
  {\bf Proof:} Proposition~\ref{maincons} asserts that if $  g  $ is a conformal Ricci flow with conformal pressure $p$, then 
  $   p =2    L ^{-1} _g( |{\rm
Ric} (g)|_g^2  )- {\textstyle\frac{2}{n}}  
$  with $L_g=   (n-1)\Delta_g +1$. Substituting   $p$ into (\ref{parabolicnew}) shows that  a   conformal Ricci flow $g$  is also a 
solution to  the reduced conformal Ricci equation (\ref{reducednew}). 

To show the converse, let $g$ be a solution to the reduced conformal Ricci flow equation (\ref{reducednew}) with initial value
$g(0)=g_0\in {\cal M}_{-1} $. Set
   $ p =2    L ^{-1} _g( |{\rm
Ric} (g)|_g^2  )- {\textstyle\frac{2}{n}}  
$  so that $g$ and   $p$   satisfy (\ref{parabolicnew}).
The task   is  to show that the constraint $R(g)=-1$  is maintained.
From the definition of $p$ and noting that  $L_g( {\textstyle\frac{2}{n}}  )= {\textstyle\frac{2}{n}}  $, we see that  $p$ satisfies
\begin{equation}
\label{stuff422}    L_gp=   ( n-1)\Delta_g  p +
p  =
   2(
|{\rm Ric} (g)|^2_g  - {\textstyle\frac{1}{n}} )\,.
\end{equation} 
 Let $\tilde R(g)=1+R(g)$.  
Then using the expression for $  D\!R(g)\frac{ 
\partial
 g }{\partial t}$ from (\ref{last4}) 
 (which is valid for any solution of (\ref{parabolicnew})  whether or not the constraint equation  
 $R(g)=-1$ is satisfied) and using  (\ref{stuff422}),
 \begin{eqnarray} \fl \qquad\qquad  \frac{  \partial\tilde R(g)}{\partial t}&= &\frac{  \partial R(g)}{\partial t}=D\!R(g)\frac{ 
\partial
 g }{\partial t} =D\! R(g)   \left(-2 {\rm Ric}(g)- ( {\textstyle\frac{2}{n}}+p)   
g \right )  \nonumber \\  &=& 
 -    \Delta_g R(g) +2
|{\rm Ric}(g)|^2_g 
- ( n-1)\Delta_g p + (p+{\textstyle\frac{2}{n}}) R(g)  \nonumber \\ 
  \nonumber  &=& 
  -    \Delta_g  (1+ R(g))  +\left(2|{\rm Ric} (g)|^2_g  - 
( n-1)\Delta_g p\right ) +    ( p +{\textstyle\frac{2}{n}})R(g)
\\ &=&  \nonumber -    \Delta_g\tilde R(g)    + (p+ {\textstyle\frac{2}{n}} )  +  (p+ {\textstyle\frac{2}{n}} )  R(g)    
\\
 &=&-    \Delta_g\tilde R(g)   +(  p+{\textstyle\frac{2}{n}}  )  \tilde
R(g)    \,. \label{rtt}
  \end{eqnarray} Since  $g=g(t)$ is given as a solution to the reduced conformal Ricci system,  
$ p(t)  =  2  L ^{-1} _{g(t)}(|{\rm Ric}(g(t))|_{g(t)}^2)- {\textstyle\frac{2}{n}}  $ is  a known function of $t$.  Thus  
the function
$\tilde R(t)=\tilde R(g(t))$  
 satisfies a linear   heat  equation on
$M$,
\begin{equation}
\label{kac} \fl\qquad\quad  \frac{  \partial\tilde R(t)}{\partial t}  =-    \Delta_{g(t)}\tilde R(t)   +\left
(p(t)+{\textstyle\frac{2}{n}}\right) 
\tilde R(t)  \,, \quad \tilde R(0)=\tilde R(g_0)=1+R(g_0)\,,
\end{equation}with time-dependent Laplacian and time-dependent coefficients.
Uniqueness of  solutions to  (\ref{kac})  
  follows  by modifying standard arguments to take into account the time-dependent Laplacian.   
 For example,  although it is more than is needed, 
uniqueness of solutions to 
  quasi-linear parabolic systems  with time-dependent coefficients 
is shown in  Taylor \cite{tay96c}, Section~15.7 (see Proposition~7.3,  p.\ 332).

Since $R(g_0)=-1$,    $\tilde R(0)=\tilde R(g_0)=1+R(g_0)=0$    and thus
$\tilde R(t)=\tilde R(g(t))\equiv0$ on $[0,T)$. Thus the constraint $ R(g(t))=-1$ is maintained on  $[0,T)$ and thus  the
reduced unconstrained equation (\ref{reducednew}) with initial values in $ {\cal M}_{-1} $ is equivalent to the  unreduced
constrained system  (\ref{parabolicnew}--\ref{ellipticnew}).
  \hfill$\rule{2mm}{3mm}$ \\ 

\noindent {\bf Remarks: }
\begin{enumerate}
\item
 The  unreduced conformal Ricci   flow equations  (\ref{parabolicnew}--\ref{ellipticnew})  require  the   constraint
$R(g )=-1$   on  the entire
flow whereas the  reduced conformal Ricci flow equation (\ref{reducednew}) requires  the constraint  $R(g_0)=-1$  
  only on the initial metric.  The constraint is then automatically  maintained on the entire flow as   a consequence of the
reduced evolution equation.
  This situation is exactly analogous to the relation of the unreduced classical Ricci flow equations 
 (\ref{parabolic22}--\ref{elliptic22}) 
  to the fully reduced
classical Ricci flow equation (\ref{usual0}) (see the remarks after   (\ref{cancel0}) and (\ref{ode})). In this analogy, the linear
non-autonomous {\it partial} differential equation  (\ref{kac}) that ensures that the constraint  $R(g)=-1$ is   preserved  if it
holds initially is    analogous to the linear non-autonomous {\it ordinary}
differential equation   (\ref{ode}) that insures that
$ {\rm vol}(M,g) =1$ is preserved if it holds initially. 
\item
 The reduced conformal Ricci equation uses the ``simpler" operator\\ $L_g=   (n-1)\Delta_g +1$,
 rather than $ L_g= (n-1)\Delta_g -R(g)$, even though the curve  $ 
 g:[0,T) \rightarrow {\cal M} 
 $ is at first only known to be in $ {\cal M} $ (with initial value $ g_0\in {\cal M}_{-1} $).
This is analogous to using the ``simpler"  $c=\frac{2}{n}R_{\rm total}(g)$ 
in the fully reduced classical Ricci equation
$(\ref{usual0})$ (with initial value $ g_0\in {\cal M}^{ 1} $)  rather than $ {{\textstyle\frac{2}{n}}  } \frac{ R_{\rm total}(g)}{
{\rm vol} (M,g)}\,$. In both cases the   task  is  then to show that the constraint is maintained on the  entire domain of
definition of any
 solution of the   
reduced equation that initially satisfies the constraint.
  \end{enumerate}  

\subsection{Increasing scalar curvature for the classical Ricci flow equation}
\label{increasing5} In this section   we  discuss  the somewhat opposite behavior
of
   the classical and conformal Ricci flow equations with respect  to  scalar curvature 
(see   Remark~(iv) in  Section~\ref{summary9}).  
  For similar somewhat opposite behavior with respect to volume, see Proposition~\ref{volel}.
Both results are summarized in Table~1 below. 

For    a curve of metrics   $t \rightarrow g_t$, let 
  $R_{\rm min} (g_t  )=\min\limits_{x\in M}
\{R(g_t)(x)\}$ denote  the time-dependent spatial minima of the scalar curvatures  $R(g_t)$. Augmenting  a result of     Hamilton 
\cite{ham99} on
increasing scalar curvature 
 to
manifolds of Yamabe type $-1$ or 0  gives the following.
 \begin{proposition}[Increasing scalar curvature]
\label{increasescale} Let $g:[0,T) \rightarrow {\cal M} ^1$ be a solution of the 
  classical Ricci  flow equation \begin{eqnarray}\nonumber  \frac{\partial g}{\partial t}\quad&=&-2\,{\rm Ric}(g)+ 
{{\textstyle\frac{2}{n}}  } 
\bar  R (g) g     \;,\qquad g_0\in {\cal M} ^1 \,,     
\end{eqnarray}such that for each $t\in [0,T)$, $R_{\rm min} (g _t)\le0$.
Then \begin{equation}
\label{increasescalar3} \frac{d}{dt}R_{\rm min} (g_t )\ge0\,,
\end{equation} so $R_{\rm min} (g_t )
$ is an increasing   function of $t$. 

If $M$ is of Yamabe type $-1$ or $0$, then   $(\ref{increasescalar3})$ holds for   all classical Ricci flows on $M$
(without any additional assumptions on $R_{\rm min} (g _t)$).
 \end{proposition} {\bf Proof:}    From  (\ref{drg}), (\ref{ricci2}),  and (\ref{ricidentity}),   the scalar
curvature
$R(g)$ satisfies  
 \begin{eqnarray}
\label{first1}\nonumber 
\frac{ \partial R(g)}{\partial t} & = &D\! R(g) \frac{\partial g}{\partial t}= D\! R(g)  \left (-2\,{\rm Ric}(g)+ 
{{\textstyle\frac{2}{n}}  }\bar  R (g) g 
\right)       
\\&=&
  -    \Delta_g R(g) +2|{\rm Ric}(g)|^2_g     -
 {\textstyle\frac{2}{n}}\bar R(g )R(g) 
\nonumber  \\&=&
  -    \Delta_g R(g) +2|{\rm Ric}^T\!(g)|^2_g + {\textstyle\frac{2}{n}} R^2(g)    -
 {\textstyle\frac{2}{n}}\bar R(g )R(g) \,.
\label{minima4}
\end{eqnarray}Note that $\bar R(g)$ is a (time-dependent) real number (but is constant in $x$)  whereas $R(g)$ is a (time-dependent)
scalar function. Since $|{\rm Ric}^T\!(g)|^2_g\ge0$ and since at a minimum point $x_{\rm min}(t)\in M$ of $R(g_t)$,   
$ -    \Delta_{g _t}R(g_t)(x_{\rm min}(t))\ge0$,  (\ref{minima4})  
 yields
 \begin{eqnarray}\fl\quad 
\frac{ d}{d t} R_{\min}(g )& \ge &  
      {\textstyle\frac{2}{n}} R_{\rm min}^2(g )    -
 {\textstyle\frac{2}{n}}\bar R(g ) R_{\rm min}(g )
=  {\textstyle\frac{2}{n}} R_{\rm min} (g ) \left( 
R_{\rm min} (g ) -\bar R(g  )    \right) \ge0\,,
\end{eqnarray} where the last inequality follows
   by the assumption    that
$R_{\rm min} (g )
\le0$ and from the inequality
   $
R_{\rm min} (g ) - \bar R(g  )  \le0$.

If $M$ is of Yamabe type $-1 $ or 0, then   $R_{\rm min} (g ) \le0$ for any $g$ on $M$.
 For if $R_{\rm min} (g ) >0$, then $R  (g ) >0$  
 and thus by Yamabe's Theorem there would exist  a pointwise conformally equivalent metric $g' =pg $, $p\in {\cal P} $,  with  scalar
curvature
$ R(g' )=$ constant $
 >0$ (see Aubin \cite{aub76b} and Schoen \cite{sch84}), contradicting the assumption that $M$ is of Yamabe type $-1$ or $0$ (see
Definition~\ref{defyam}). 
 Thus the   assumption 
$R_{\rm min} (g ) \le0$ and the  conclusion $ \frac{d}{dt}R_{\rm min} (g_t )\ge0$ are  automatically true for {\it any} classical
Ricci flow on a  manifold 
$M$ of Yamabe type $-1$ or $0$. \hfill$\rule{2mm}{3mm}$ \bigskip

 Table~1 compares the classical and conformal Ricci flows with respect to volume and scalar curvature. 
 \begin{table} [htbp] 
\caption[]{A comparison of the volume and scalar curvature for the classical and conformal Ricci flows.  The scalar
curvature inequality for  the classical Ricci flow assumes that either $R_{\rm min}(g_t)\le0$ or that $M$ is of Yamabe type $-1$
or
$0$.}
\label{table:isotropy}
{\small 
 \begin{center}
\begin{tabular}{ l l l l l l l}
\hline 
\\[-2.5mm]
     &  &
  \multicolumn{2}{c}{Volume}      & &     \multicolumn{2}{c}{ Scalar curvature } 
\\  [1.5mm]  
\hline  
\\[-5mm]
 \\
Classical Ricci  Flow            && $ {\rm vol}(M,g_t) =1$ &    &   &$\frac{d}{dt}R_{\rm min}(g_t)\ge0$\\ 
            &                &     \\
Conformal Ricci    Flow     &$\qquad$&  $\frac{d}{dt}{\rm vol}(M,g_t)\le0  $&         &$\qquad$&      $ R(g_t)=-1$&    
\\
          &$\qquad$&        &         & &          &  \\
\hline
\end{tabular}
 \end{center}
}
\end{table}

In Section~\ref{lochom} we shall consider the locally homogeneous case for both the classical and conformal Ricci flow. 
 In that case
the scalar curvatures  are constant and so 
 $R(g)=R_{\rm min} (g )  =
 \bar R(g  )= R_{\rm max}(g)$. Thus under the assumption that $R(g)\le0$ (or that $M$ is of Yamabe type $-1$ or $0$), we can
conclude that the flow of spatially constant scalar curvatures is increasing as a function of time.

 \section{Local existence and uniqueness for the conformal Ricci flow equations}
\setcounter{equation}{0}
\label{existenceCRF}In this section we shall need information regarding  Sobolev spaces of Riemannian
metrics and tensors. We refer to        Ebin
\cite{ebin70}, Berger-Ebin
\cite{be69},    Fischer-Marsden \cite{fm75b},   Marsden-Ebin-Fischer \cite{mef72}, and Palais (\cite{pal65}, \cite{pal68}) for the
necessary background.
 For $n\ge3$ and   $s>\frac{n}{2}  $, let
 $ S_2^s$  denote the space of  symmetric 2-covariant tensor fields on $M$ of Sobolev class  $H^s$,   $ {\cal M} ^s$ denote the
space of
$H^s$ metrics on
$M$,
   $ {\cal F} ^s=H^s(M, \mbox{\boldmath$R$} )$  denote  the space of
$H^s$ real-valued functions, $ {\cal P} ^s=H^s(M,\mbox{\boldmath$R$} ^+)$ denote the multiplicative Abelian group of $H^s$ positive
functions, and
$ {\cal D} ^{s+1}$ denote the group of 
$H^{s+1}$ diffeomorphisms of $M$. When a superscript is   omitted, we shall continue to mean     $C^\infty$.

There are two issues involved in  proving local existence and uniqueness for the conformal Ricci flow.
 The
first   is concerned with the non-ellipticity of the Ricci tensor when viewed as  a nonlinear partial differential
operator on the space of Riemannian metrics. This non-ellipticity arises solely because of the  covariance of the Ricci tensor with
respect to the group of  diffeomorphisms of $M$, which, locally, is the pseudo-group of coordinate transformations of $M$. Modulo
this group, the Ricci operator  is elliptic. We call   an operator
$E$ that is elliptic modulo the pseudo-group of coordinate transformations a   {\bf
quasi-elliptic operator} and the resulting ``heat equation" derived from that operator,
 \begin{eqnarray}\label{unnormalized} \frac{\partial
\psi}{\partial t}&=&E(\psi) \,,
\end{eqnarray} 
 a {\bf quasi-parabolic equation}. 
 Thus  the   unnormalized classical Ricci flow  equation,  where  the solutions are no
longer volume preserving, \begin{eqnarray}\label{unnormalized}
\frac{\partial g}{\partial t}&=&-2 {\rm Ric}(g)\,, 
\end{eqnarray} is  {\bf quasi-parabolic} and  the problem of the non-ellipticity of the Ricci tensor has been overcome by Hamilton
(\cite{ham82},\cite{ham95})  and DeTurck
\cite{det83} who have shown that there exists a   unique short-time
        solution  to the initial value problem of (\ref{unnormalized})      for any closed
$n$-manifold,
$n\ge2$,  with arbitrary initial metric
$g_0\in {\cal M} $.

The second issue  regarding existence and uniqueness is that     the conformal Ricci flow equations as well as being quasi-parabolic
are also elliptic because   the constraint
$R(g)=-1$ leads to a linear inhomogeneous elliptic equation (\ref{elliptic6}) for the conformal pressure $p$.
Thus overall the conformal Ricci flow equations are       a    {\bf linear-elliptic quasilinear  quasi-parabolic system}. It is
this  latter    issue of  the  additional elliptic equation  for the conformal pressure that we address here.

We proceed as follows.
 We first show that the (fully) reduced conformal Ricci flow equation defined on $ {\cal M}_{-1} $  is a smooth   bounded
perturbation of the unnormalized classical Ricci equation and thus is a vector field sum of a smooth and a densely defined vector
field on 
$ {\cal M}  $.
We then   consider an evolution equation that involves  only the smooth
term.  That this auxiliary evolution equation has a smooth flow then follows from the usual Picard 
iteration method for smooth vector fields on Hilbert (or Banach) manifolds.
 We then ``add" back in the   Ricci term to this auxiliary equation to get the  reduced conformal Ricci flow equation.  That this
equation has a flow then follows from a nonlinear Trotter product formula.

This bootstrap method from a smooth flow to the flow of interest    is inspired by a similar method developed by Ebin-Marsden
\cite{em70} that shows    existence and uniqueness of solutions to the   Navier-Stokes equations.  In this analogy,  the
Navier-Stokes equation is best  viewed as the   sum of the  Euler equation  and the heat equation  taken on the tangent
bundle $T {\cal D} _ \mu $  of the group
$ {\cal D} _ \mu $ of volume-preserving diffeomorphisms, where $ \mu $ is the volume element of a fixed Riemannian metric (here we
denote volume elements by $ \mu $, rather than $ d   \mu $).   On
$T {\cal D} _
\mu
$ there is no loss of derivatives for the Euler  
  equation   and in fact it is a    smooth vector field      and hence has
a smooth flow by the usual methods of ordinary differential equations applied to Hilbert manifolds. To this smooth flow is added
the heat equation, or more accurately, the heat equation and a Ricci curvature term (see Ebin-Marsden \cite{em70}, p.\ 162).  The
resulting system
 then  has a flow  by a nonlinear
Trotter product formula and,    pulling back this flow  to the space of
divergence-free vector fields on
$M$, proves existence and uniqueness for  the classical Navier-Stokes equation.

  In this analogy, the constraint phase space $T {\cal D} _ \mu $ corresponds to the constraint space $ {\cal M}_{-1} $, 
the  smooth  Euler equation on $T {\cal D} _ \mu $ corresponds to the smooth auxiliary equation 
 (\ref{auxheat}) below,   the heat equation on $T {\cal D} _ \mu $ corresponds to the unnormalized Ricci flow equation, and the 
Navier-Stokes equation on $T {\cal D} _ \mu $  corresponds to the reduced conformal Ricci  flow equation on $ {\cal M}_{-1} $. 
However, a difference   between these two systems is that whereas the components of the Navier-Stokes equation, namely, the 
Euler and heat equations, both live on the constraint space $T {\cal D} _ \mu $, the components of the reduced
conformal Ricci flow equation, the auxiliary equation  (\ref{auxheat}) below  and the unnormalized Ricci flow equation,
do not separately flow in the constraint space $ {\cal M}_{-1} $ (only their sum does; see Figure~1 in Section~\ref{nonortho}).

The reduced conformal Ricci flow equation can also be viewed as a reaction-diffusion equation.
   In this view the reaction part of the equation  corresponds to the 
smooth auxiliary equation  (\ref{auxheat})   and the diffusion part corresponds to the   unnormalized Ricci equation.

 For $ s>\frac{n}{2}+1$, let $   {\cal M}_{-1} ^s=\{\,g\in {\cal M} ^s\mid R(g)=-1\,\}$. 
 
\begin{theorem}{\bf(Local existence and uniqueness of solutions for the conformal Ricci flow equations)} 
 \label{thmMain}Let $M$ be a closed connected   oriented $n$-manifold, $n\ge3$,    $ 
s>\frac{n}{2}+2 $,  and   $g_0\in {\cal M}^s_{-1} $. Then there exists a $T>0$, a   
 $C^0$  curve
 $ 
 g\,\colon[0, T) \rightarrow {\cal M}^{s } 
 $, 
 $C^1$  as a  curve
  $ 
 g\,\colon[0,T) \rightarrow {\cal M}^{s-2} 
 $,   
and a   $C^0$ curve $p\,\colon[0,T) \rightarrow {\cal F}^s$, 
such that 
$g$ has initial condition $g(0)=g_0$ and  such that 
$g$ and $p$ satisfy 
    the   conformal Ricci flow  equations, \begin{eqnarray}\label{parabolic65} \frac{\partial
g}{\partial t}+2\big({\rm Ric}(g)+{\textstyle \frac{1}{n}}g\big)&=&-p   g          \\
\label{elliptic65} \;\;\;\qquad\qquad\quad  R(g)&=&-1 \,.
\end{eqnarray} 

If $ s>\frac{n}{2}+3$, then    $g$ and $p$ are   unique.
\end{theorem} {\bf Proof:} We sketch a proof that views the reduced conformal
Ricci flow equation as a bounded perturbation of the  unnormalized Ricci equation.
 Our starting point is to
note that   Proposition~\ref{elimination1}  asserts that the
 conformal Ricci initial value problem is   equivalent to the reduced conformal Ricci initial value problem, 
 \begin{equation}
\label{reduced121}   \frac{\partial g}{\partial t} +2  {\rm Ric}(g)=-  
      2\big( L ^{-1} _g( |{\rm Ric} (g)|_g^2 )\big) g\;,\qquad g_0\in {\cal M}_{-1} ^s\,,
\end{equation} 
 where 
 associated with the yet-to-be-determined  $C^0$ curve  $g:[0,T ) \rightarrow {\cal M} ^s$  
is the  $C^0$  curve of   invertible   elliptic operators $ L_g 
 = (n-1)\Delta_g +1 : {\cal F}^s \rightarrow {\cal F}^{s-2}
$
  with inverses
$  L_{g }^{-1}: {\cal F} ^{s-2}\rightarrow {\cal F}^s $.
 Note that if such a curve $g$ exists, then from Proposition~\ref{elimination1} (adapted to the Sobolev setting), 
$g$ maps to 
$ {\cal M}^{s }_{-1}
$ as a $C^0$ curve  and    to $ {\cal M}^{s-2}_{-1} $ as a $C^1$ curve.

The reduced conformal Ricci flow equation    (\ref{reduced121})
differs from the 
 unnormalized Ricci flow equation (\ref{unnormalized})  by 
the term $-  
      2 ( L ^{-1} _g( |{\rm Ric} (g)|_g^2 ) ) g $.
 Thus we define the difference vector field on $ {\cal M} ^s$, 
 \begin{equation}
\label{aux}  Z
  : {\cal M} ^s \longrightarrow S_2 ^s\;,\qquad g \longmapsto Z(g) = -2\left(L_g ^{-1}  (    |{\rm Ric}(g) |_g^2   ) \right)g
\end{equation} and consider the corresponding ``heat" equation \begin{equation}
\label{auxheat}  \frac{\partial g}{\partial t}  =Z(g)=-  
      2\big( L ^{-1} _g( |{\rm Ric} (g)|_g^2 )\big) g=- ({\textstyle\frac{2}{n}} +p)g\,,
\end{equation} where     $p=  {   2\, ( L_g ^{-1} (      |{\rm Ric}(g) |^2   )  -  \textstyle\frac{1}{n}})
   $. 
  The key fact is  that $ Z$ 
does not lose derivatives, i.e., $Z$  actually maps $H^s$  to $H^s$, 
and 
 is    a smooth vector field on the indicated spaces.
 
To show that 
$Z$    maps $H^s$  to $H^s$, 
we first make the following remarks. Let 
  $g\in {\cal M} ^s$, $ s>\frac{n}{2}+2 $. Then the components of the inverse $g ^{-1} $    are a rational 
combination of the components of $g$
with non-zero denominator $\det g$. Thus
$g ^{-1} $ as a function of $g$ is a smooth map. 
   From the
Schauder ring  property of
Sobolev functions,
$ g ^{-1} $ is   an $H^s$ contravariant metric and then from the local formula (\ref{riccilocal3}) for the Ricci tensor   and the
multiplicative properties of Sobolev functions, 
$ {\rm Ric}(g)
\in S^{s-2}_2$ and   
 $  
    |{\rm Ric}(g) |_g^2   \in {\cal F} ^{s-2}$.  Note that the  Schauder ring property
states that 
$H^s$ is a ring if
$ s>\frac{n}{2} 
$.  Thus to be able to conclude that    $  |{\rm Ric}(g) |_g^2 $ is of class  $H^{s-2}$     requires $ s>\frac{n}{2}+2$ 
 since $ {\rm Ric}(g)
\in S^{s-2}_2$.
 
Since
 $  |{\rm Ric}(g) |_g^2\in {\cal F} ^{s-2}$, from 
   ellipticity of $L_g$,
$L_g ^{-1}  (    |{\rm Ric}(g) |_g^2   )      \in {\cal F} ^s
$ and thus $Z(g)=-(L_g ^{-1}  (    |{\rm Ric}(g) |_g^2   )   )g   \in S_2 ^s
$.
 Thus,
remarkably,   $Z$
does not lose  derivatives relative to $g\in  {\cal M} ^s$.
 Although
this does not prove that     $Z$ is smooth, it does make it plausible since now the range space  does not have a weaker
topology than the domain space. 

To show that $Z $ is smooth,  
 we proceed as in Fischer-Marsden \cite{fm75b} and  consider the Ricci tensor as a map, the {\bf Ricci map} (or {\bf
operator}), 
  \begin{equation}
\label{LABEL} {\rm Ric}: {\cal M} ^s \longrightarrow S^{s-2}_2\;,\qquad g \longmapsto   {\rm Ric}(g) \,,
\end{equation} 
which can  also be interpreted as  a    vector field on $ {\cal M} ^{s-2}$ defined on the dense domain $ {\cal M}
^s\subset {\cal M} ^{s-2}$.  This vector field  loses derivatives and is only  an everywhere defined  smooth   vector field  when
$s=\infty$, in which case the   theorem of ordinary differential equations on Banach manifolds fails.

First we show that for $ s>\frac{n}{2}+1 $,  $ {\rm Ric}$ is a smooth function on the indicated spaces.  Because differentiation is
a continuous linear map between the spaces indicated, the smoothness of
$ {\rm Ric}$ depends on the multiplications that occur in computing $ {\rm Ric}(g) $.  The second-order derivatives appear linearly
with components of $g^{-1} $ as coefficients 
whose terms are various  natural (i.e., non-metric) contractions of  the homothetically invariant expression
$ g^{ij}\frac{\partial^2g_{kl}}{\partial x^m\partial x^n}$
 (i.e., invariant under the
homothetic transformation $ g \mapsto cg$, $c>0$)
 which we write symbolically as 
 $g ^{-1}  \otimes
D^2g$, to yield a symmetric 2-covariant expression. 
  Since   $g ^{-1} $ is a 
  smooth function of   
$g$,    by the multiplicative properties  for Sobolev spaces, for $ s>\frac{n}{2}$,    the pointwise bilinear map $g ^{-1} \otimes
D^2g$ induces a multiplication
$H^s\times H^{s-2} \rightarrow H^{s-2}$ which  is continuous bilinear  and hence smooth. Thus    the second order terms are
smooth functions of
$g$ and we note that $ s>\frac{n}{2}$ suffices for smoothness of the second-order terms.

The first order terms  are 
 various natural contractions of the  homothetically invariant generic form 
$g_{ij}g^{kl}g^{mn} g^{pq}\frac{\partial g_{ab}}{\partial x^c}  \frac{\partial  g_{de}}{\partial x^f}$ 
to yield a symmetric 2-covariant expression. 
  Thus  the first order
terms  are rational combinations of the components of 
$g$ with non-zero denominator $\det g$ and quadratic functions of
$Dg$. Thus again
 by the multiplicative properties  for Sobolev spaces,  for   $
s>\frac{n}{2}+1$, $H^{s-1}$ is a ring  and    the pointwise multilinear map  $g\otimes  g ^{-1}\otimes  g
^{-1} \otimes  g ^{-1} \otimes D g\otimes D g$ induces a multiplication
$H^{s }\times H^{s }\times H^{s }\times H^{s }\times H^{s-1} \times H^{s-1} \rightarrow H^{s-1}$ which  is continuous multilinear 
and hence smooth. Thus    the first order terms are smooth functions of
$g$ and we note that since $Dg$  appears quadratically,  $ s>\frac{n}{2}+1$ is necessary for smoothness of the first-order terms. 
  Thus $ {\rm Ric}$ is $C^\infty$ for the indicated spaces   under the
condition $ s>\frac{n}{2}+1 $. 

For $s>\frac{n}{2}+2$, $H^{s-2}$ is a ring and   
  again by the multiplicative properties of Sobolev spaces, the 
pointwise squared Ricci norm map
\[\fl\qquad\quad {\cal M} ^s
\rightarrow S^{s-2}_2 \times S^{s-2}_2 \rightarrow {\cal F} ^{s-2}\,,\qquad  g \longmapsto  ({\rm Ric}(g),{\rm Ric}(g)) \longmapsto 
   | {\rm Ric}(g) |_g^2\]
 is a composition of the smooth map   $g \mapsto  {\rm Ric}(g) \times {\rm Ric}(g) $ and  a multiplication
$H^{s-2}\times H^{s-2} \rightarrow H^{s-2}$ which  is continuous bilinear  and hence smooth. 
Thus $ g \mapsto  
   | {\rm Ric}(g) |_g^2$ is smooth.  We note that since  $ {\rm Ric}(g) \in S_2^{s-2}$
appears quadratically in $| {\rm Ric}(g) |_g^2\,$,  
   $ s>\frac{n}{2}+2$ is necessary for smoothness of $ g \mapsto  | {\rm Ric}(g) |_g^2\,$.

 In a similar manner, one shows that $ g
\mapsto  L  _g  $ is smooth as a function of $g$ and that the map    $g \mapsto  L_g ^{-1} (  | {\rm
Ric}(g) |_g^2)$  is smooth.  
 Thus the vector field  $ Z
  : {\cal M} ^s \longrightarrow S_2 ^s,$ $g \mapsto  -2( L_g ^{-1} (  | {\rm
Ric}(g) |_g^2))g $
  is smooth. Since $ {\cal M} ^s$ is a smooth  Hilbert manifold, the existence and uniqueness of local  flows follows by the
fundamental theorem of ordinary  differential equations on  Hilbert manifolds (see, for example, Lang \cite{lan95}, Theorem 2.6).

If we now consider the   sum of the two vector fields  $-2  {\rm Ric} $ and  $Z$,    \begin{equation}
\label{LABEL}\fl\quad   - 2  {\rm Ric}+Z
  : {\cal M} ^s \longrightarrow S_2 ^{s-2}\;,\qquad g \longmapsto  
 -2  {\rm Ric}(g) -2\left(L_g ^{-1}  (    |{\rm Ric}(g) |_g^2   ) \right)g\,,
\end{equation}  the corresponding
 ``heat" equation  
 \begin{equation}
\label{redcrf}  \frac{\partial g}{\partial t}  =- 2{\rm Ric}(g) -
      2\big( L ^{-1} _g( |{\rm Ric} (g)|_g^2 )\big) g 
\end{equation} is the reduced conformal Ricci flow equation.
 This equation does lose derivatives but  this can be overcome by  using a  nonlinear Trotter product formula (see
Ebin-Marsden~\cite{em70}, pp.\ 142, 157 and Taylor  \cite{tay96c}, Section~15.5). To do so, one approximates solutions of
(\ref{redcrf})  by successively solving the two   component equations 
\begin{eqnarray}
\label{reduced99}   \frac{\partial g}{\partial t}  &=& -  
      2\big( L ^{-1} _g( |{\rm Ric} (g)|_g^2 )\big) g
\\
 \label{reduced100}   \frac{\partial g}{\partial t}& =&-  
      2 {\rm Ric}(g) 
\end{eqnarray} over small time intervals  and then composing the resulting solution operators.

As we have seen,   the first equation is smooth and the second equation is the unnormalized
Ricci flow  equation which, as discussed earlier, is known to have unique  local semi-flows $F_t: {\cal M} ^s \rightarrow {\cal M}
^s$,
$t\in [0, \tau )$, $ \tau >0$, generated  by the     densely defined    vector field $  - 2  {\rm Ric}  : {\cal M} ^s\subset
{\cal M} ^{s-2}
\rightarrow S_2^{s-2}  $  (for   existence and uniqueness results of the classical  Ricci flow equation in the $H^s$-setting, see
Fischer~\cite{fis04}).  Thus the reduced conformal Ricci flow  equation is the vector field sum of a smooth vector field and a
densely defined
  vector field, the former of which has  unique local flows by  classical methods of ordinary differential equations  and the
latter of which  is   known to    have unique local semi-flows by quasi-parabolic methods.
The  resulting  vector field sum then fits 
  the setting for  the nonlinear 
  Trotter product formula as it appears in Ebin-Marsden \cite{em70}. After
checking the     hypotheses of that formulation as they apply to (\ref{reduced99}--\ref{reduced100}),   we can conclude  
 that the reduced conformal Ricci equation  also has unique local semi-flows.

That uniqueness requires $ s>\frac{n}{2}+3$ (instead of $  s>\frac{n}{2}+2$ as    for existence) is discussed in
Fischer~\cite{fis04}.
 \hfill$\rule{2mm}{3mm}$ 
\medskip

\noindent {\bf Remarks: }
\begin{enumerate}
\item
Each term on the left hand side of (\ref{parabolic65})  is separately a $C^0$ curve in
$S_2^{s-2} $ whereas the right hand side (and thus the {\it sum} of the terms on the left) is a
$C^0$ curve in $S_2^s$ (actually in $- {\cal M} ^s$ for non-static curves since $p>0$).
However, in the rearrangement $\frac{\partial g}{\partial t} =  {\rm
Ric}(g)-(p +  \frac{2}{n} g) g        $, each side is a $C^0$ curve in $S_2^{s-2} $. 
\item
Since $ s>\frac{n}{2}+2 $, the solutions are genuine, i.e.,  the equations are satisfied in the classical
sense.
\item
Although the auxiliary equation (\ref{aux}) is smooth and in particular does not lose derivatives, 
the conformal Ricci flow
equations themselves do lose derivatives  since  if $g _t $ is $H^s$, $ {\rm Ric}(g_t)$ is only  $H^{s-2}$. 
\hfill$\rule{2mm}{3mm}$ 
 \end{enumerate}  
 
\medskip By adding  $   {\textstyle\frac{2}{n}}g $  
to the   right hand side of  (\ref{reduced99}) and subtracting $   {\textstyle\frac{2}{n}}g $  
 from  
(\ref{reduced100}),  (\ref{reduced99}--\ref{reduced100}) can be re-written  in the more
geometrical (but analytically equivalent) form
  \begin{eqnarray}
\label{equivalentform1}   \frac{\partial g}{\partial t}  &=&Z_1(g)= -  2\big( L ^{-1} _g( |{\rm Ric} (g)|_g^2 )\big) g +
{\textstyle\frac{2}{n}}g = - pg 
 \\
\label{equivalentform2}   \frac{\partial g}{\partial t}& =&Z_2(g)=-  
   2 (   {\rm Ric}(g)+ {\textstyle\frac{1}{n}} g)\,, 
\end{eqnarray}
where
 $p=p(g)=  2 L_g ^{-1}  (
|{\rm Ric}(g)|^2_g )
-{\textstyle\frac{2}{n}}  $ and where the vector field sum of the right hand sides  is again the reduced conformal Ricci
flow
 equation.
Analytically,  the first vector field is smooth (see   Proposition~\ref{summarizing}) and 
thus has a local flow. 
The second vector field is densely defined on $ {\cal M} ^s$  but modulo the bounded perturbation term  $- {\textstyle\frac{2}{n}}
g$ 
  is      the  unnormalized Ricci flow
equation    and thus     also has  unique  local semi-flows.

The vector fields (\ref{equivalentform1}--\ref{equivalentform2}) have the geometrical feature that they are pointwise orthogonal on
$ {\cal M}_{-1} $ since \begin{equation}
\label{onsince}  Z_1(g)\cdot Z_2(g)=(-pg)\cdot-2( {\rm Ric}(g) + {\textstyle\frac{1}{n}} g)=2p(R(g)+1)=0\,.
\end{equation} We shall further explore this feature   in Section~\ref{nonortho}.   Here we   note that 
neither vector field  $Z_1$ nor $Z_2$ preserves   the volume  nor the scalar curvature of the initial metric $g(0)=g_0$  but that
 $Z_1$       preserves the pointwise conformal class of  the initial
metric. Thus  if $g$ is a solution to  (\ref{equivalentform1}), 
then $  \frac{\partial g(t)}{\partial t}   =  - p(g(t))g(t)$ so that  $g$ must satisfy \begin{equation}
\label{LABEL} g(t)=g_0e^{-\int_0^tp(g(t'))dt'}\in Pg_0\,.
\end{equation}   Thus the flow of (\ref{equivalentform1})   lies in the  $ {\cal P} $-orbit of the
initial metric  $ g_0$ and thus preserves the pointwise
conformal class of  
$g_0$.

We also remark that a 
  split-step algorithm  based on the pair of equations (\ref{equivalentform1}--\ref{equivalentform2}) can also be used as a basis
for a numerical solution of the reduced conformal Ricci flow equation.

Lastly, we note that   (\ref{equivalentform2})
is similar  in structure to the historically   first Ricci flow  equation  (see Hamilton \cite{ham82}, p.\ 256),
  \begin{equation}
\label{andusing}  \frac{\partial g}{\partial t} =- 2( {\rm Ric}(g)-{\textstyle\frac{1}{n}} R(g)g) \,,
\end{equation} 
 which  has    a backwards heat equation in $R(g)$. Indeed, 
calculating    as in  (\ref{minima4}), the scalar curvature of a solution  $g $ to (\ref{andusing}) must  
satisfy  
  \begin{eqnarray}\nonumber   \frac{\partial}{\partial t}R(g) & = &
 D\! R(g) \frac{\partial g}{\partial t}= D\! R(g)  \left (-2\,{\rm Ric}(g)+ 
{{\textstyle\frac{2}{n}}  }  R (g) g\right) \\
\nonumber 
 & = & -    \Delta_g R(g) +2
|{\rm Ric}(g)|^2_g    
+{\textstyle\frac{2}{n}}D\!R(g) ( R(g)g)          \\   & = &  -    \Delta_g R(g) +2
|{\rm Ric}(g)|^2_g 
 + {\textstyle\frac{2}{n}} (  n-1)\Delta_g R(g)-{\textstyle\frac{2}{n}}  R^2(g)  \nonumber \\  & = & {\textstyle\frac{n-2}{n}} 
   \Delta_g R(g) +2
|{\rm Ric}(g)|^2_g 
 - {\textstyle\frac{2}{n}}R^2(g)\,,  
\label{mustsatisfy}
\end{eqnarray}
  where,    because of the change in sign of the Laplacian term from 
$ - \Delta_g$ to  $  {\textstyle\frac{2}{n}}  ( n-1)\Delta_g$,     is a backwards heat equation in $R(g)$. Thus
 (\ref{andusing})  {\it in general} cannot have short-time solutions. However, in the locally homogeneous case (i.e., where there
exists a local isometry between neighborhoods $U_x, U_y$  of every pair of points $ x,y $ in $(M,g)$),
$R(g)=\bar R(g)$ and so
(\ref{andusing})  is equal to the classical Ricci flow equation 
and thus is well-posed (note that in this case 
 $\Delta_g R(g)=0$; see also     
  Proposition~\ref{yam21} and the remark following that Proposition). 

 Although       (\ref{andusing}) 
is similar to 
  (\ref{equivalentform2})
(and the equations are  identical   when $R(g)=-1$), 
   (\ref{equivalentform2})
 does not suffer from a backwards heat equation in $R(g)$, since 
  \begin{eqnarray}\nonumber   \frac{\partial}{\partial t}R(g) & = &
 D\! R(g) \frac{\partial g}{\partial t}= D\! R(g)  \left (-2\,{\rm Ric}(g)-
{{\textstyle\frac{2}{n}}  }   g\right) \\
\nonumber 
 & = & -    \Delta_g R(g) +2
|{\rm Ric}(g)|^2_g    
-{\textstyle\frac{2}{n}}D\!R(g) ( g)          \\   & = &  -    \Delta_g R(g) +2
|{\rm Ric}(g)|^2_g 
  -{\textstyle\frac{2}{n}}   \,,  
\label{mustsatisfy}
\end{eqnarray}
 which is a proper heat equation for $R(g)$.

Similarly, we remark in passing that for  the conformal Ricci flow  equations, since $R(g)=-1$, 
   there is no issue of whether
or not 
$R(g)$ solves a backwards or proper heat equation.

\section{The geometry of the conformal Ricci flow}
\label{nonortho}
To put the fact that the conformal Ricci flow equations  do not lose derivatives {\it relative} to the unnormalized Ricci equation 
into further geometrical perspective,  we consider a  non-orthogonal $L_2$-splitting of $S^s_2$.
  We then apply this splitting  to $ {\rm Ric}(g) $ and show why this splitting  considerably simplifies when the
    constraint $R(g)=-1$ is satisfied.   Indeed, as we shall see, when $ R(g)=-1$,     the map   $g \mapsto -pg$    does not lose
derivatives   and is smooth as a function of
$g$.  Thus
 the reduced conformal Ricci  flow equation is a bounded perturbation of the 
unnormalized Ricci equation.

 For $n\ge3$ and   $s>\frac{n}{2}+1 $, we
  continue to let \begin{equation}
\label{LABEL} R: {\cal M} ^s \longrightarrow {\cal F} ^{s-2}\;,\quad g \longmapsto R(g)
\end{equation}  denote the smooth
 scalar curvature map  in the Sobolev setting,
with derivative at $ g\in {\cal M} ^s$  
\begin{equation}
\label{LABEL}D\!R(g):S^s_2 \longrightarrow {\cal F} ^{s-2}\,,
\end{equation}  and with   
$\bar S^s_2(g)=\ker D\!R(g)\subset S^s_2$    the kernel of $D\!R(g)$.

 If
$g\in {\cal M} ^s$ and
$p\in {\cal P} ^s$, then the product $p g\in {\cal M} ^s$. Let 
  \begin{equation}
\label{LABEL} {\cal P} ^sg=\{\,p g\mid p\in {\cal P} ^s\,\}\subset {\cal M} ^s
\end{equation} denote the space of {\bf (pointwise) conformal deformations of} $g$,
or   the {\bf orbit  under the group action of $ {\cal P} ^s$ on} $ {\cal M}
^s$.  The orbit $ {\cal P} ^sg$ is a closed submanifold of $ {\cal M} ^s$.

  If
$g\in {\cal M} ^s$ and
$\varphi\in {\cal F} ^s$, then the product $\varphi g\in S_2^s$. Then  the tangent space
  \begin{equation}
\label{LABEL} T_g( {\cal P} ^sg)\approx{\cal F} ^sg=\{\,\varphi g\mid \varphi\in {\cal F} ^s\,\}\subset S_2^s
\end{equation} is the space of {\bf infinitesimal (pointwise) conformal deformations of} $g$.
Both $\bar S^s_2(g)$ and $ {\cal F} ^sg$ are closed subspaces of 
$S^s_2$.

Let \[S_2^{T}\!(g)^s=\{\,h\in S^s_2\mid {\rm tr}_gh=0\,\}\]
denote the space of traceless (with respect to $g\in {\cal M} ^s$) symmetric 2-covariant  tensor fields on $M$.  Thus we have the
pointwise orthogonal (and hence $L_2$-orthogonal) splitting \begin{equation}
\label{henceorthogonal} S_2^s= {\cal F} ^sg\oplus S_2^{T}\!(g)^s\;,\qquad h= {\textstyle\frac{1}{n}}( {\rm tr}_g h)g+h^T
\end{equation} where $h^T=h-{\textstyle\frac{1}{n}}({\rm tr}_g h)g$ is the traceless part of $h$. This splitting can be written more
geometrically as \begin{equation}
\label{geometricallyas} T_g   {\cal M} ^s=T_g( {\cal P} ^sg)\oplus T^\perp_g( {\cal P} ^sg)\,,
\end{equation}where $T^\perp_g( {\cal P} ^sg)\approx S_2^{T}\!(g)^s$ represents the directions orthogonal to the 
infinitesimal pointwise conformal deformations of $g$. We also remark that  at $g\in {\cal M} ^s$, the orthogonal space  $
T^\perp_g( {\cal P} ^sg)$ can be  expressed as the  tangent space  of the {\bf manifold of metrics ${\cal M} _ {d\mu_g}$ with fixed
volume element} 
$
 d \mu _g$, the volume element of $g$ (see Figure~1). Thus, if \begin{equation}
\label{LABEL}  {\cal M} _ {d\mu_g} =\{\, g'\in {\cal M}^s \mid d \mu _{g'} =  d \mu _{g }\,\}\,,
\end{equation}   then  $T_g {\cal M}^s _ {d\mu_g} \approx S^T_2(g)\approx  T^\perp_g(
{\cal P} ^sg)$ since $D( d \mu _g)h= {\textstyle\frac{1}{2}} ({\rm tr}_gh)d \mu _g=0$.

   Now we combine the two orthogonal splittings (\ref{eqn:splitting4}) and
(\ref{henceorthogonal}) 
to give a third non-orthogonal splitting.
 \begin{proposition}{\bf (A non-orthogonal   $L_2$-splitting of $S^s_2$)}
\label{non4}  For $n\ge3$, $ s>\frac{n}{2}+1 $,  and $g\in {\cal M} ^s$,  assume that the scalar
curvature
$R(g)$ is such that  the operator \begin{equation}
\label{operator2}L_g: {\cal F}^s \longrightarrow {\cal F} ^{s-2}\;,\quad \varphi \longmapsto   L_g\varphi=
(n-1)\Delta_g\varphi-R(g)\varphi
\end{equation}
 is an isomorphism.  Let $ L_g ^{-1} : {\cal F}^{s-2} \longrightarrow {\cal F} ^{s }$ denote its inverse.
Then  $S_2^s$ splits into a
non-orthogonal  direct sum 
 \begin{equation}
\label{directsum0} S_2^s=\ker D\!R(g)\oplus {\cal F}^s g=\bar S^s_2(g)\oplus {\cal F}^s g\;,\qquad 
 h=\tilde h+\varphi g
\end{equation} 
   where 
 $ \varphi =L_g ^{-1} (D\!R(g)h)\in {\cal F} ^s$,   $\varphi g\in {\cal F} ^sg$, and   \begin{equation}
\label{LABEL} \tilde h=h-\varphi g=h- (L_g ^{-1} (D\!R(g)h) ) g\in \bar S^s_2(g)\,.
\end{equation} 
 
 Equivalently, but more geometrically, if $ \rho =R(g)\in {\cal F} ^{s-2}$ is such that $ {\cal M} _ \rho =R ^{-1} (\rho  )=\{\,g\in
{\cal M} ^s\mid R(g)= \rho \,\}$ is a submanifold of
$ {\cal M} ^s$, the tangent space
$T_g {\cal M} ^s$  of $ {\cal M} ^s$ at $g$  splits  
\begin{equation}
\label{directsum1}T_g {\cal M} ^s   =T_g {\cal M}_{\rho } ^s\oplus  T_g({\cal P}^s  g )
\end{equation} into two non-orthogonal tangent spaces corresponding to the  two   closed submanifolds $ {\cal M}^s_{\rho }
$ and
${\cal P} ^s g$ of $ {\cal M} ^s$ intersecting transversally at $g$ (see Figure~1).

\end{proposition}\noindent
 {\bf Notational remark}:
When we  let $
\bar S^s_2(g)=
\ker D\!R(g) $, we are considering the kernel of the map $ D\!R(g):  S^s_2\rightarrow  {\cal F} ^{s-2}$.
Also, for
   $h\in S^s_2$,
we use $\tilde h\in{\bar  S}^s_2$ to denote the component of $h$ in $\bar S^s_2$ with respect to this non-orthogonal splitting, in
contrast to 
$\bar h\in  \bar S^s_2$ which denotes the component of $h$ in $\bar S^s_2$ using the $L_2$-orthogonal splitting given by
$\ref{eqn:splitting4}$.

\medskip
\noindent
 {\bf Proof:} For $g\in {\cal M} ^s$, we first     note that \begin{equation}
\label{LABEL} \bar S^s_2(g)\cap {\cal F} ^s g=\{0\}\,,
\end{equation} since if $\tilde h=\varphi g\in  \bar S^s_2(g)\cap {\cal F} ^s g$, then since
$s>\frac{n}{2}+1 $,  $ R: {\cal M} ^s \rightarrow {\cal F} ^{s-2}$ is a smooth mapping and
from
(\ref{last4d}),
  \begin{eqnarray}  D\!R(g)(\varphi g) =
  (n-1)\Delta_g \varphi  - R(g)\varphi
 =L_g\varphi 
 =0      \,.  
  \label{123}
\end{eqnarray} 
Since the curvature assumption on $R(g)$ is that $L_g$ is an isomorphism,  (\ref{123})  has unique solution $\varphi=0$.

To show that the direct sum $\bar S_2^s\oplus {\cal F}^s g$  of  the    indicated closed subspaces exhausts $S_2^s$, 
for $h\in S^s_2$  let 
$ \varphi=L_g ^{-1} (D\!R(g)h)$. Since $ D\!R(g)h \in {\cal F} ^{s-2}$,  by  the ellipticity of $L_g$,  $\varphi=L_g
^{-1} (D\!R(g)h)\in {\cal F} ^s$,
  $\varphi g=L_g ^{-1} (D\!R(g)h)g\in {\cal F} ^sg$, and thus
  \begin{equation}
\label{LABEL} \tilde h=h-\varphi g  \in S_2^s
\end{equation} 
We need to show that  $\tilde h\in \ker D\!R(g)=\bar S_2^s\,$, which follows from \begin{eqnarray}  D\!R(g)\tilde h & = & 
D\!R(g)h-D\!R(g)(\varphi g)  =    D\!R(g)h -L_g\varphi
\nonumber \\   & = & D\!R(g)h -L_g(L_g ^{-1} (D\!R(g)h)   )=0 \nonumber 
\end{eqnarray} 
 Thus
   $h=\tilde h+\varphi g$ splits $h$ according to (\ref{directsum0}).
  \hfill$\rule{2mm}{3mm}$
\bigskip

Note that although neither the $L_2$-orthogonal splitting (\ref{eqn:splitting4}) nor the pointwise orthogonal splitting
(\ref{henceorthogonal})   requires a curvature condition on
$g$, the  $L_2$-non-orthogonal splitting  (\ref{directsum0})   does  since it is required that  $g$ is such that    
$L_g=(n-1)\Delta_g-R(g)$ is an isomorphism.

Note also that the summands in the $L_2$-non-orthogonal splitting (\ref{directsum1}) of $S_2^s$ can be interpreted geometrically as
a  splitting of the tangent space
$ T_g {\cal M} ^s$ into the tangent
spaces of two closed transversally intersecting submanifolds of $ {\cal M} ^s$, namely,  $ {\cal M}_ {\rho }  ^s$ and  $ {\cal P}
^sg$.

Associated with the splitting (\ref{directsum0})  $h=\tilde h+\varphi g$,  we let  \begin{equation}
\label{directsum01}\fl\qquad \widetilde P_g: S_2^s \longrightarrow \bar S_2^s\;,\quad h \longmapsto\widetilde P_g(h)=\tilde
h=h-\varphi g  =h-
\big(L_g ^{-1} (D\!R(g)h)\big)g
\end{equation} denote the projection onto $\bar S^s_2(g)=\ker D\!R(g)$ and refer to $\tilde h$ as the 
{\bf   tangential component} (or {\bf part}) of $h$ and $ \varphi g $ as the {\bf (infinitesimal)
conformal component} (or {\bf part}) of
$h$. We use this  terminology   since, generically, for $\rho  \in {\cal F} ^{s-2}$, $ {\cal M} _ \rho ^s=\{\,g\in {\cal M} ^s\mid
R(g)=
\rho
\,\}$ is a closed submanifold of
$ {\cal M} ^s$ with tangent space at $ g\in {\cal M}  _ \rho ^s$ given by $T_g {\cal M}_{\rho } \approx \ker D\!R(g)$. Thus these
  directions infinitesimally preserve the scalar curvature.  Similarly,      the space of (pointwise) conformal deformations  $
{\cal P} ^sg$ of
$g$ is a closed submanifold of
$ {\cal M} ^s$
 with tangent space $T_g( {\cal P} ^sg)\approx {\cal F} ^sg$ so that the (infinitesimal) conformal directions preserve the pointwise
conformal class of $g$. Although the conformal directions are also tangential directions (to $ {\cal P} ^sg$),  
by the tangential part of $h\in S_2^s$ we shall mean the tangential directions to $ {\cal M}_{-1} ^s$  since this is the constraint
space for the conformal Ricci flow.

\begin{example}
\label{splitRicci1}
 {\bf(Splitting of ${\rm Ric}(g) $)}  {\rm  (Compare Example~\ref{ex2}.)  As an important example, for $n\ge3$, $
s>\frac{n}{2}+3$, 
$ g\in {\cal M} ^s$, we split the Ricci tensor
$h= {\rm Ric}(g) \in S_2^{s-2}$ of $g$ according to (\ref{directsum0}) under the    assumption that 
 $L_g = (n-1)\Delta_g -R(g) $ is an isomorphism. Note  that $ s-2 >  \frac{n}{2}+1 $ so the hypothesis on $s$ of  
Proposition~\ref{non4} 
 is   satisfied. Thus  $ {\rm Ric}(g) $ splits as \begin{equation}
\label{ricsplit}  {\rm Ric}(g) =\widetilde {{\rm Ric}(g) }+\varphi g\,,
\end{equation} where  from (\ref{ricci2}), 
$ D\!R(g){\rm Ric}(g )=
  {\textstyle\frac{1}{2}} \Delta_g R(g)    - |{\rm Ric}(g) |^2   \in {\cal F} ^{s-4}$
since  $  |{\rm Ric}(g) |^2$ and  $ R(g)$ are  of class $  H^{s-2}$ so that $\Delta_g R(g) \in {\cal F} ^{s-4}$.
Then, by the ellipticity of $L_g$,
\begin{equation}
\label{splitric}  \varphi= L_g ^{-1} ( D\!R(g){\rm Ric}(g))=
 L_g ^{-1} \Big(  {\textstyle\frac{1}{2}} \Delta_g R(g)    - |{\rm Ric}(g) |^2  \Big )\in {\cal F} ^{s-2}\,,
\end{equation}so that $\varphi g\in S_2^{s-2}$.   Thus $ {\rm Ric}(g) $ splits as   \begin{equation}
\label{riccidec}   {\rm Ric}(g) =\widetilde {{\rm Ric}(g)} +\Big( L_g ^{-1} \Big(  {\textstyle\frac{1}{2}} \Delta_g R(g)    - |{\rm
Ric}(g) |^2  \Big ) \Big)g\,,
\end{equation}  where
\begin{equation}
\label{LABEL} \widetilde P_g({\rm Ric}(g) )= \widetilde {{\rm Ric}(g) } ={\rm Ric}(g) - \varphi g\in\bar S^{s-2}_2\,. 
\end{equation}  
Note that each summand
of  $ {\rm Ric}(g) \in S_2^{s-2}$   is also  of class $H^{s-2}$ as required by the splitting.\hfill$\rule{2mm}{3mm}$  
 }
 \end{example}

We remark that   $ {\rm Ric}(g) $ splits non-orthogonally into two non-trivial parts (see
also Figure~1 below). This is in contrast with the $L_2$-orthogonal splitting (\ref{eqn:splitting4})   where   $ {\rm
Ric}(g) \in {\rm range}\,D\!R^*(g)$ and thus does not split non-trivially (see Example~\ref{ex2}). 

 In the splitting  (\ref{ricsplit}), since
$ {\rm Ric}(g) 
$ is
$H^{s-2}$,  each term in the splitting is, as expected,  also $H^{s-2}$. However, relative to $g\in {\cal M} ^s$, this splitting
loses two derivatives, which is also expected since $ {\rm Ric}(\,\cdot\,) $ is a second-order operator. What is completely
unexpected is the following.

\begin{example} {\bf(Splitting of ${\rm Ric}(g) $ when $R(g)=-1$)} 
\label{splittingexample}
{\rm 
Suppose $g\in {\cal M}_{-1} ^s$, $ s>\frac{n}{2}+3 $.
Then $ L_g = (n-1)\Delta_g +1$ is an isomorphism   and thus the necessary  curvature  condition on $g$ for the splitting is
satisfied. Thus 
  we can    split $ {\rm Ric}(g) \in S^{s-2}_2$ as above.

When the constraint is not satisfied,  the 
  argument   $ D\!R(g){\rm Ric}(g )={\textstyle\frac{1}{2}} \Delta_g R(g)    - |{\rm Ric}(g) |^2  $ of $L_g ^{-1} $ 
in (\ref{riccidec})
consists of two  terms, the first of which  involves fourth order derivatives of $g$ and the second of which involves only  second
order derivatives of
$g$.  Thus the loss of derivatives comes from the first  term  involving the
fourth order derivatives of $g$. 

Now suppose that $g$ satisfies the constraint equation $ R(g)=-1$. Then the  term $\Delta_g R(g)
=0$  involving  the fourth order derivatives of $g$ drops out
  and thus the remaining  argument $- |{\rm Ric}(g) |^2 $ of $L_g ^{-1} $ is $H^{s-2}$ rather than  $H^{s-4}$ as before. 
Thus   
from the ellipticity of $L_g$,
 $ \varphi= L_g ^{-1} ( D\!R(g){\rm Ric}(g))=
- L_g ^{-1} (    |{\rm Ric}(g) |^2  )\in {\cal F} ^{s }
$
 rather than ${\cal F} ^{s-2}$ as before 
and
 $\varphi g\in  {\cal F} ^{s }g$  rather than ${\cal F} ^{s-2 }g$ as before.
Thus
 (\ref{riccidec}) reduces to 
\begin{equation}
\label{riccidec1}   {\rm Ric}(g) =\widetilde { {\rm Ric}(g)  }-( L_g ^{-1}  (       |{\rm
Ric}(g) |^2   )) g\,.
\end{equation} 
 Thus,
remarkably,  if the constraint $R(g)=-1$ is satisfied,
 $
 \varphi g=-(L_g ^{-1}  (       |{\rm Ric}(g)
|^2   ) )g=  {\rm Ric}(g) -\widetilde {{\rm Ric}(g) }= (I-\widetilde P_g) {\rm Ric}(g)\in {\cal F} ^sg
$
  does not lose derivatives relative to $g\in  {\cal M} ^s$.

Note also that since $g\in {\cal M} ^s$, we do not expect that $\varphi g$ can in  general be any smoother than $H^s$ as it is in
this case.
\hfill$\rule{2mm}{3mm}$ 
}
\end{example}

It is important to remark that this improvement in differentiability of 
 ${\rm Ric}(g) -\widetilde{ {\rm Ric}(g)}$
only takes place in the special case  when  $ \Delta_gR(g)$  in (\ref{splitric})  is of  class 
 $H^{s-2}$ (instead of
$H^{s-4}$), as occurs for example when 
 the constraint equation  $R(g)=-1$ is satisfied.   
 We   also note  that the individual terms $ {\rm Ric}(g) $ and 
$\widetilde {\rm Ric}(g)= {\rm Ric}(g) - \varphi g  $ are still $H^{s-2}$ as before. It is only the difference 
$ \varphi g={\rm Ric}(g) -\widetilde {{\rm Ric}(g)}$ that is $H^s$ rather than $H^{s-2}$.

Thus, summarizing, if  we compare  the map $ g \mapsto  (I-P_g) {\rm Ric}(g) = {\rm Ric}(g) -\widetilde{{\rm Ric}(g)}
$ when the constraint equation is not satisfied with when it is satisfied, we have\begin{eqnarray} \fl
  {\cal M} ^s
\longrightarrow S^{s-2}_2\;,\quad
\label{mapa} g \longmapsto   {\rm Ric}(g)- \widetilde{{\rm Ric}(g)}= 
\left( L_g ^{-1} \left(  {\textstyle\frac{1}{2}} \Delta_g R(g)    - |{\rm Ric}(g) |^2  \right)\right)g=\varphi g       \\ \fl
  {\cal M}_{-1}    ^s \longrightarrow S^s_2\;,\quad
\label{mapb} \;\,g \longmapsto  {\rm Ric}(g)- \widetilde{{\rm Ric}(g)}= -
\left( L_g ^{-1} (      |{\rm Ric}(g) |^2   )\right)g=\varphi g\,.
\end{eqnarray} Thus, in general, if the constraint equation is not satisfied, 
  (\ref{mapa}) 
maps $H^s$ Riemannian metrics to $H^{s -2}$ tensors  with loss of derivatives whereas
if the constraint $R(g)=-1$     is satisfied, 
(\ref{mapb}) 
maps $H^s$ Riemannian metrics to $H^{s }$ tensors   $\varphi g\in {\cal F} ^sg\subset S_2^s$ without loss of derivatives.

This improvement when $R(g)=-1$  is the basis for the fact that the 
  conformal Ricci flow does not lose derivatives  relative to the unnormalized Ricci flow. This discussion can be rephrased in terms
of the conformal pressure as follows.

\begin{example} {\bf(Splitting of $ -2({\rm Ric}(g) + {\textstyle\frac{1}{n}} g)$ when $R(g)=-1$)} 
\label{splittingexample1}
{\rm
Continuing with the assumptions of Example~\ref{splittingexample}, we split $ -2({\rm Ric}(g) + {\textstyle\frac{1}{n}} g)$
 when $R(g)=-1$. From  (\ref{global}) and     (\ref{riccidec1}),
 \begin{eqnarray}
\nonumber 
-  2{  ({\rm Ric}(g) +  \textstyle\frac{1}{n}} g)   & = &  -  2 {   \widetilde{{\rm
Ric}(g) }  
+ {  \big(2   L_g ^{-1} (      |{\rm Ric}(g) |^2   )  -  \textstyle\frac{2}{n}}
 \big) g        
 }
 \nonumber \\ 
\label{spliteqn}  & = & - 2 {   \widetilde{{\rm
Ric}(g) } 
+pg    \,,    
 }
\end{eqnarray}  where 
  $  -2    \widetilde{{\rm Ric}(g)} $ is the tangential component,
      the   pressure
term  
$ pg$ is the
  (infinitesimal)   conformal    component,
and
 $p=  {   2  ( L_g ^{-1} (      |{\rm Ric}(g) |^2   ) 
-  \textstyle\frac{1}{n}})
   $.
  \hfill$\rule{2mm}{3mm}$ 
}
\end{example}

By rearranging (\ref{spliteqn}),  
 we can   write the splitting of 
  $ -2  {\rm Ric}(g) $ as \begin{equation}
\label{splittingof} - 2 {\rm Ric}(g) = -  2 \widetilde {{\rm Ric}(g)}+
(p+ {\textstyle\frac{2}{n}} )g\,,
\end{equation} with   
  tangential component 
$- 2 \widetilde {{\rm Ric}(g)}$ and  conformal component  
$
 (p+ {\textstyle\frac{2}{n}} )g 
$. Also, for the purpose of Figure~1, we re-write (\ref{spliteqn})   as  \begin{equation}
\label{writtenas} -  2 {   \widetilde{{\rm
Ric}(g) }  = -  2{  ({\rm Ric}(g) +  \textstyle\frac{1}{n}} g) -pg       
 }
\end{equation}
Thus, from   (\ref{writtenas}),   the    conformal Ricci equation can   be written
in terms of the projection map $\widetilde P_g$  as follows,
\begin{eqnarray}\nonumber  0&=&  \frac{\partial g}{\partial t} +   2{  ({\rm Ric}(g) +  \textstyle\frac{1}{n}} g) +pg  \\&=&
\frac{\partial g}{\partial t}+ 2 \widetilde{\rm Ric}(g) \nonumber  \\ & = & \frac{\partial g}{\partial t} +2\widetilde P_g({\rm
Ric}(g)) \,. 
\label{red430} 
\end{eqnarray} 
 Since   $\widetilde P_g(g)=0$,   (\ref{red430})   can also be written as  
$ \frac{\partial g}{\partial t}+ 2\widetilde P_g( 
  {  {\rm Ric}(g) +  \textstyle\frac{1}{n}} g)=0$.
  
The ``vector"  relationships  (\ref{spliteqn}), (\ref{splittingof}), and  (\ref{writtenas}) are  depicted in Figure~1 (see page 51).
Note in particular that $-2 {\rm Ric}(g) $ is orthogonal to $T_g {\cal M}_{-1} $ and that 
on the right hand side of (\ref{writtenas}), 
   the two summands $  -  2{  ({\rm Ric}(g) +  \textstyle\frac{1}{n}} g)  \in S_2^T\!(g)  $
and 
$  -pg\in {\cal F}  g\subset S_2 $ are (pointwise) orthogonal to each other
 since $R(g)=-1$ (see (\ref{onsince})).

Thus we have the final geometrical picture. The non-orthogonal  splitting  (\ref{directsum0})    decomposes 
the nonlinear restoring force
$ -2{  ({\rm Ric}(g) +  \textstyle\frac{1}{n}} g)$
into a conformal component $pg$ and   a tangential component
$-2 \widetilde{ {\rm Ric}(g) }$ 
(see (\ref{spliteqn})).
Equivalently, (\ref{writtenas}) resolves the inertial ``vector" $\frac{\partial g}{\partial t} =-2 \widetilde{ {\rm Ric}(g) }=-  
2{  ({\rm Ric}(g) + 
\textstyle\frac{1}{n}} g) -pg \in \ker D\!R(g)
$ into   orthogonal components where the
  constraint force $-pg\in {\cal F}   g\approx T_g( {\cal P}  g)$ acts  (pointwise) orthogonally to the nonlinear restoring force
$ -2{  ({\rm Ric}(g) +  \textstyle\frac{1}{n}} g)\in S_2^T(g) \approx T_g^\perp ( {\cal P}  g)$ and counter-balances the conformal
component  $pg$ of this force (in the non-orthogonal splitting).
   The scalar curvature of the resulting 
    conformal Ricci flow is therefore    preserved. The resolution of $ \dot g=\frac{\partial g}{\partial t}$ into these two
components is shown in  Figure~1 (see page 51). 

We also 
remark that 
  in Figure~1  we have represented the 
 ``thick" spaces  $ {\cal M}_{-1} $ and $ {\cal M} _ {d\mu_g} $,  which are codimension-$  C^\infty(M, \mbox{\boldmath$R$} )$ in $
{\cal M} $,  as 2-dimensional spaces, whereas we have represented the ``thin spaces"  $ {\cal P} g$ and $T_g ^\perp{\cal M}_{-1}
\approx {\rm range}\,D\!R(g)^*= D\!R(g)^*( {\cal F} )$, which have the dimensionality of $ C^\infty(M, \mbox{\boldmath$R$} )= {\cal
F} $,  as  1-dimensional spaces. These thin spaces are represented by  thin solid lines whereas the dotted lines are construction
lines.

The conformal pressure term  $-p g$     measures the ``force" or pressure that the metric
experiences by being constrained to lie in
$   {\cal M}_{-1}
$.  If there were no pressure term, the metric would follow the flow lines of the   equation
$ \frac{\partial g}{\partial t}+  2{  ({\rm Ric}(g) +  \textstyle\frac{1}{n}} g) =0 $
 and the
resulting motion would be neither volume nor scalar curvature preserving.

From the point of view of the non-orthogonal splitting (\ref{spliteqn}), the   nonlinear restoring term $- 2( {\rm Ric}(g) +
{\textstyle\frac{1}{n}} g)$ has a component
$-2 
\widetilde{ {\rm Ric}(g) }$ in the tangential direction of $ {\cal M}_{-1} $  and a transverse  component  $  pg$ in the
infinitesimal conformal direction
$T_g( {\cal P} g)\approx {\cal F} g$. The
  conformal pressure term $-p g$ then acts oppositely   
to provide the necessary   counter-force 
to  balance  the transversal component  $  pg$  of $- 2( {\rm Ric}(g) + {\textstyle\frac{1}{n}} g)$, thereby keeping the
flow in $ {\cal M}_{-1} $.

Note that  $- 2( {\rm Ric}(g) +
{\textstyle\frac{1}{n}} g)\in S_2^T(g)\approx T^\perp_g( {\cal P} g)$ and thus has no $ {\cal P} g$ component when split
orthogonally by (\ref{geometricallyas})  but when split non-orthogonally its $ {\cal P} g$ component is $pg$.
This is analogous to the fact that $-2 {\rm Ric}(g) \in {\rm range}\,D\!R(g)^*\approx T_g {\cal M}_{-1} ^\perp$ and thus has no 
$T_g {\cal M}_{-1} $ component when split orthogonally by (\ref{eqn:splitting4})  
 but when split non-orthogonally its $T_g {\cal M}_{-1} $  component is $-2\widetilde{ {\rm Ric}(g) }$.

Note that if we write (\ref{red430}) as $ 
  \frac{\partial g}{\partial t}=-  2\widetilde P_g({\rm Ric}(g) ) = -  2{  {\rm Ric}(g) -  \textstyle\frac{2}{n}} g   
-  pg$, then the right hand side has zero net conformal component, which can be viewed in two different ways. Either 
$ pg$ is the conformal part of $  -  2({  {\rm Ric}(g) +  \textstyle\frac{1}{n}}  g) $ or 
$  (p+\textstyle\frac{2}{n} ) g $  is the conformal part of $-   2  {\rm Ric}(g)$, as shown in Figure~1.
Since we  view the system without the constraint as being $ \frac{\partial g}{\partial t}+  2{  ({\rm Ric}(g) + 
\textstyle\frac{1}{n}} g) =0 $ (and not the unnormalized Ricci equation $ \frac{\partial g}{\partial t}+  2  {\rm Ric}(g)  
  =0 $), we adhere to  the first interpretation.

If we consider the 
 conformal pressure as a map, also denoted $p$, 
\begin{equation}
\label{LABEL} p: {\cal M}_{-1} ^s \longrightarrow  {\cal F}  ^s\;,\qquad g \longmapsto p(g)=
2 L_g ^{-1} \big(
|{\rm Ric} (g)|^2_g \big)-  {\textstyle\frac{2}{n}}\,,
\end{equation}    then  $p$ is a  smooth map which is  everywhere defined and   does not lose derivatives. 
 
 Summarizing, 
 \begin{proposition}{\bf(The constraint force for the conformal Ricci flow does not lose derivatives)}
\label{summarizing} If $g \in {\cal M}_{-1} ^s$, $s>\frac{n}{2}+3$,  then $ 2({\rm
Ric}(g)+ {\textstyle\frac{1}{n}} g)
$ and $ 2\widetilde P_g( {\rm Ric}(g)+ {\textstyle\frac{1}{n}} g )$ are    
$H^{s-2}$ but their difference  
\begin{eqnarray}\nonumber \fl\qquad\qquad   2(I-\widetilde P_g) ({\rm Ric}(g) +  {\textstyle\frac{1}{n}} g)  & = & 2( {\rm Ric}(g) +
{\textstyle\frac{1}{n}} g-\widetilde P_g( {\rm Ric}(g) ))       =  -  pg \nonumber  
\end{eqnarray} 
    is $H^s$ and not merely $H^{s-2}$, the same differentiability class of $g$. 
The conformal pressure 
\begin{equation}
\label{LABEL} p
 :{\cal M}_{-1}  ^s\longrightarrow {\cal F} ^s   \;,\quad g \longmapsto  2\,L_g ^{-1} (| {\rm Ric}(g)
|^2) - {\textstyle\frac{2}{n}} 
\end{equation} is  a smooth map on the indicated spaces and the associated conformal pressure constraint force
\begin{equation}
\label{LABEL}  {\cal M}_{-1}  ^s\longrightarrow   S_2 ^s \;,\quad g \longmapsto- pg
\end{equation}   is a smooth   vector field on $ {\cal
M}^s_{-1} $ in the conformal direction.

The reduced conformal Ricci equation can be written as \begin{equation}
\label{proj3}  \frac{\partial g}{\partial t}+2 \widetilde  P_g({\rm
Ric}(g))   =0\,,
\end{equation} where $\widetilde P_g\,\colon S^s_2 \rightarrow \bar S_2^s(g)=\ker D\!R(g)$ is the projection onto the $\bar S_2^s$
component of the non-orthogonal   splitting $S_2^s=\bar S_2^s\oplus {\cal F} ^sg\approx T_g   {\cal M}_{-1}
^s\oplus T_g( {\cal P} ^sg)$.
 \hfill$\rule{2mm}{3mm}$ 
 \end{proposition} 

Lastly, we remark on the general structure of the conformal Ricci flow equations.
One might try to   use the $L_2$-orthogonal projection
(\ref{eqn:splitting4}) instead of  the non-orthogonal projection of (\ref{proj3})  
 in order to constrain the Ricci flow to $ {\cal M}_{-1} $. Thus, for example, one might consider the system
 \begin{eqnarray}\label{parabolic88} \frac{\partial g}{\partial t}+2\,{\rm
Ric}(g) &=&-2D\!R(g)^*\phi              \\
\label{elliptic88}  \;\;\qquad\quad   R(g)&=&-1   
\end{eqnarray} where the function $\phi$ is to be solved   so as to maintain the   flow   in $ {\cal M}_{-1} $ (in other words, one
might consider the equation 
$ \frac{\partial g}{\partial t}+2 \bar P_g({\rm
Ric}(g))=0          $ where $\bar P_g$ is $L_2$-orthogonal projection onto $\ker D\!R(g)\approx T_g {\cal M}_{-1} $). 
However,   the Ricci tensor is orthogonal to $ {\cal M}_{-1} $, $ {\rm Ric}(g) \in {\rm range}\,D\!R(g)^*\approx T_g^\perp {\cal
M}_{-1}  $, so that   $\bar P_g({\rm
Ric}(g))=0        $. Thus one cannot non-trivially constrain the Ricci flow to $ {\cal M}_{-1} $ by orthogonal 
projection of
$ {\rm Ric}(g) $. To put this another way, for every $g \in {\cal M}_{-1} $, if $\phi=1$, then $-2 {\rm Ric}(g) -2  D\!R(g)^*(1)=
-2 {\rm Ric}(g) - ( - 2{\rm Ric}(g) )=0$ so that  every point of $ {\cal M}_{-1} $ is a critical point for  
(\ref{parabolic88}--\ref{elliptic88}).
 Thus to get a non-trivial {\it conformal}  Ricci flow in $ {\cal M}_{-1} $, we have had to use  a non-orthogonal  
projection   onto $T_g {\cal M}_{-1} $.

 \section[A variational approach to the conformal Ricci flow equation]
{A variational approach to the conformal Ricci flow equations} In this section we discuss a variational approach to the
conformal Ricci flow  equation.
 
 Let
$ v_g= {\rm vol}(M,g)
$ denote the volume of the Riemannian manifold $(M,g)$.   The 
 {\it Yamabe functional} is defined by 
\[
\label{yamabe1}\fl\qquad\qquad
Y : {\cal M} \longrightarrow \mbox{\boldmath$R$}  \;,\quad g \longmapsto \frac{\int_MR(g) d\mu _g}{(\int_M  d\mu
_g)^{(n-2)/n}}= v_g^{2/n}\frac{R _{\rm total} (g)}{v_g }=  v_g   ^{2/n} \bar R  (g)
\]where $ \bar R  (g)= v_g ^{-1}  R_{\rm total} (g)$ is the {\bf volume-averaged total scalar curvature}. Thus the Yamabe
functional
$Y$ is  a {\bf volume-normalized  total scalar curvature functional}  weighted to be invariant under homothetic transformations
of
$g$, i.e., if
$c\in
\mbox{\boldmath$R$} ^+ $, then
\[ Y(cg)=\frac{\int_MR(cg) d\mu _{cg}}{ 
v_{cg} ^{(n-2)/n}}=\frac{c ^{-1} c ^{n/2}}{(c^{n/2})^{(n-2)/n}} \frac{\int_MR( g) d\mu _{ g}}{ 
v_{ g} ^{(n-2)/n}}=Y(g) 
\]
 
Let  $ {\cal G} $ denote the
 natural   weak
  $L_2$-Riemannian metric   (see $(\ref{naturalmetric})$) on  $ {\cal M} $ and let 
  $ {\cal G} _{-1}\equiv{\cal G} _{ {\cal M}_{-1} }$ denote the intrinsic  Riemannian metric naturally induced  on  the submanifold $
{\cal M}_{-1} $ from  
$ {\cal G}
$. Thus  for $g\in {\cal M}_{-1} $, 
\[ {\cal G} _{-1}(g):T_g {\cal M}_{-1}
\times T_g {\cal M}_{-1}
\longrightarrow
\mbox{\boldmath$R$} \,,\]
where if $\bar h, \bar k\in \ker D\!R(g)\approx T_g {\cal M}_{-1} $, then 
 ${\cal G} _{-1}(g)(\bar h,\bar k)= {\cal G}  (g)(\bar h,\bar k)$.
In the next proposition, we  first compute the intrinsic gradient of the
  Yamabe functional restricted to $  {\cal M}_{-1} $ and see that this  does not give the conformal Ricci flow equation. However, by
using the non-orthogonal splitting
 $(\ref{directsum0})$  and the non-orthogonal projection 
    (\ref{directsum01}),  
we  can formulate the conformal Ricci flow  equations in a  
 gradient-like manner, which we refer to as a {\bf quasi-gradient}. 
\begin{proposition}[A quasi-gradient form of  the conformal Ricci flow equations]
\label{vari} The gradient of the Yamabe functional 
\[
Y : {\cal M} \longrightarrow \mbox{\boldmath$R$}  \;,\quad g \longmapsto  v_g^{(2-n)/n  }R_{\rm total} (g)=  v_g   ^{2/n} \bar
R  (g)
\]
in the natural $L_2$-Riemannian metric  $ {\cal G} $   on $ {\cal M} $ is given
by
 \begin{eqnarray}   \nonumber\fl\quad 
 {\rm grad}\,Y:{\cal M} \longrightarrow S_2\;, \quad  g \longmapsto ({\rm grad}\,Y)(g) =- v_g^{(2-n)/n  } \big (     
{\rm Ein}(g)
  + 
 {\textstyle\frac{n-2}{2n}}    \bar R  (g)g   \big)  \nonumber 
\\\qquad\qquad\quad= - v_g^{(2-n)/n  }\Big(   \big(  
{\rm Ric}(g)
  - {\textstyle\frac{1}{n}} \bar R  (g) g \big)+ 
 {\textstyle\frac{1}{2}} (  \bar R  (g) - R(g) \big )g\Big) \,.
\nonumber 
\end{eqnarray}  
 When
restricted to $ {\cal M}_{-1} $,
   $ {\rm grad}\,Y$  simplifies to
   \begin{eqnarray}\fl \;
\nonumber 
 ({\rm grad}\,Y)  _{|{\cal M}_{-1} }: {\cal M}_{-1} \longrightarrow S_2\;,\quad g \longmapsto 
({\rm grad}\,Y)  _{|{\cal M}_{-1} }(g)\\\qquad\qquad\qquad =
- v_g^{(2-n)/n  }\big(     {\rm Ric}(g)
  + {\textstyle\frac{1}{n}}   g\big)  \nonumber=- v_g^{(2-n)/n  }      {\rm Ric}^T\!(g)\,.
 \nonumber 
\end{eqnarray} 

 The Yamabe functional $Y$  when restricted to $ {\cal M}_{-1} $ simplifies to
\[
  Y_{|{\cal M}_{-1} }: {\cal M}_{-1} \longrightarrow  \mbox{\boldmath$R$} \;,\qquad g \longmapsto Y_{|{\cal M}_{-1}
}(g)=Y(g)=-v_g^{2/n}\,.
\] 
 Let 
 \[{\rm grad}_   {-1  } Y   _{|{\cal
M}_{-1} }: {\cal M}_{-1} \rightarrow T {\cal M}_{-1}\] 
  denote the gradient of 
the restricted functional $Y   _{|{\cal
M}_{-1} } $   in the   Riemannian manifold    $( {\cal M}_{-1}, {\cal G}  _{-1} ) $.
Then with  $T_g {\cal M}_{-1}$ identified with $\ker D\!R(g)$,
  for $ g\in {\cal M}_{-1} $,
\[\fl\qquad
\label{LABEL}  ({\rm grad}_{ {-1} }   Y   _{|{\cal
M}_{-1} })(g) = \bar P_g\big( ({\rm grad}\,Y)  _{|{\cal M}_{-1} }(g)\big)=-   {\textstyle\frac{1}{n}}v_g^{(2-n)/n  }\bar
P_g (     
   g ) \in \ker D\!R(g)\,,
\] where   
$\bar P_g\,\colon S_2 \rightarrow \ker D\!R(g)$  is  
  $L_2$-orthogonal   projection   (see  $(\ref{projectionortho})$)
  and  where $\bar P_g(g)= \bar g= g -  D\!R(g)^*\big( (D\!R(g) D\!R(g)^* ) ^{-1} ( 
 1
 )   
\big)$  
(see $(\ref{splitortho33})$).

Using the non-orthogonal   projection 
$ \widetilde P_g:S_2 \rightarrow \ker D\!R(g)$ (see $(\ref{directsum01})$),   
the non-orthogonal   projection of  $({\rm grad}\,Y)  _{|{\cal M}_{-1} }$ to $T_g {\cal M}_{-1} $ 
 is given by 
  \begin{eqnarray}\widetilde  P_g \big(({\rm grad}\,Y)  _{|{\cal M}_{-1} }(g)\big) & = & - v_g^{(2-n)/n}  \big(     {\rm Ric}(g)
  + {\textstyle\frac{1}{n}}   g+ {\textstyle\frac{1}{2}}p g\big) \nonumber         \\   & = &- v_g^{(2-n)/n}  \big(     {\rm Ric}(g)
  +  L_g ^{-1}  (
|{\rm Ric} (g)|^2_g  )\big)  \,,\nonumber  
\end{eqnarray} 
 where
$
  p=  2 L_g ^{-1}  (
|{\rm Ric} (g)|^2_g  )-  {\textstyle\frac{2}{n}}
 $  is  the conformal pressure. 
Thus the reduced  conformal Ricci flow equation  can   be written in the  quasi-gradient form,  
 \begin{eqnarray}   \frac{\partial g}{\partial t} 
&=&2 v_g^{(n-2)/n}\widetilde P_g \left(({\rm grad}\,Y)  _{|{\cal M}_{-1} }(g) \right)  
 \label{red4344}  \,.
\end{eqnarray} 
 \end{proposition} {\bf Proof:} First we recall  the   facts that
$({\rm grad}\,\, {\rm vol}  )(g)= {\textstyle\frac{1}{2}}\, g$  and  $({\rm grad}\,\, R_{\rm total})(g)=-{\rm Ein}(g)$.
 That $({\rm grad}\,\, {\rm vol}  )(g)= {\textstyle\frac{1}{2}} 
g$ follows from  the fact that $({\rm grad}\,\, {\rm vol}  )(g) $ is defined as  the ``vector" in $S_2$ such that  for all $h\in
S_2\,$, \[
(d \,{\rm vol})(g)h=\int_M ({\rm grad}\,\, {\rm vol}  )(g)\cdot h\,d \mu _g\,.
\]Thus, since
 \begin{eqnarray}\fl (d \,{\rm vol})(g)h= \int_M D(d   \mu _g)h= {\textstyle\frac{1}{2}}  \int_M{\rm tr}_g h
\,d
\mu _g          
\nonumber = {\textstyle\frac{1}{2}}  \int_M g \cdot h \,d \mu
_g =\int_M ({\rm grad}\,\, {\rm vol}  )(g)\cdot h\,d \mu _g 
\end{eqnarray}(where  ``$\,\cdot\,$" is the metric contraction)   for all $h\in S_2$, it follows that  $ ({\rm grad}\,\, {\rm
vol})(g)    ={\textstyle\frac{1}{2}}   g$.

Similarly,  
since for all $h\in S_2$,  \begin{eqnarray} (dR_{\rm total})(g)  h &
= &
\int_M D(R(g)d  
\mu _g)h=
 \int_M
D \!R(g )h \,d  
\mu _g + \int_M
R(g)\,D(d  
\mu _g)h\nonumber \\ \nonumber &= &
\int_M
( \Delta_g{\rm
tr}_gh+\delta_g\delta_gh- {\rm Ric}(g)  \cdot h)
+
{\textstyle\frac{1}{2}}  \int_MR(g)\,{\rm tr}_g h \,d \mu _g          
\nonumber   \\  & = & \int_M
( - {\rm Ric}(g)  
+
{\textstyle\frac{1}{2}}  R(g) g )\cdot h \,d \mu _g 
\nonumber \\ \nonumber &=&- \int_M
  {\rm Ein}(g)  
 \cdot h \,d \mu _g   
= \int_M({\rm grad}\,R_{\rm total})(g)\cdot h \,d \mu _g   \,, \nonumber 
\end{eqnarray} we have
 $({\rm grad}\,\, R_{\rm total})(g)=-{\rm Ein}(g)$. With these two basic gradients in hand, we find
 \begin{eqnarray}\fl \qquad \nonumber 
 ({\rm grad}\,Y)(g)  
&=&  ({\rm grad}\,
{\rm vol}  ^{(2-n)/n  }R_{\rm total} )(g) \\\nonumber 
&= &({\rm grad}\,
{\rm vol}  ^{(2-n)/n  })(g)R_{\rm total}   (g) +  
v_g ^{(2-n)/n  }({\rm grad}\,R_{\rm total}  )(g) \\ \nonumber 
&=& {\textstyle\frac{2-n}{n}}   ( 
v_g ^{(2-n)/n  })v_g ^{-1}({\textstyle\frac{1}{2}}\, g)R_{\rm total}  (g) - 
v_g ^{(2-n)/n  } {\rm Ein}(g)   
\nonumber  \\ \nonumber 
&=&    
- v_g^{(2-n)/n  } \big(     
{\rm Ein}(g)
  + 
 {\textstyle\frac{n-2}{2n}}    \bar R   (g)g  \big)  
\\&=& - v_g^{(2-n)/n  } \big(     
({\rm Ric}(g)- {\textstyle\frac{1}{2}} R(g)g)
  + 
 ({\textstyle\frac{1}{2 }}- {\textstyle\frac{1}{n}}  )   \bar R   (g)g  \big)  \nonumber    \\ \nonumber  &   =&  - v_g^{(2-n)/n 
}\big(  
 (   {\rm Ric}(g)
  - {\textstyle\frac{1}{n}} \bar R   (g) g  )+ 
 {\textstyle\frac{1}{2}}  (  \bar R   (g) - R(g)   )g\big) \,.
\end{eqnarray}

If $ g\in {\cal M}_{-1} $, then  $\bar R   (g)
  =  v_g ^{-1} R_{\rm total} (g)  =-1=R(g)$ and thus   $ {\rm grad}\,Y$ restricted to $ {\cal M}_{-1} $ simplifies to
   \begin{eqnarray} \fl\qquad\qquad
 ({\rm grad}\,Y)  _{|{\cal M}_{-1} }(g)=
- v_g^{(2-n)/n  } (     {\rm Ric}(g)
  + {\textstyle\frac{1}{n}}   g ) 
=- v_g^{(2-n)/n  }     {\rm Ric}^T\!(g)\,.
  \label{simp1}
\end{eqnarray} 

The   Yamabe functional restricted to $ {\cal M}_{-1} $  is given by $ Y   _{|{\cal
M}_{-1} } (g) =v_g^{(2-n)/n  }R^T  (g)= v_g   ^{2/n} \bar R^T  (g)=-v_g^{ 2  /n  } $.
Denoting its gradient      in the weak Riemannian manifold $( {\cal M}_{-1} , {\cal G}_{-1} )$ as ${\rm grad}_{  {-1} }$, we have,
using (\ref{simp1}) and $\bar P_g (   {\rm Ric}(g) )=0$, 
\[
\label{LABEL}  ({\rm grad}_{  {-1} }   Y   _{|{\cal
M}_{-1} })(g) = \bar P_g\Big( ({\rm grad}\,Y)  _{|{\cal M}_{-1} }(g)\Big)=-   {\textstyle\frac{1}{n}}v_g^{(2-n)/n  }\bar
P_g (     
   g ) \,,
\] 
  where $\bar P_g(g)= \bar g= g -  D\!R(g)^*\big( (D\!R(g) D\!R(g)^* ) ^{-1} ( 
 1
 )   
\big)$ as given in 
 $(\ref{splitortho3})$.  

Using  $\widetilde P_g(g)=0$ and   (\ref{spliteqn}),
the non-orthogonal   projection of  $({\rm grad}\,Y)  _{|{\cal M}_{-1} }$ to $T_g {\cal M}_{-1} $ 
 is  
\begin{eqnarray} \fl\qquad\widetilde P_g \big(({\rm grad}\,Y)  _{|{\cal M}_{-1} }(g)\big) & = & - v_g^{(2-n)/n}\widetilde
P_g\big(     {\rm Ric}(g)
  + {\textstyle\frac{1}{n}}   g\big)       \nonumber    =  - v_g^{(2-n)/n} \widetilde{  {  {\rm Ric}(g)}} \nonumber \\   &
= &
 - v_g^{(2-n)/n}  \big(     {\rm Ric}(g)
  + {\textstyle\frac{1}{n}}   g+ {\textstyle\frac{1}{2}}p g\big)
\nonumber 
 \\   & = &- v_g^{(2-n)/n}  \big(     {\rm Ric}(g)
  +  L_g ^{-1}  (
|{\rm Ric} (g)|^2_g  )\big) \,, \nonumber  
\end{eqnarray} 
  where the conformal pressure 
$
  p=  2 L_g ^{-1}  (
|{\rm Ric} (g)|^2_g  )-  {\textstyle\frac{2}{n}}
\,$.
 Thus
  \[
    2\big(   {\rm Ric}(g)
  + {\textstyle\frac{1}{n}}   g\big)+ p g =
- 2 v_g^{(n-2)/n}\widetilde P_g (({\rm grad}\,Y)  _{|{\cal M}_{-1} }(g) )  
\]
and so the reduced conformal Ricci equation can be written in the quasi-gradient form 
 \begin{eqnarray}\nonumber   \frac{\partial g}{\partial t}  &=&- 2{  ({\rm Ric}(g) +  \textstyle\frac{1}{n}} g)   
- pg
\nonumber  =
  2 v_g^{(n-2)/n}\widetilde P_g (({\rm grad}\,Y)  _{|{\cal M}_{-1} }(g) )  \,.
\nonumber 
\end{eqnarray} 
 \hfill$\rule{2mm}{3mm}$ \\ 
\medskip

Comparing  the two projections of $  ({\rm grad}\,Y)  _{|{\cal M}_{-1} }$,   the $L_2$-orthogonal projection
\begin{eqnarray} \nonumber\fl\qquad\quad \bar P_g\Big( ({\rm grad}\,Y)  _{|{\cal M}_{-1} }(g)\Big)&=&
 ({\rm grad}_{  {-1} }   Y   _{|{\cal
M}_{-1} })(g)\\& =& -   {\textstyle\frac{1}{n}}v_g^{(2-n)/n 
}\big(g -  D\!R(g)^*\big( (D\!R(g) D\!R(g)^* ) ^{-1} ( 
 1
 )   
\big) \big) \nonumber 
 \end{eqnarray} 
 and the non-orthogonal projection,
 \begin{eqnarray} \fl\qquad\quad\widetilde  P_g \big(({\rm grad}\,Y)  _{|{\cal M}_{-1} }(g)\big) & = &
 - v_g^{(2-n)/n}   \big(     {\rm Ric}(g)
  +  L_g ^{-1}  (
|{\rm Ric} (g)|^2_g  )\big)  \,,
\nonumber 
\end{eqnarray} 
 we see that the reduced conformal Ricci flow equation is not the intrinsic  gradient    of the restricted functional $  Y   _{|{\cal
M}_{-1} }$    but is closely related to  the non-orthogonal    projection     
  of the ambient gradient
 $ ({\rm grad}\,Y)  _{|{\cal M}_{-1} }$  of $Y$ restricted to  $ {\cal M}_{-1} $.

Lastly, we remark that as
  we shall see in Proposition~\ref{volel}, for non-static flows  the volume of the conformal Ricci flow is a strictly
monotonically decreasing real-valued function of time.  Since $n\ge3$,   the coefficient $  2 v_g^{(n-2)/n}$ of
equation~(\ref{red4344})   is  also strictly monotonically decreasing    and thus can be used to 
renormalize time so that  in the new time variable,      (\ref{red4344}) can be written in the quasi-gradient form  $\dot
g=\widetilde P_g (({\rm grad}\,Y)  _{|{\cal M}_{-1} }(g) ) $.

\section{Decaying global and local volume for the conformal Ricci flow}

 In this section we show that for non-static flows,  the global and local volumes   are strictly monotonically
decreasing under the conformal Ricci flow. We then consider some applications of these results.

\begin{proposition}[Global volume decay  for the conformal Ricci flow]
\label{volel}  Let\\     $g\,\colon[0,T)\rightarrow {\cal M}_{-1}  $ be a  
conformal Ricci flow
with conformal pressure
  $p\,\colon [0,T) \rightarrow {\cal F} $.
  Then on $ [0,T)$, the volume element $d \mu _g$  satisfies
\begin{equation}
\label{voldensity}  
  \frac{\partial}{\partial t}\,d\mu _{g } = - {\textstyle\frac{n}{2}} p\, d \mu _g
\end{equation} and the  volume $ {\rm vol}(M ,g )  =
\int_M d\mu _{g  }     $ satisfies
\begin{equation}
\label{globalvol} \frac{d}
{dt}\,{\rm vol}(M ,g ) 
= - {\textstyle\frac{n}{2}} 
 \int_{M }  p\,
 d\mu _{g  } = - n
 \int_{M }  |{\rm Ric}^T \!(g)|^2_g    \,
 d\mu _{g  } 
  \,. 
\end{equation}Thus if the flow $g$ is non-static, ${\rm Ric}^T \!(g)\not=0$, $p>0$, and \begin{equation}
\label{LABEL} \frac{d}
{dt}\,{\rm vol}(M ,g ) 
<0
\end{equation}and so  the curve of volumes  ${\rm vol}(M,g  ) $   is a strictly monotonically decreasing function along the non-static flow lines
of the conformal Ricci flow.
 \end{proposition}
  {\bf Proof: }Computing  the time derivative of the volume element $d \mu _g$     yields  
\begin{eqnarray} \frac{\partial}{\partial t}\,d\mu _{g } = {\textstyle\frac{1}{2}}{\rm tr}_{g }\left(\frac{\partial g }{\partial
t}
\right)d \mu _{g }  \nonumber   = {\textstyle\frac{1}{2}}{\rm tr}_{g }\left( -2{\rm Ric}(g)- 
({\textstyle\frac{2}{n}}+p)  g \right) d\mu _{g }         \\   = \left(- R(g)-
{\textstyle\frac{n}{2}}({\textstyle\frac{2}{n}}+p)  \right)d \mu _{g } \nonumber   
  =
 (1-1  - {\textstyle\frac{n}{2}}p ) d\mu
_{g } 
   =    - {\textstyle\frac{n}{2}} p\, d \mu _g\,.
     \label{end} 
\end{eqnarray}

Integrating   (\ref{elliptic6}),  $
   ( n-1)\Delta_g  p    +
p = 
   2
(|{\rm Ric}^T \!(g)|^2_g  ) 
$,   over
$ M$,  
  using
$\int_M\Delta_g  p\,d \mu _g=0$, we find 
the following  relationship between the integrated conformal pressure and the integrated curvature,
 \begin{equation}
\label{relation5} \int_M
p \,d \mu _g= 
   2\int_M
(|{\rm Ric}^T \!(g)|^2_g  ) d \mu _g\,.
\end{equation} 
 Integrating
  (\ref{voldensity}) 
    over $ M$ then  yields 
\begin{eqnarray} \fl\quad
\nonumber  \frac{d }{dt}\,{\rm vol}(M,g )&=&\frac{d  }{dt}    \int_M d\mu _{g } 
=\int_M \frac{\partial }{\partial t}\,d\mu
_{g }=    - {\textstyle\frac{n}{2}}
\int_M  p \,d \mu _g 
=
 - n
 \int_{M }  |{\rm Ric}^T \!(g)|^2_g    \,
 d\mu _{g  } 
  \,.  
\end{eqnarray}  
Thus for non-static flows, ${\rm Ric}^T \!(g)\not=0$, $p>0$   (see Propositions~\ref{equilprop} and \ref{maxprin}), and 
 $ \frac{d }{dt}\,{\rm vol}(M,g )<0$.
  \hfill$\rule{2mm}{3mm}$ 
\bigskip

 As an application of Proposition~\ref{volel} we have the following non-existence result.
 \begin{proposition}
\label{quasi2}
There exist no non-static periodic   conformal Ricci flows.
 \end{proposition}  {\bf Proof:}  If  $g\,\colon[0,\infty)\rightarrow {\cal M}_{-1}  $ is a non-static periodic solution of the
conformal Ricci flow, then $g(t_0)=g(t_1)$ for some $t_0,t_1\in [0,\infty)$, $t_0<t_1$, and so
${\rm vol}(M,g (t_0) ) ={\rm vol}(M,g (t_1) )$, contradicting
${\rm vol}(M,g (t_0) ) >{\rm vol}(M,g (t_1) )  $ for non-static flows.
\hfill$\rule{2mm}{3mm}$ 
\bigskip

Adopting terminology from Perelman \cite{per02}  for the classical Ricci flow, we define the following.
 \begin{definition}
 [Conformal Ricci breathers and solitons]
\label{homothetic}  A conformal Ricci flow     $g\,\colon  [0,T)\rightarrow {\cal M}_{-1}  $   
 with conformal pressure
  $p\,\colon[0,T) \rightarrow {\cal F} $ and initial condition $ g(0)=g_0$  is a {\bf conformal Ricci breather}
 if for some pair $t_0, t_1\in [0,T)$, $t_0<t_1$, there exists a  $ c_1 >0$ and a diffeomorphism $f_1\in {\cal D} $  
   such that \begin{equation}
\label{LABEL}  g(t_1)= c_1f^*_1(g(t_0))\,.
\end{equation}
  For $c_1<1$, $c_1=1$, and $c_1>1$, the breather is called {\bf shrinking}, {\bf steady}, and
{\bf expanding}, respectively.

If, moreover, there exists a 
  one-parameter family of diffeomorphisms, $f\,\colon [0,T) \rightarrow 
{\cal D} $, $t \mapsto f(t)=f_t$, and a time-dependent homothetic map $c \,\colon [0,T) \rightarrow 
\mbox{\boldmath$R$} ^+$, $t \mapsto c(t)=c_t$, such that 
\begin{equation}
\label{LABEL} g(t)=c_tf^*_tg_0\,, 
\end{equation} then $g$ 
is a {\bf conformal Ricci soliton}.
  \hfill$ \rule{2mm}{3mm} $
 \end{definition}

 Thus, modulo isometry and homothety,  a conformal Ricci breather is a periodic solution and a conformal Ricci soliton is a static
solution.

Now we have the following additional  non-existence result.\begin{proposition}
\label{quasi2}
There exist no non-static       conformal Ricci breathers or conformal Ricci solitons.
 \end{proposition}  {\bf Proof:}  Let  $g\,\colon[0,T)\rightarrow {\cal M}_{-1}  $ be a    conformal Ricci breather and let
  $t_0, t_1\in [0,T)$, $t_0<t_1$,    $ c_1 >0$, and   $f_1\in {\cal D} $  
be   such that \begin{equation}
\label{LABEL}  g(t_1)= c_1f^*_1(g(t_0)) 
\end{equation}
  as given in Definition~\ref{homothetic}.
Since $g(t_0),g(t_1)\in {\cal M}_{-1} $, from (\ref{coveqn1}), \begin{eqnarray}
\nonumber
-1 & = & R(g(t_1))=R(c_1f^*_1(g(t_0)))=c^{-1} _1R(
f^*_1(g(t_0)))         \\  & = & c^{-1} _1(R( g(t_0))\circ f _1) = c^{-1} _1((-1)\circ
f _1)=-c^{-1} _1 \,, \nonumber 
\end{eqnarray} 
  so that $c_1=1$ and so $ g(t_1)=  f^*_1(g(t_0)) $. Thus  from the constraint $R(g)=-1$, any conformal Ricci breather must be
steady.

From the invariance of the volume
functional with respect to diffeomorphism, 
 ${\rm vol}(M, g(t_1) ) ={\rm vol}(M,f^*_1(g(t_0)) )={\rm vol}(M, g(t_0) ) $. Since $t_1>t_0$ 
and the volumes are not strictly monotonically decreasing, from
Proposition~\ref{volel}, $g(t)=g_0$ must be a static flow.

Lastly, since a conformal Ricci soliton 
   $g\,\colon  [0,T)\rightarrow {\cal M}_{-1}  $  is a  special case of a conformal Ricci breather (by letting $t_0=0$,   $t_1\in
(0,T)$,
$c_1=c(t_1)$, and $f_1=f(t_1)$), there can exist no non-static conformal Ricci solitons.
\hfill$\rule{2mm}{3mm}$ 
\medskip

Thus in essence the constraint $R(g )=-1$ prohibits non-trivial homothetic transformations     and the strictly
monotonically decreasing volume result for non-static conformal Ricci flows prohibits isometric transformations.

We remark that in Perelman's \cite{per02} consideration of classical Ricci flow, the proof that there are no non-static
Ricci breathers or Ricci solitons is considerably more involved, necessitating as it does  the introduction of a   non-scale and a
scale invariant functional.

Using the strong maximum principle as contained  in Proposition~\ref{maxprin}, we  can   strengthen the global volume decay  
of conformal Ricci flow to a local   
  volume decay result.  The local result is motivated by a similar result that occurs for the reduced
Hamiltonian of 
  general relativity (see Fischer-Moncrief~\cite{fm02b},   p.~5585, for
details).

 \begin{theorem}[Local volume decay  for the conformal Ricci flow]
\label{quasilocal}
 Let \\ $g\,\colon[0,T)\rightarrow {\cal M}_{-1}  $ be a non-static  
conformal Ricci flow with conformal pressure
  $p\,\colon[0,T) \rightarrow {\cal P} $. Let
$D\subseteq M$ be a   domain in $M$, possibly all of $M$,  and let $ {\rm vol}(D ,g )  =
\int_{D} d\mu _{g  }     $ denote the volume of the domain using the time-dependent metric $ g $.
  Then
on $[0,T)$, 
      \begin{equation} 
\frac{d}
{dt}\,{\rm vol}(D ,g ) 
= - {\textstyle\frac{n}{2}} 
 \int_{D }  p\,
 d\mu _{g  }  
 <0  \,.  
\label{increasing}
\end{equation} 
 Thus  ${\rm vol}(D,g  )$ is a strictly monotonically decreasing function along the non-static flow lines of the conformal Ricci
flow.
   In particular, if $D=M$,   we recover the global volume result for non-static conformal Ricci flows,
\begin{equation}
\label{LABEL} \frac{d}
{dt}\,{\rm vol}(M ,g ) 
= - {\textstyle\frac{n}{2}} 
 \int_{M }  p\,
 d\mu _{g  } 
 <0  \,.
\end{equation} 

 \end{theorem}  {\bf Proof:} From Proposition~\ref{maxprin}, if the flow $g$ is non-static, then  $p>0$ on $M$ so that the curve of
conformal pressures lies in the space of positive functions $ {\cal P} $ on $M$. Integrating (\ref{voldensity})  over
$D\subseteq M$   yields
\begin{eqnarray} 
\nonumber  \frac{d }{dt}\,{\rm vol}(D,g )&=&\frac{d  }{dt}    \int_D d\mu _{g } 
=\int_D \frac{\partial }{\partial t}\,d\mu
_{g }=    - {\textstyle\frac{n}{2}}
\int_D  p \,\,
 d\mu _{g }  <0   
\,.\hspace{0.5in}\rule{2mm}{3mm}
\end{eqnarray} 
 \medskip

Note that  unlike   the case for $D=M$,  the relationship (\ref{relation5})  between the integrated conformal pressure and the
integrated curvature is not maintained  in general 
 since when  integrating
     (\ref{elliptic6})   over $D\not=M$, from   Gauss'  theorem, the integrated Laplacian term $\int_D\Delta_g  p\,d \mu _g=-\int_D
\mbox{div}_g (\mbox{grad}_g p)d
\mu _g=-\int_{\partial D}(\mbox{grad}_g p)^\perp d
\mu _{i^*g} =
-\int_{\partial D}i^*(*dp)  
  \not=0$ need not drop out.  Here $(\mbox{grad}_g p)^\perp=g(\mbox{grad}_g p,n)=dp\cdot n $ denotes the perpendicular projection
of   
$\mbox{grad}_g p$ at the boundary $\partial D$, where $n$ denotes the unit outward normal vector field to the boundary,  
 $i\,\colon \partial D \rightarrow D$ is the inclusion
map,
$i^*g$ is the induced metric on the boundary, $d \mu _{i^*g}$ is the    volume element 
of
the boundary induced by
  $i^*g$ and the orientation of the boundary, 
  and $*$ is the usual Riemannian star operation from $k$-forms to $(n-k)$-forms.

A rather surprising consequence of this analysis is that the local volume of an arbitrary domain $D$ cannot remain constant,
even for a short time interval, unless the global flow is static.
This local result is of potential importance for applications of the conformal Ricci flow to 
geometrization of 3-manifolds since even if a singularity is forming within a domain $D$  of a manifold, the volume must still be
decreasing on that domain (see Section~\ref{end0}).

If   $g\,\colon[0,\infty)\rightarrow {\cal M}_{-1}  $ is an all-time
conformal Ricci flow, by Proposition~\ref{volel} the volume functional  ${\rm vol}(M ,g (t))$ along the flow  must converge,
\begin{equation}
\label{vconverges}   {\rm vol}(M ,g (t))~\longrightarrow~v_\infty\ge0\qquad\mbox{as}\qquad t \longrightarrow \infty\,.
\end{equation} 
If $v_\infty=0$, then we say that the conformal Ricci flow $g$ {\bf volume collapses}.

As an application of the relationship between volume and curvature as reflected in the results of this section, we have the
following curvature condition that ensures that the volume does not collapse.
 \begin{theorem}[Curvature condition for no volume collapse]
\label{nocol}
 Let   $g\,\colon[0,\infty)\rightarrow {\cal M}_{-1}  $ be an all-time non-static
conformal Ricci flow with conformal pressure
  $p\,\colon[0,\infty) \rightarrow {\cal P} $. 
For   $t\in [0,\infty)$, let 
$A(t)=
\max_{x\in M}|{\rm Ric}^T\! (g(t,x))|^2_{g(t,x)}>0$. Assume that  
$ \int_{ 0}^tA(t')dt'$
 converges as $t \rightarrow \infty$.
   Let $B=
\lim_{t \rightarrow \infty} \int_{ 0}^tA(t')dt'>0$. Then for $t_0\in [0,\infty)$, 
      \begin{equation} 
 {\rm vol}(M ,g (t)) \ge {\rm vol}(M ,g (t_0))e^{-n\int_{t_0}^tA(t')dt' }> {\rm vol}(M ,g (t_0))e^{-nB}
\,.
\end{equation} 
  In particular,     the conformal Ricci flow $g$ does not volume
collapse, i.e., $ {\rm vol}(M ,g (t))~\rightarrow~v_\infty>0$   as $ t \rightarrow \infty$.

 \end{theorem}  {\bf Proof:} From  (\ref{globalvol}) and the definition of $A(t)$,  \[
\label{globalvol2} \fl-\frac{d}
{dt}\,{\rm vol}(M ,g(t)) 
  =   n
 \int_{M }  |{\rm Ric}^T \!(g(t))|^2_{g(t)}   
 d\mu _{g(t)  }  \le
   nA(t) \int_{M } 
 d\mu _{g  } =
   n 
 A(t)
 {\rm vol}(M ,g(t) ) 
  \,. 
\]Thus
\[
\label{globalvol2} \frac{d}
{dt}\ln { {\rm vol}(M ,g(t)) } 
  \ge
 - n 
 A(t)
  \,. 
\]Integrating this differential inequality from $t_0$ to $t$ and exponentiating yields
\[
\label{LABEL} 
 {\rm vol}(M ,g (t)) \ge {\rm vol}(M ,g (t_0))e^{-n\int_{t_0}^tA(t')dt' }> {\rm vol}(M ,g (t_0))e^{-nB}
\,.
\]
 \hfill \hfill$\rule{2mm}{3mm}$ 
\medskip

This theorem gives a sufficient curvature condition to prevent volume collapse, namely, 
  if the curvature invariant $A(t)=\max_{x\in M}|{\rm Ric}^T\! (g(t,x))|^2_{g(t,x)}$  falls off sufficiently fast as $ t
\rightarrow \infty$, then  volume collapse cannot occur.
  As an example,
 if
$A(t)$ is asymptotic to $  1/t^{1+ \epsilon }$ for some  $ \epsilon >0$,   then the hypothesis of the theorem  is
met and volume collapse cannot occur.

  However,      the theorem gives only  a sufficient condition to prevent volume
collapse. Thus, for example,  a conformal Ricci flow  
 $ g \rightarrow g_e$ that   converges at any rate  as $ t \rightarrow \infty$  to an Einstein
metric
$ g_e \in {\cal M}_{-1} 
$,
${\rm Ric}  (g _e)=-{\textstyle\frac{1}{n}} g_e$,
 has ${\rm vol}(M,g) \rightarrow {\rm vol}(M,g_e) >0$ as $ t \rightarrow \infty$ and so the volume does
not collapse.
 Similarly,  if 
  $A(t)=\max_{x\in M}|{\rm Ric}^T\! (g(t,x))|^2_{g(t,x)}  =A=$ constant, or more generally, if $A(t)\ge \epsilon >0$ is  
 is bounded away from zero,  then the hypothesis of the theorem is not met and so, as far as the theorem goes, 
  the volume may or may not collapse.

\section{Locally homogeneous    conformal Ricci flows}
\label{lochom}
 \setcounter{equation}{0}
 In  \cite{ij92}, Isenberg and Jackson consider  locally homogenous solutions to the 
classical Ricci flow equation.
Local homogeneity implies that all scalar geometric
quantities  such as  the scalar curvature and  the Riemann and  Ricci curvature norm         are
spatially constant. 
 We take advantage of this fact  to use the scalar curvature to   rescale in space 
 and   reparameterize in time  locally homogenous   classical Ricci  flows  to give   locally homogenous  
  conformal Ricci flows. Thus  the locally homogenous case      provides somewhat of an intersection between
the classical and conformal Ricci flows.  Moreover, by rescaling and reparameterizing locally homogenous classical Ricci flows,
one can   
   find    explicit  locally homogeneous   conformal Ricci flows.

The following result   relates   locally homogeneous   classical Ricci flows to  locally homogeneous
   conformal Ricci flows.

\begin{proposition}{\bf (Locally homogeneous   classical  and conformal   Ricci flows)}
\label{yam21} Let
\[  g\,\colon[0,T) \longrightarrow {\cal M}^1 \;,
\qquad t \longmapsto  g(t)
\] be a  locally homogeneous  solution of the classical     Ricci equation  $(\ref{usual})$ 
\begin{eqnarray}  \quad\quad\frac{\partial g}{\partial t} & = &-2 {\rm Ric}(g)+  {{\textstyle\frac{2}{n}}  }
\bar  R (g) g         \;,\qquad g_0\in {\cal M} ^1 \,.      \nonumber 
\end{eqnarray} 
Then $g$ is also a solution of  the historical first Ricci flow equation $(\ref{andusing})$, 
\begin{eqnarray}  \quad\quad\frac{\partial g}{\partial t} & = &-2 {\rm Ric}(g)+  {{\textstyle\frac{2}{n}}  }
  R (g) g         \;,\qquad g_0\in {\cal M} ^1  \,,     \nonumber 
\end{eqnarray} 
and both the volume and  the volume element of the flow are constant, \begin{equation}
\label{LABEL} d \mu _g= d \mu _{g_0}\,.
\end{equation}

 Assume that the spatially constant scalar curvatures satisfy  
$R(g_t)=c_t<0$. (If $M$ is of Yamabe type $-1$, this condition is automatically satisfied, as in Proposition~\ref{increasescale}.) 
Define a new  time parameter
\begin{equation}
\label{LABEL} {s}(t)=\int_0^t|R(g_{t'})|dt'
\end{equation}for $t\in [0,T)$  
and let $ {\bar T}=\lim_{t \rightarrow T} {s}(t)
$ so that $0<{\bar T}\le \infty$. Then $ \frac{d {s}(t)}{dt}=|R(g(t))|>0$ so that  $ {s}(t)$ is
strictly monotonically increasing.
 Let  $t= t( {s})$ denote the  inverse of $ {s}(t)$ so that $t:[0,{\bar T}) \rightarrow [0,T)$.
Define 
a {\bf rescaled
reparameterized  flow}   \begin{equation}
\label{curve1} \label{LABEL} {\bar g} :[0,{\bar T}) \longrightarrow {\cal M}_{-1}  \;,
\qquad {s}\longmapsto  {\bar g} ({s})
   =|R(g(t( {s}))| g(t( {s})) \,.
\end{equation} 
  Then $ {{\bar g}} ( {s}) $ is a locally homogeneous solution to  the conformal Ricci equations,  
 \begin{eqnarray} \frac{\partial {{\bar g}} ( {s})}{\partial {s}}   +2  \left({\rm
Ric}({{\bar g}} ( {s}))+{\textstyle\frac{1}{n}}  g  \right) & = &-   p( {s})  {{\bar g}} (  {s})
\nonumber \\ \nonumber  
\qquad\qquad \quad\qquad R( {{\bar g}} (  {s})) 
&=&-1\,,
\end{eqnarray}  with 
 initial value $ {{\bar g}} _0=|R(g_0)|g_0\in {\cal M}_{-1} $, with
spatially  constant conformal pressure \[
\label{LABEL} p( {s})=
 2\, |{\rm Ric}^T \!({{\bar g}} (  {s}))|^2_{{{\bar g}} (  {s})}  \,,
\] and if $g$ is non-static,  with strictly monotonically decreasing  
  volume 
\begin{equation}
\label{LABEL} \label{strict44}{\rm vol} (M,{{\bar g}} ( {s})) =|R(g(t( {s}))|^{n/2}\,.
\end{equation} 
  \end{proposition}

\noindent  {\bf Proof:} Since the flow   $ t \mapsto g_t$ is locally homogeneous, the scalar curvatures are constant, $R(g_t)= c_t
=$   constant. Thus 
$\bar R(g_t)=\frac{\int_MR(g_t)d
\mu _{g_t}}{{\rm vol}(M,g_t)}=R(g_t)
$ and thus $g_t$ satisfies    (\ref{andusing}),
 \begin{equation}  \frac{\partial g}{\partial t}=-2 {\rm Ric}(g)+  {{\textstyle\frac{2}{n}}  }
   R (g) g \,,    
\end{equation} with $R(g)$ spatially constant.    
Thus  from    (\ref{mustsatisfy}) (with $R(g)$ spatially constant),   
 \begin{eqnarray} \fl\qquad
\frac{ 
\partial R(g)}{\partial t} & = &  {\textstyle\frac{n-2}{n}} 
   \Delta_g R(g) +2
|{\rm Ric}(g)|^2_g 
 - {\textstyle\frac{2}{n}}R^2(g)      
= 2 |{\rm Ric}(g)|^2_g  
   -{\textstyle\frac{2}{n}}  R^2(g) \,,  
\label{backward9}
\end{eqnarray}which  is no longer   a backwards heat equation. Thus in this case   (\ref{andusing})  is equal to
the 
classical Ricci flow equation and is well-posed (see also the remark right after this proof).

 Also,
\begin{eqnarray} \fl \frac{d}{dt} d \mu _{g }  & = & {\textstyle\frac{1}{2}} {\rm tr}_{g }\left(\frac{\partial g
}{\partial t}
\right)d \mu _{g } 
={\textstyle\frac{1}{2}} {\rm tr}_{g }\left(  -2\,{\rm Ric}(g)+  {{\textstyle\frac{2}{n}}  }
  R (g) g            
\right)d \mu _{g } \nonumber     =(
-   R (g) + R (g))d \mu _g=0\,,
\end{eqnarray} 
so that the volume element (and not just the volume) of the flow is preserved.

Now assume that the spatially constant scalar curvatures satisfy  
$R(g_t)=c_t<0$. 
 First we rescale and then we  reparameterize the classical Ricci solution
$ g:[0,T) \rightarrow {\cal M} ^1 
$.
 Since  for each $t\in [0,T)$, $R(g_t)= c_t<0$
is   constant, we can rescale $g_t$  by a   time-dependent  homothetic transformation
 to get  $\tilde g_t=|R(g_t)|g_t=- R(g_t) g_t$. Then 
\begin{equation}
\label{LABEL}R( \tilde g_t)=R(|R(g_t)|g_t)=\frac{1}{|R(g_t)|}R(g_t)=-1\,,
\end{equation} so that $\tilde g_t\in {\cal M}_{-1} $ 
 satisfies the constraint equation of
the conformal Ricci system.

Noting that for 
the rescaled curve  $\tilde g_t=|R(g_t)|g_t=- R(g_t) g_t$,
  ${\rm
Ric}(\tilde g)={\rm Ric}(g)$  and  
$|{\rm Ric}(\tilde g)|^2_{\tilde g} =|R(g)|^{-2}|{\rm Ric}(  g)|^2_{  g}$,
it follows  from   (\ref{backward9}) that $\tilde g$ satisfies  
  the evolution equation\begin{eqnarray}
 \frac{\partial \tilde g}{\partial t}& = &-\frac{\partial R( g)}{\partial t} g -    R (g)\frac{\partial   g}{\partial t}
        \nonumber \\ \nonumber & =&-  \Big(2|{\rm Ric}(g)|^2_g  
   -{\textstyle\frac{2}{n}}  R^2(g) \Big)g-R(g)\Big(-2\,{\rm Ric}(g)+  {{\textstyle\frac{2}{n}}  }  R (g) g         \Big)     
\\
\nonumber     
  & = & -   2|{\rm Ric}(g)|^2_g\,  g + 2 R (g)  \,{\rm Ric}(g)  
\\&=& \nonumber  -   2|R(g)|^2|{\rm Ric}(\tilde g)|^2_{\tilde g}\,  g - 2 |R (g)  |\,{\rm Ric}(\tilde g) 
\\\nonumber 
  & = & -   2|R(g)||{\rm Ric}(\tilde g)|^2_{\tilde g}\, \tilde g - 2 |R (g)  |\,{\rm Ric}(\tilde g) \,. 
\end{eqnarray} 
 Thus $\tilde g_t=|R(g_t)|g_t$ solves the system
\begin{eqnarray}\label{sys1}
\frac{1}{|R(g)|} \frac{\partial \tilde g}{\partial t}&= & - 2  \,{\rm
Ric}(\tilde g)  -   2\, |{\rm Ric}(\tilde g)|^2_{\tilde g}\, \tilde g\\
\qquad  R( \tilde g_t)\; &=&-1\,.\label{sys2}
\end{eqnarray} 
 
Now define a new time parameter $
\label{LABEL} s (t)=\int_0^t|R(g_{t'})|dt'
$ so  that  $ s (0)=0$ and $  \frac{d s }{dt}(t)= |R(g_t)|>0$ so that $ s (t)$ is a strictly monotonically increasing
function of $t$.  Let
$t= t(
s )$ denote the inverse of
$
s (t)$ so that
$
\label{LABEL}  \frac{dt}{d s }(s )= |R(g (t( s ))|^{-1} 
$. 
 Let  \begin{equation}
\label{LABEL}\bar g ( s )=\tilde g (t( s ))=|R(g(t(s))|g(t(s))
\end{equation} 
  denote the reparameterized rescaled flow. 
 Then $\bar g ( s )$ has initial value
$\bar g _0=\bar g (0 )=\tilde g (t( 0))=\tilde g (0)=|R(g (0))|    g (0)=|R(g _0)|    g_0$.

From (\ref{sys2}),  the reparameterized   $\bar g ( s )$
also satisfies the constraint
$R(\bar g (  s  )) =
R(\tilde g(t( s )))
=
 -1
$. Moreover, from (\ref{sys1}), $\bar g ( s )$ satisfies  
  \begin{eqnarray} \frac{\partial\bar g ( s )}{\partial s }& = &\frac{\partial\tilde g }{\partial t
}(t(
s ))
\frac{\partial t }{\partial s  }(s )
  \nonumber  =
 \frac{\partial\tilde g }{\partial t }(t( s )) \frac{1}{|R(g (t( s ))|}
  \nonumber
\nonumber \\&=& - 2  \,{\rm
Ric}(\tilde g(t( s )))  -   2\, |{\rm Ric}(\tilde g(t( s )))|^2_{\tilde g(t( s ))}\, \tilde g(t( s ))
\nonumber \\&=& - 2  \,{\rm
Ric}(\bar g)  -   2\, |{\rm Ric}(\bar g)|^2_{\bar g}\,\bar g (  s  )\,.
\nonumber  
\end{eqnarray} 

Let  $p( s )=
  2(|{\rm Ric}(\bar g(s))|^2_{\bar g(s)} - {\textstyle\frac{1}{n}}) 
=  
 2\, |{\rm Ric}^T \!({{\bar g}} (  {s}))|^2_{{{\bar g}} (  {s})}  $. Then $p(s)$ is spatially constant
since   $g(t)$ and thus $\bar g(s)$ are locally homogeneous. 
 Note that since $p(s)$ is spatially constant,  it is an explicit solution to  (\ref{elliptic6}), 
\[
    L_{\bar g(s)}p(s)=   ( n-1)\Delta_{\bar g(s)}  p(s) +
p (s) =p(s)=
   2(|{\rm Ric}(\bar g(s))|^2_{\bar g(s)} - {\textstyle\frac{1}{n}} )\,.
\]
  Moreover,  
$\bar g ( s )$ satisfies the conformal Ricci system 
 \begin{eqnarray} \frac{\partial\bar g ( s )}{\partial s }& = &   - 2  \,{\rm
Ric}(\bar g(s)) - {\textstyle\frac{2}{n}} \bar g (  s  )- p( s )\bar g (  s  )
\nonumber \\ \nonumber 
 R(\bar g (  s  )) 
&=&-1
\end{eqnarray} with initial value $\bar g _0=|R(g_0)|g_0$.
  Lastly, since ${\rm vol} (M, 
 g(t ))=1$,
\begin{eqnarray} \label{inc2}{\rm vol} (M,\bar g ( s )) & = &{\rm vol}\Big (M,|R(g(t( s ))| g(t( s ))\Big)\nonumber \\&=&
\nonumber |R(g(t(
s ))|^{n/2}{\rm vol}         \Big (M, 
 g(t( s ))\Big)=|R(g(t( s ))|^{n/2}\,. \nonumber  
\end{eqnarray} 
Thus if $g$ is non-static, that
   ${\rm vol}(M,\bar g(s) )$ is strictly monotonically decreasing   follows from Proposition~\ref{volel}.
\hfill$\rule{2mm}{3mm}$ \\
 
We remark that in the locally homogeneous case, 
   $\bar R(g_t)=  R(g_t) $  and thus 
  (\ref{andusing}) 
and   Hamilton's Ricci flow equation are equal. In this case $  \Delta_g R(g)=0$ and thus  the
differential equation for the  scalar curvature reduces to the ordinary differential equation given by (\ref{backward9})  and thus 
is not  a backwards heat equation and is well-posed.

 As a corollary of Proposition~\ref{yam21},   we can deduce the following information about   locally homogeneous
solutions of the classical Ricci flow equation by using  information about     conformal Ricci flow.

\begin{corollary} 
\label{yamabeCor}Let
 \begin{equation}\label{curve1}   g:[0,T) \longrightarrow {\cal M} \;,
\qquad t \longmapsto  g(t)\,,
 \end{equation}with initial value $g(0)=g_0$  be a non-static  locally homogeneous  solution of the classical     Ricci flow equation
\begin{eqnarray}\nonumber   \frac{\partial g}{\partial t} & = &-2\,{\rm Ric}(g)+  {{\textstyle\frac{2}{n}}  }
\bar  R (g) g              \;,\qquad g_0\in {\cal M} ^1 \,, \nonumber 
\end{eqnarray}  with constant negative scalar curvatures $R(g_t)=c_t<0$. (If $M$ is of Yamabe type $-1$, then this condition is
automatically satisfied.) Then\[ \frac{d}{dt}R(g_t)>0\,,\] so that
  the  curve of spatially constant negative scalar  curvatures $R(g_t)$  is strictly monotonically increasing in $t$. 
 \end{corollary}
 {\bf Proof: }Let
$\bar g ( s )=\tilde g (t( s ))=|R(g(t(s))|g(t(s))$ denote the rescaled reparameterized solution of the conformal Ricci system
so that   $\bar g ( s(t) )= |R(g(t )|g(t )$.
From (\ref{strict44}), 
${\rm vol} (M,{{\bar g}} ( {s})) =|R(g(t( {s}))|^{n/2}
$   
is strictly monotonically decreasing
(Proposition~\ref{yam21}). 
From the assumption  
$R(g(t))=c_t<0$,
 \begin{equation}
 R(g(t))=-  {\rm vol}(M,\bar g ( s (t)) )^{2/n} 
\end{equation} 
and is thus strictly monotonically increasing.

If $M$ is of Yamabe type $-1$, then from the definition of Yamabe type (Definition~\ref{defyam}),   the scalar curvature of any  
constant scalar curvature metric must be  negative.
\hfill$\rule{2mm}{3mm}$
\\ 

  Thus a  non-static locally homogeneous  classical Ricci flow $g(t)\in {\cal M} ^1$   
   with $R(g_t)=c_t<0$   has constant volume ($= 1$) but    strictly monotonically  increasing scalar
curvature, $ \frac{d}{dt}R(g_t)>0$,   whereas    the transformed locally homogeneous conformal Ricci solution
$\bar g (
s )\in {\cal M}_{-1} $  has constant scalar curvature ($=-1$) but 
  strictly monotonically decreasing  volume,   $ \frac{d}{dt}{\rm vol}(M,g_t) <0$.

 Since in the locally homogeneous case $R(g_t)$
is spatially constant, 
 under the stronger condition $R(g_t)<0$ (rather than  $R(g_t)\le0$ as occurs  in Proposition~\ref{increasescale}), we get the
stronger conclusion 
$ \frac{d}{dt}R(g_t)>0$ (rather than  $ \frac{d}{dt}R(g_t)\ge0$); see also  the remark following  
Proposition~\ref{increasescale} and also Table~1 in Section~\ref{summary9}. The current situation is summarized in Table~2.
 \begin{table} [htbp] 
\caption[ ]{The volume and scalar curvature for non-static locally homogeneous solutions to the  classical and conformal Ricci
flow equation under either  the assumption that  $R(g_t)<0$ or that $M$ is of Yamabe type $-1$.}
\label{table:isotropy}
{\small 
 \begin{center}
\begin{tabular}{ l l l l l l l}
\hline 
\\[-2.5mm]
     &  &
  \multicolumn{2}{c}{Volume}      & &     \multicolumn{2}{c}{ Scalar curvature } 
\\  [1.5mm]  
\hline  
\\[-5mm]
 \\
Classical Ricci  Flow            && $ {\rm vol}(M,g_t) =1$ &    &   &$\frac{d}{dt}R (g_t)>0$\\ 
            &                &     \\
Conformal Ricci    Flow     &$\qquad$&  $\frac{d}{dt}{\rm vol}(M,g_t)<0  $&         &$\qquad$&      $ R(g_t)=-1$&    
\\ 
          &$\qquad$&        &         & &          &  \\
\hline
\end{tabular}
 \end{center}
}
\end{table}

 \section{The conformal Ricci flow and the $ \sigma $-constant of $M$}
\label{conclusions}
An
important   question  regarding the conformal Ricci flow equation  is under what conditions  does a 
non-equilibrium  initial condition $g_0\in {\cal M}_{-1}
$, ${\rm Ric}(g_0)\not =-{\textstyle\frac{1}{n}} g_0$,
 have an   all-time solution  
 \begin{equation}
\label{LABEL} g:[0,\infty) \longrightarrow  {\cal M}_{-1} \,.
\end{equation} 

 Since the  conformal Ricci flow $g:[0,T) \rightarrow {\cal M}_{-1} $  takes place in $ {\cal M}_{-1} $,   
the curvature norm    $ |{\rm Riem}(g) |_g$ is always  bounded away from zero  (see Proposition~\ref{kinematical2}).
A
 conformal Ricci flow   with  curvature norm   $|{\rm Riem}(g)
|_g$ which   blows up in either  finite  or infinite time   is called {\bf singular}.

  A
conformal Ricci flow $g:[0,\infty) \rightarrow {\cal M}_{-1} $ is {\bf non-singular} if  it is an 
  all-time flow and if     $|{\rm
Riem}(g) |_g$ is   bounded above, i.e., if  there exists a real constant  $ B$  such that for all $t\in [0,\infty)$,
\begin{equation}
\label{LABEL}  {\textstyle\frac{2}{n(n-1)}} \le |{\rm Riem}(g) |_g \le B<\infty\,.
\end{equation}

If $ g:[0,\infty) \rightarrow  {\cal M}_{-1} 
$ is a  non-static  all-time conformal Ricci flow, then  ${\rm vol}(M,g) $
is a strictly monotonically decreasing function of $t$ (see Proposition~\ref{volel}) and from (\ref{vconverges}), 
$ {\rm vol}(M ,g (t))~\rightarrow~v_\infty\ge0$ as $ t \rightarrow \infty $.
  Thus of interest and importance is,
firstly, the value of   
 \begin{equation}
\label{volmap4}  \inf 
 {\rm vol}_{-1} \equiv \inf\limits_{g \in{\cal
M}_{-1} } {\rm vol}(M, g  ) \,, 
\end{equation}and, secondly,  whether $ {\rm vol}(M,g) $ asymptotically approaches this infimum, 
 \begin{equation}
\label{volmap}{\rm vol}(M,g) \longrightarrow  \inf 
 {\rm vol}_{-1}   \qquad\mbox{as}\qquad t \longrightarrow \infty\,.
\end{equation} 

  If $M$ is of Yamabe type $-1$, then the value of $  \inf 
 {\rm vol}_{-1}$ is shown in 
  Fischer-Moncrief \cite{fm00} to be 
\begin{equation}
\label{LABEL}   \inf 
 {\rm vol}_{-1} = (- \sigma (M))^{n/2}\,,
\end{equation} 
where $ \sigma (M)$ denotes the $ \sigma $-constant of $M$, an important topological invariant of $M$ defined by a minimax process
involving the Yamabe functional
 (see Anderson \cite{and97} and Fischer-Moncrief \cite{fm02a}, Section~7).
 For manifolds of Yamabe type $-1$,   a metric $g$  of constant scalar curvature must have $R(g)=$ constant $<0$ and thus  from the
definition of the $ \sigma
$-constant, 
$
\sigma (M)\le0$.

If we assume that    $M$ is of Yamabe type $-1$, that  $ \sigma (M)=0$, 
  that  there exists   a non-singular     conformal Ricci flow 
  $ g\,\colon[0,\infty) \rightarrow  {\cal M}_{-1} 
$,  and that for this flow 
 ${\rm vol}(M,g)
\rightarrow 
\inf 
 {\rm vol}_{-1} =(- \sigma (M))^{n/2} =0 $ as $ t \rightarrow \infty$, then 
  in this case  the conformal Ricci flow volume-collapses $M$
with bounded  curvature norm $|{\rm Riem}(g)|_g$.  Under
these circumstances, $g$ cannot converge to a limit metric
$g_\infty\in {\cal M}_{-1} $ for if it did, by the continuity of the volume functional, $ {\rm vol}(M,g) \rightarrow {\rm
vol}(M,g_\infty) >0$, contradicting ${\rm vol}(M,g) \rightarrow \inf{\rm vol}_{-1} =0$.  Thus the flow must degenerate and the
nature of the degeneration is of importance.

 A closed  orientable 3-manifold $M$ is a {\bf graph manifold}  
  if there is a finite collection $ {\cal T} =\{\mbox{\boldmath$T$}^2 _i\}$ of
disjoint  embedded   tori $\mbox{\boldmath$T$}^2 _i\subset M$ such that each component $M_j$ of
$M\verb+\+\cup\mbox{\boldmath$T$}^2 _i	$ is a Seifert fibered space, i.e., a space that admits a foliation by circles 
 with the property that a foliated tubular neighborhood 
$D^2\times \mbox{\boldmath$S$}^1 $  of each leaf is either the trivial foliation of a solid torus $D^2\times \mbox{\boldmath$S$}^1 $ 
or its quotient by a standard action of a cyclic group 
(see  Anderson \cite{and97}   for more information about graph manifolds).

Thus a graph manifold is a
  union of Seifert fibered
spaces  glued together by toral
automorphisms along  toral boundary components. In  particular, a Seifert fibered manifold
is a graph manifold.

The $ \sigma $-constant of a graph manifold $M$  satisfies $ \sigma (M)\ge0$  (Anderson \cite{and97}).  On the other hand, if
$M$ is of 
 Yamabe type
$-1$, a metric $g$  of constant scalar curvature must have $R(g)=$ constant $<0$ and thus  from the definition of the $ \sigma
$-constant, 
$
\sigma (M)\le0$. Thus 
   a 3-manifold $M$ that is both  a     graph manifold and of
Yamabe type
$-1$ satisfies  $ \sigma (M)=0$.
 Thus if there exists a non-singular conformal Ricci flow on such a manifold  whose volume converges to 
$\inf 
 {\rm vol}_{-1} =(- \sigma (M))^{3/2} =0
$, 
 one expects the 
volume collapse of
$M$  to occur  along either 
$ \mbox{\boldmath$S$} ^1$ or $   \mbox{\boldmath$T$} ^2$-fibers with the embedded   graph manifold 
structure of $M$ describing  how
the degeneration occurs.  
Indeed, for the reduced Einstein field equations,  cosmological models with closed spatial hypersurfaces
 of this type are known
whose conformal volumes collapse to zero 
while the conformal metric collapses 
 with bounded Ricci (or equivalently, Riemann) curvature norm
along either 
$ \mbox{\boldmath$S$} ^1$ or $   \mbox{\boldmath$T$} ^2$-fibers
 (see Fischer-Moncrief \cite{fm00}, \cite{fm02b}, Section~9, and
\cite{fm04})   and it is   expected  that
the conformal Ricci flow behaves similarly    under similar assumptions.

If   $ \sigma (M)<0$, then Anderson {\cite{and97}  has conjectured that either $M$ admits a hyperbolic
metric  or is the union along incompressible tori of finite volume hyperbolic manifolds  and graph manifolds.  In this case, 
  if there exists a non-singular conformal Ricci flow whose volume converges to 
$\inf 
 {\rm vol}_{-1} 
$, 
then  
  ${\rm vol}(M,g)
\rightarrow 
\inf 
 {\rm vol}_{-1} =(- \sigma (M))^{3/2}>0 $ as $ t \rightarrow \infty$ so no volume collapse can occur. Nevertheless, if
   $M$ does not admit a hyperbolic metric, then the conformal Ricci flow cannot converge to a fixed
point,  which must be hyperbolic by Proposition~\ref{equilprop} (see also Remark (iii) following that Proposition) but is
nevertheless seeking to attain the
$
\sigma
$-constant asymptotically 
 insofar as the volume is strictly monotonically decaying to its infimum. 
Thus again the conformal Ricci flow must degenerate and we expect the degeneration to occur   by volume collapsing   
the 
 graph manifold  pieces  but not the hyperbolic pieces. 

If,  on the other hand, 
 $M$ does admit  a hyperbolic metric $g_h\in {\cal M}_{-1} $, then there exists a
  static
 equilibrium solution $g =g_h$ which by Mostow rigidity is  unique up to isometry (Mostow rigidity applies because the scalar
curvature is normalized to $-1$; see remark  about hyperbolic metrics in Section~\ref{crs}). We   conjecture that this static
solution, as it is for    the classical Ricci flow  equation (Ye
\cite{ye93}), the   Einstein equations (Andersson-Moncrief
\cite{am03}), and the reduced Einstein equations    (Fischer-Moncrief
\cite{fm02a})  is {\bf asymptotically stable} under   conformal Ricci flow; i.e., given any neighborhood
 $ {\cal U} _{g_h}\subset
{\cal M}_{-1} $ of $g_h$ there exists a 
{\it capture neighborhood}
 $ {\cal U}' _{g_h}\subseteq {\cal U}  _{g_h} $
 such that if $g_0\in {\cal U}' _{g_h}$ and $g$ is the solution of the conformal Ricci system with initial value $g_0$, 
 then (i) $
  g:[0,\infty) \rightarrow  {\cal M}_{-1} 
$   is non-singular, (ii) for $t\in [0,\infty)$, $g_t\in  {\cal U} _{g_h}$, and (iii) $g_t \rightarrow g_h$ as $ t \rightarrow
\infty$.

\section{The conformal Ricci flow and geometrization of 3-manifolds}
 \label{end0} 
 Using the classical Ricci flow,    Hamilton (\cite{ham95},\cite{ham99}) has 
proposed a program   to prove Thurston's geometrization conjectures for closed 3-manifolds \cite{thu97}. 
The main idea is to try to  show that the Ricci flow evolves  any given initial Riemannian manifold $(M,g_0)$ to a geometric
structure after performing a finite number of surgeries as curvature singularities arise in finite time during the Ricci flow. 
The hope is that if  the surgeries are done before the singularities arise,
the Ricci flow can   be continued for all time 
and
 the initial manifold $(M,g_0)$ will naturally    
decompose into  geometrical pieces.   In this picture the  classical
Ricci flow causes the  initial Riemannian manifold  $(M,g_0)$ to   self-geometrize to a 
limit space $(M_\infty, g_\infty)$  which is   the natural geometrical decomposition of     $(M,g_0)$.

Recently, 
  Perelman (\cite{per02},\cite{per03a},\cite{per03b}) has attempted to further implement Hamilton's program. In  Perelman's work,
much like in the conformal Ricci flow, an important role is played by manifolds of Yamabe type
$-1$, although he doesn't call them that; see \cite{per03a}, p.\ 17, where Perelman
    assumes that the initial manifold does not admit a Riemannian metric with non-negative scalar curvature.
 Since manifolds of Yamabe type 0 and 1 admit metrics with scalar curvature 0 and 1, respectively, these manifold types are
excluded and so Perelman's manifolds    must be of Yamabe type $-1$ (see Definition~\ref{defyam}).

In Perelman's work, 
after each surgery,  the manifold may become disconnected in which case each component is dealt with separately. 
If  
  the initial manifold  develops a component   with nonnegative scalar curvature, then the Ricci flow on that component becomes
extinct so that that  component is removed from further consideration 
  and the Ricci flow then continues on the new manifold   (\cite{per03a}, pp.\ 6, 17).
Note    that the extinct components  which are left behind  are no longer of Yamabe type $-1$,  reflecting again the importance
of Yamabe type
$-1$ manifolds in Perelman's work.

Some of Perelman's techniques   have a natural counterpart for the conformal Ricci flow. For example, for  the unnormalized
Ricci flow, the volume $ {\rm vol}(M,g) $ satisfies $ \frac{d}{dt}{\rm vol}(M,g)=-\int_MR(g) d \mu _g$ and Perelman shows that  
 the rescaled     volume functional $ {\rm vol}(M,g(t)) (t+ {\textstyle\frac{1}{4}}  )^{-3/2}$
is a   non-increasing function  in $t$ (\cite{per03a}, p.\ 17). However, under the conformal Ricci flow, the volume functional
itself  without rescaling  is automatically non-increasing  and  is in fact   strictly monotonically decreasing   
if the flow is non-static (Proposition~\ref{volel}).

As another example of some of these parallel techniques, 
  Perelman  
  identifies an   ``entropy" that increases as the space flows and helps 
move the space toward geometrization. 
For the conformal Ricci flow, the strictly monotonically decreasing   volume $ {\rm vol}(M,g) $ naturally has the same effect
(although the monotonicity is in the opposite direction). For example, globally, the   strictly monotonically decreasing   volume
acts as Perelman's entropy and  prevents periodic orbits in the conformal Ricci flow (Corollary~\ref{quasi2}).   Potential
advantages of this decreasing volume over entropy are that it occurs naturally in the conformal Ricci flow  and is geometrical 
simpler.

Locally, and somewhat  more subtly, the strictly monotonically decreasing {\it local}  volume also
has the potential to      control   the size and shape of collapsing regions under the conformal Ricci flow  since,     if a
singularity is forming within a domain
$D$  of a manifold, the volume must still be decreasing on that domain.  Thus one can reasonably hope to  contain the regions
undergoing    singularity development in the conformal Ricci flow by using the decreasing local volume. 

Thus, altogether,  we believe that the conformal Ricci flow has
the potential to further our understanding of geometrization of 3-manifolds.
  
Lastly, we remark that the
  use of evolution equations to understand and control changes in topology actually has a history longer than is usually
recognized.  Hamilton's  program   has a
precursor in Wheeler's geometrodynamical formulation of general relativity
  (\cite{mw57}, 
reprinted in   \cite{whe62}, and
 \cite{whe68}, pp.\ 279--281, 290--291) which dates back to 1957.  
In  \cite{whe62}, p.\ 306,
   Misner and Wheeler ask 
{\it When the  deterministic evolution  of the metric with time  leads at a certain moment to fission or   coalescence of wormhole
mouths  or to any  other change of topology, what new phenomena occur?}

This question is partially answered in (\cite{whe68}, pp.\ 279--280) where 
 it is recognized that space evolving according to the dynamical formulation  of Einstein's equations can,
  quoting Wheeler, 
 {\it signal its ``intention" to change topology by developing somewhere a curvature that increases without limit
(``gravitational collapse").  To  go further with the analysis of the collapse phenomenon and treat changes in topology forces one
to go outside the framework of classical theory.} A foundation for this new framework, quantum geometrodynamics (QGD), is then
discussed (\cite{whe68}, p.\ 284), with the hope that   one day  QGD might  provide the theoretical basis  for the sequence of
surgeries that an initial Riemannian manifold $(M,g_0)$, together with an initial second fundamental form $k_0\in S_2$,  evolving
under the Einstein  evolution equations,    must undergo  to geometrize itself.

 \ack\addcontentsline{toc}{section}{Acknowledgments}
 The concept of a      Ricci flow of conformal type   was developed jointly by  Professor Vincent  Moncrief    and myself as
an outgrowth of our work on the reduced Einstein equations (\cite{fm94b}--\cite{fm04}) and I would like to acknowledge the many
conversations with him that led to the development of this work.
 I would   like also to thank   the {\it Max-Planck-Institut f\"{u}r
Gravitationsphysik,   Albert-Einstein-Institut}, in Golm, Germany,      {\it The  Erwin Schr\"{o}dinger International Institute
for Mathematical Physics} in Vienna, Austria, and the {\it American Institute of Mathematics} in Palo Alto, California,   for
their gracious  hospitality and    financial support for several periods  during  which  
 this  research was   carried out.

\newpage

\begin{figure}  [htbp] 
\vspace{0.5in} 
\begin{center}
\hspace{0.3in} 
\resizebox{4.5in}{!}{\includegraphics{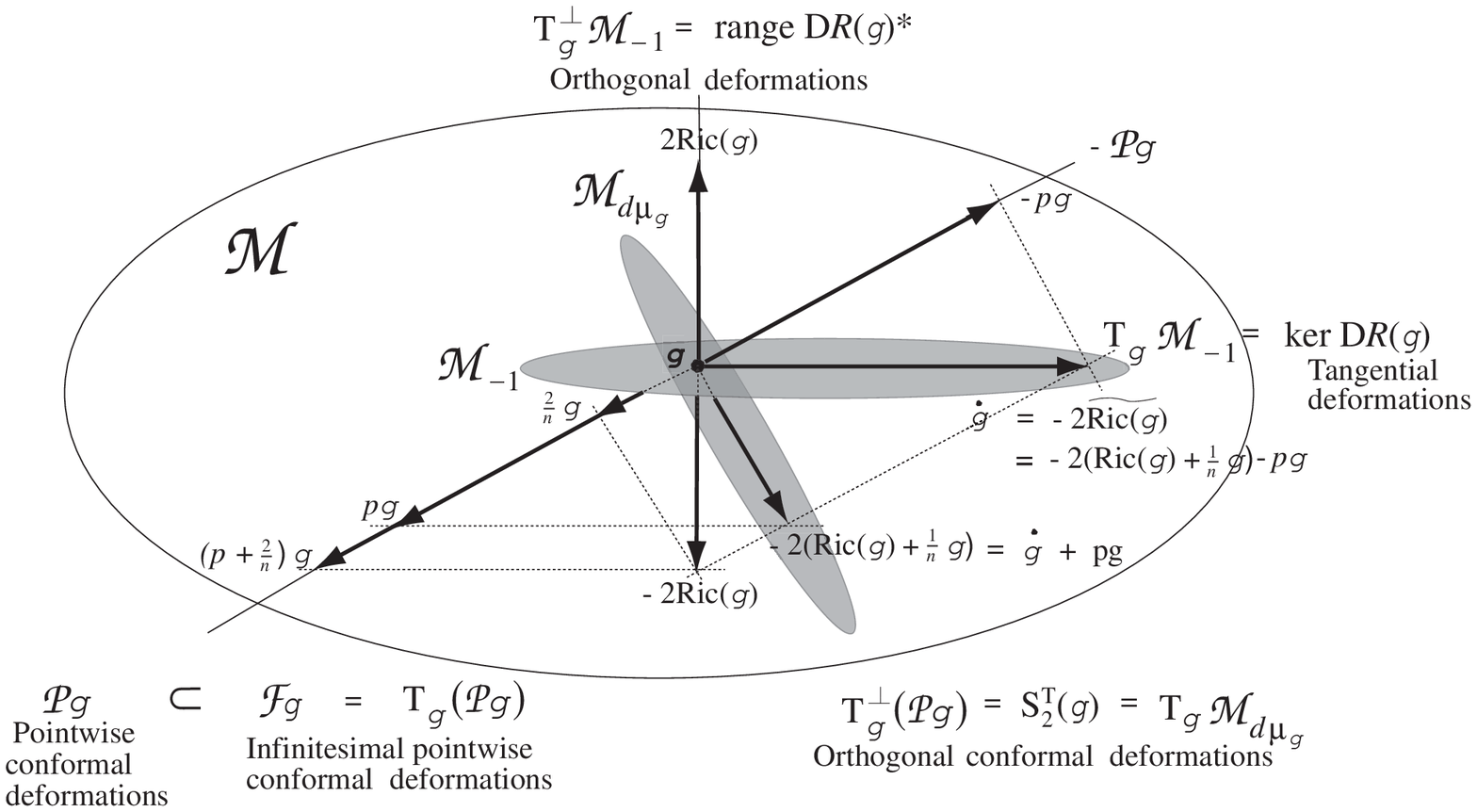}}
\end{center}

\caption {\small 
 A depiction of the  four  ``vector" sums \begin{eqnarray}\;\;\;\quad -  2\,    \widetilde{{\rm Ric}(g) }  &= &-  2{  ({\rm Ric}(g)
+ 
\textstyle\frac{1}{n}} g) -pg   =\frac{\partial g}{\partial t}    
 \nonumber  
 \\  \;\;\;\quad - 2\, {\rm Ric}(g)& =& -  2\,\widetilde {{\rm Ric}(g)}+
(p+ {\textstyle\frac{2}{n}} )g \nonumber
  \\  \;\;\;\quad - 2\, {\rm Ric}(g)& =& -  2 ( {{\rm Ric}(g)}+
  {\textstyle\frac{1}{n}} )g + {\textstyle\frac{2}{n}} g\nonumber \\\nonumber  -2{  ({\rm Ric}(g) +  \textstyle\frac{1}{n}} g)&=& -2
\widetilde {{\rm Ric}(g)}+ pg \,,
\end{eqnarray}

\medskip
\noindent
 where $p=  {   2  L_g ^{-1} (      |{\rm Ric}(g) |^2   )  -  \textstyle\frac{2}{n}} 
   $.      
 The  tangential direction of $ {\cal M}_{-1} $  is     $T_g {\cal M}_{-1} \approx \ker D\!R(g) $ and 
the $L_2$-orthogonal direction is  $ T_g^\perp {\cal M}_{-1} \approx {\rm range}\, D\!R(g)^*$. The infinitesimal conformal
direction is
 $ T_g(Pg)\approx {\cal F} g $ with pointwise orthogonal direction $ T^\perp_g(Pg)\approx S_2^T(g) $.
Note that
the  ``vector"  $ -2  {\rm Ric}(g) \in {\rm range}\, D\!R(g)^*$ 
is orthogonal to $   {\cal M}_{-1} $,  that  $-2{  ({\rm Ric}(g) +  \textstyle\frac{1}{n}} g)$ 
is orthogonal to the conformal orbit $ {\cal P} g$,  and 
 that 
 $  -2  {\rm Ric}(g)   $ and  $  -2{  ({\rm Ric}(g) +  \textstyle\frac{1}{n}} g) $  have the same tangential component
$- 2 \widetilde{ {\rm Ric}(g) }\in \ker D\!R(g)$. 
The first equation above  resolves the inertial ``vector"
$  \frac{\partial g}{\partial t} \in \ker D\!R(g) $ into orthogonal directions 
so that
 the constraint force $-pg$ compensates for the conformal component $pg$ of the nonlinear restoring force
  $-2{  ({\rm Ric}(g) +  \textstyle\frac{1}{n}} g)$ in the  non-orthogonal splitting. The resulting flow then lies in $ {\cal
M}_{-1}
$. 
 
 Lastly,  we note that since $ {\rm Ric}(g) $ is orthogonal to $ {\cal M}_{-1} $ (see Example~\ref{ex2}), the classical unnormalized
Ricci flow (\ref{unnormalRicci}), at a metric $g\in {\cal M}_{-1} $, is orthogonal to
$  {\cal M}_{-1} $, whereas  the conformal Ricci flow  is
tangential to $ {\cal M}_{-1} $. 
Thus in some sense the  classical
unnormalized Ricci flow acts    like a Riemannian  gradient  as it is orthogonal to the level hypersurface $ {\cal M}_{-1} $, 
whereas the conformal Ricci flow acts like a symplectic gradient as it  is tangent  to the level hypersurface $ {\cal M}_{-1} $.
 }  

\label{figure1}
\end{figure}




\section*{References}
 \begin{thebibliography}{9}
 \addcontentsline{toc}{section}{References}

 

\bibitem{and97}
Anderson, M (1997),
{\it Scalar curvature and geometrization conjectures for $3$-manifolds},
in {\it Comparison geometry} (Berkeley, CA 1993-94),
Math.\ Sci.\ Res.\ Inst.\ Publ.\ {\bf30}, Cambridge University Press, Cambridge.

\bibitem{am03}
Andersson, L, and Moncrief, V  (2003),
{\it On the  global evolution problem in $3+1$ dimensional  general relativity},
to appear.

 \bibitem {aub76} 
Aubin, T (1976), 
{\it Equations diff\'erentielles non lin\'eaires et probleme de Yamab\'e concernant la courbure scalaire}, 
J.\ Math.\ Pures Appl.\ {\bf 55}, 269--296.

 \bibitem {aub76b} 
Aubin T (1976), 
{\it The scalar curvature}, in 
{\it Differential geometry and relativity}, Mathematical Phys.\ and Appl.\ Math., Vol. 3, 
Reidel, Dordrecht, 5--18.


\bibitem{be69}
Berger, M, and Ebin, D
(1969),
{\it Some decompositions of the space of symmetric tensors on a Riemannian manifold},
J. Diff. Geom. {\bf 3}, 379-392.

\bibitem {cc99}
Cao, H, and  Chow, B  (1999), 
{\it Recent developments on the Ricci flow},
Bull.\ Am.\ Math.\ Soc.\ {\bf 36}, 59--74.



 
 

 


\bibitem{det83}
DeTurck, D  (1983),
{\it Deforming metrics in the direction of their Ricci tensors},  
J.\ Differential Geometry {\bf 18}, 157--162. 

\bibitem{det03}
DeTurck, D  (2003),
{\it Deforming metrics in the direction of their Ricci tensors (Improved version)},  in
{\it Collected papers on Ricci flow}, H Cao, B Chow, S Chu, and S T Yau, editors,
International Press, Somerville, Massachusetts, 
  163--165. 



 \bibitem{ebin70}
Ebin, D 
(1970),
{\it The space of Riemannian metrics},
Proc. Symp. Pure Math., Amer. Math. Soc. {15},  
11--40.




\bibitem{em70}
Ebin, D, and Marsden, J (1970),
{\it Groups of diffeomorphisms and the motion of an incompressible fluid},
Annals of Mathematics, {\bf 92}, 102-163.


\bibitem{fis04}
Fischer, A (2004),
{\it Solutions to the Ricci flow equation in the $H^s$-setting},
in preparation.


\bibitem {fm72}
Fischer, A, and Marsden, J  (1972), 
{\it The Einstein equations of evolution - a geometric approach},
Journal of Mathematical Physics, {\bf 28}, 1-38.

\bibitem{fm75b}
Fischer, A, and Marsden, J (1975), 
{\it Deformations of the scalar curvature},
Duke Mathematical Journal  {\bf42}, 519--547.  
  
\bibitem {fm94b}
Fischer, A, and Moncrief, V  (1994),
{\it Reducing Einstein's equations to an unconstrained Hamiltonian system on the cotangent bundle of Teichm\"uller space},
in {\it Physics on Manifolds, Proceedings on the International Colloquium in honour of Yvonne Choquet-Bruhat},
edited by M~Flato, R~Kerner, and A~Lichnerowicz,
Kluwer Academic Publishers, Boston,  111--151.


\bibitem {fm96}
Fischer, A, and Moncrief, V  (1996),
{\it The structure of quantum conformal superspace},
in {\it Global Structure and Evolution in General Relativity},
edited by S~Cotsakis and G\ W\ Gibbons,
Springer-Verlag, Berlin, 111--173.


\bibitem {fm97}
Fischer, A, and Moncrief, V  (1997),
{\it Hamiltonian reduction of Einstein's equations of general relativity},
Nuclear Physics B (Proc.~Suppl.) {\bf 57}, 142--161.


\bibitem {fm99a}
Fischer, A, and Moncrief, V  (1999),
{\it The Einstein flow, the $\sigma$-constant and the geometrization of $3$-manifolds},
Classical and Quantum Gravity, {\bf 16}, L79--L87.


\bibitem {fm00}
Fischer, A, and Moncrief, V  (2000),
{\it Hamiltonian reduction, the Einstein flow, and collapse of $3$-manifolds},
Nuclear Physics B (Proc.\ Suppl.) {\bf 88}, 83-102.



\bibitem {fm01}
Fischer, A, and Moncrief, V  (2001),
{\it The reduced Einstein equations and the conformal volume collapse of $3$-manifolds},   Classical and Quantum Gravity, {\bf 18},
4493--4515.



\bibitem {fm02a}
Fischer, A, and Moncrief, V  (2002), {\it Collapse of $3$-manifolds and the reduced Einstein flow}
in {\it Geometry,   Mechanics, and Dynamics: Volume in Honor of the 60th Birthday  
of J.\ E.\  Marsden}, edited by  P   Holmes, P  K  Newton, A   Weinstein,
Springer-Verlag, 463--522.




 


\bibitem {fm02b}
Fischer, A, and Moncrief, V  (2002), {\it Hamiltonian reduction and perturbations of continuously self-similar
$(n+1)$-dimensional Einstein vacuum spacetimes}, Class.\ Quantum Grav.\ {\bf 19}, 5557--5589.


\bibitem {fm04}
Fischer, A, and Moncrief, V  (2004), {\it Hidden structures in closed cosmological models}, to appear.



\bibitem{ft84a}
Fischer, A, and Tromba, A (1984),
{\it On a purely ``Riemannian" proof of the structure
and dimension of the unramified moduli space of a compact Riemann surface},
Mathematische Annalen {\bf 267}, 311--345.


\bibitem{ham82}
Hamilton, R (1982),
{\it Three manifolds with positive Ricci curvature}, Jour.\ Diff.\ Geometry {\bf 17}, 255--306.
 

\bibitem{ham95}
Hamilton, R (1995),
{\it The formation of singularities in the Ricci flow}, Surveys in  Differential Geometry {\bf 2}, 7--136,
International Press, Boston.

 
\bibitem{ham99}
Hamilton, R (1999),
{\it Non-singular solutions of the Ricci flow on three-manifolds}, Commun.\ Anal.\ Geom.\ {\bf 7}, 695--729.

\bibitem{ij92}
Isenberg, J, and Jackson, M  (1992),
{\it Ricci flow of locally homogeneous geometries on closed manifolds},
 J.\ Differential Geometry {\bf 35}, 723-741.
 
\bibitem{lan95}
Lang, S
(1995),
{\it Differentiable and Riemannian manifolds},  
Springer-Verlag, New York.
 

\bibitem{mef72}
Marsden, J, Ebin, D, and Fischer, A (1972), {\it Diffeomorphism groups, hydrodynamics, and relativity}, 
in {\it Proceedings of the Thirteenth Biennial Seminar of the Canadian Mathematical Congress on Differential 
Topology; Differential
Geometry and Applications, Volume~{\bf 1}},   edited by J  R  Vanstone, Canadian Mathematical Society, Montreal, Canada,   pp.\
135--279.

\bibitem{mw57}
Misner, C W,  and  Wheeler, J A  (1957),
{\it Classical physics as geometry},
Ann.\ of Phys.\ {\bf2}, 525--603.
 

\bibitem{omo70}
Omori, H (1970),
{\it On the group of diffeomorphisms of a compact manifold},
Proc. Symp. Pure Math., Amer. Math. Soc. {\bf 15}, 167-183.



\bibitem{pal61}
Palais, R
(1961),
{\it Letter to Serge Lang}, dated November 10, 1961,
unpublished.

\bibitem{pal65}
Palais, R
(1965),
{\it Seminar on the Atiyah-Singer index theorem},
Annals of Mathematics Studies, Number 57,
Princeton University Press, Princeton, New Jersey.


 

\bibitem{pal68}
Palais, R
(1968),
{\it Foundations of global non-linear analysis},
W. A. Benjamin, Inc., New York.

\bibitem{per02}
Perelman, G  (2002),  
{\it The entropy formula for the Ricci flow and its geometric applications},
http://arxiv.org/abs/math/0211159, 1--39.  
 
\bibitem{per03a}
Perelman, G  (2003),  
{\it  Ricci flow with surgery on three-manifolds},\\
http://arxiv.org/abs/math/0303109, 1--22. 

\bibitem{per03b}
Perelman, G  (2003),  
{\it  Finite extinction time for the solutions of the  Ricci flow on certain three-manifolds},\\
http://arxiv.org/abs/math/0307245, 1--7. 


\bibitem{pw84}
Protter, M, and Weinberger, H (1984),
{\it Maximum Principles in Differential Equations}, Springer-Verlag, New York.

\bibitem{sch84}
Schoen, R (1984),
{\it Conformal deformation of a Riemannian metric to constant scalar curvature},
J.\   Differential \ Geom.\ {\bf20},  479--495.
 

  \bibitem{tay96c}
Taylor, M
(1996)
{\it Partial Differential Equations III, Nonlinear Equations},
Springer, New York.

\bibitem{thu97}
Thurston, W (1997),
{\it Three-dimensional geometry and topology}, volume 1, 
edited by S~Levy,
Princeton University Press, Princeton, New Jersey.

\bibitem{whe62}
Wheeler, J A
(1962),
{\it Geometrodynamics},
Academic Press, New York.

\bibitem{whe68}
Wheeler, J A  
(1968),
{\it Superspace and the nature of quantum geometrodynamics},
in {\it Battelle Rencontres -- $1967$ Lectures in Mathematics and Physics}, C.\ DeWitt and
J.\ A.\ Wheeler,   editors, W. A. Benjamin, Inc., New York. 

\bibitem{ye93}
  Ye, R (1993),  
{\it Ricci flow, Einstein metrics and space forms},
Trans.\ Amer.\ Math.\ Soc.\ {\bf338}, 871--896. 
 \end{thebibliography}
 \end{document}